%% file: Quasitilings_and_Applications.tex
\pdfoutput=1
\documentclass[11pt,a4paper,twoside]{book}
\usepackage{TFG}
\usepackage{enumitem}
\usepackage[colorlinks=true,linkcolor=blue]{hyperref}
\hypersetup{
	pdfinfo={
		Title={Quasitilings and Applications},
		Author={Lander Guerrero Sanchez},
		Director1={andrei.jaikin},
		Director2={},
		Ndirectores={1},
		Tipo={TFM},
		Curso={2018-2019},
		MSC={43A07},
		Palabrasclave={amenable, groups, quasitilings, sofic},
	}
}
\usepackage{tikz}
\usetikzlibrary{arrows,decorations.pathmorphing,automata}
\usepackage{tikz-cd}
\usepackage{color}

\usepackage{pdfpages}

\DeclarePairedDelimiter{\abs}{\lvert}{\rvert}
\DeclarePairedDelimiter{\norm}{\lVert}{\rVert}
\DeclarePairedDelimiter{\scalar}{\langle}{\rangle}
\DeclarePairedDelimiter{\ceil}{\lceil}{\rceil}

\newcommand{\Oo}{\mathcal{O}}

\newcommand{\m}{\mathfrak{m}}

\newcommand{\ninN}{_{n\in\N}}

\newcommand{\isom}{\cong}

\DeclareMathOperator{\Mat}{Mat}
\DeclareMathOperator{\GL}{GL}

\newcommand{\PM}{\mathcal{PM}}
\newcommand{\M}{\mathcal{M}}

\newcommand{\I}{\mathcal{I}(\Oo)}
\newcommand{\OO}{\mathcal{O}}
\newcommand{\Imax}{\mathcal{I}(\Oo)_{\operatorname{max}}}

\newcommand{\SSS}{\mathbb{S}}
\newcommand{\B}{\mathbb{B}}

\newcommand{\detp}{{\det}_+}

\DeclareMathOperator{\SO}{SO}
\DeclareMathOperator{\E}{E}

\DeclareMathOperator{\Cay}{Cay}

\DeclareMathOperator{\id}{id}
\DeclareMathOperator{\tr}{tr}

\DeclareMathOperator{\rk}{rk}
\DeclareMathOperator{\im}{im}
\DeclareMathOperator{\proj}{proj}
\newcommand{\dd}{\mathrm{d}}
\DeclareMathOperator{\normaleq}{\unlhd}

\DeclareMathOperator{\spec}{spec}

\newcommand{\tors}{_{\operatorname{tors}}}

\makeatletter
\def\moverlay{\mathpalette\mov@rlay}
\def\mov@rlay#1#2{\leavevmode\vtop{%
		\baselineskip\z@skip \lineskiplimit-\maxdimen
		\ialign{\hfil$\m@th#1##$\hfil\cr#2\crcr}}}
\newcommand{\charfusion}[3][\mathord]{
	#1{\ifx#1\mathop\vphantom{#2}\fi
		\mathpalette\mov@rlay{#2\cr#3}
	}
	\ifx#1\mathop\expandafter\displaylimits\fi}
\makeatother

\newcommand{\cupdot}{\charfusion[\mathbin]{\cup}{\cdot}}
\newcommand{\bigcupdot}{\charfusion[\mathop]{\bigcup}{\cdot}}

\begin{document}
\frontmatter

%
%
\begin{titlepage}
	\includepdf[pages=-]{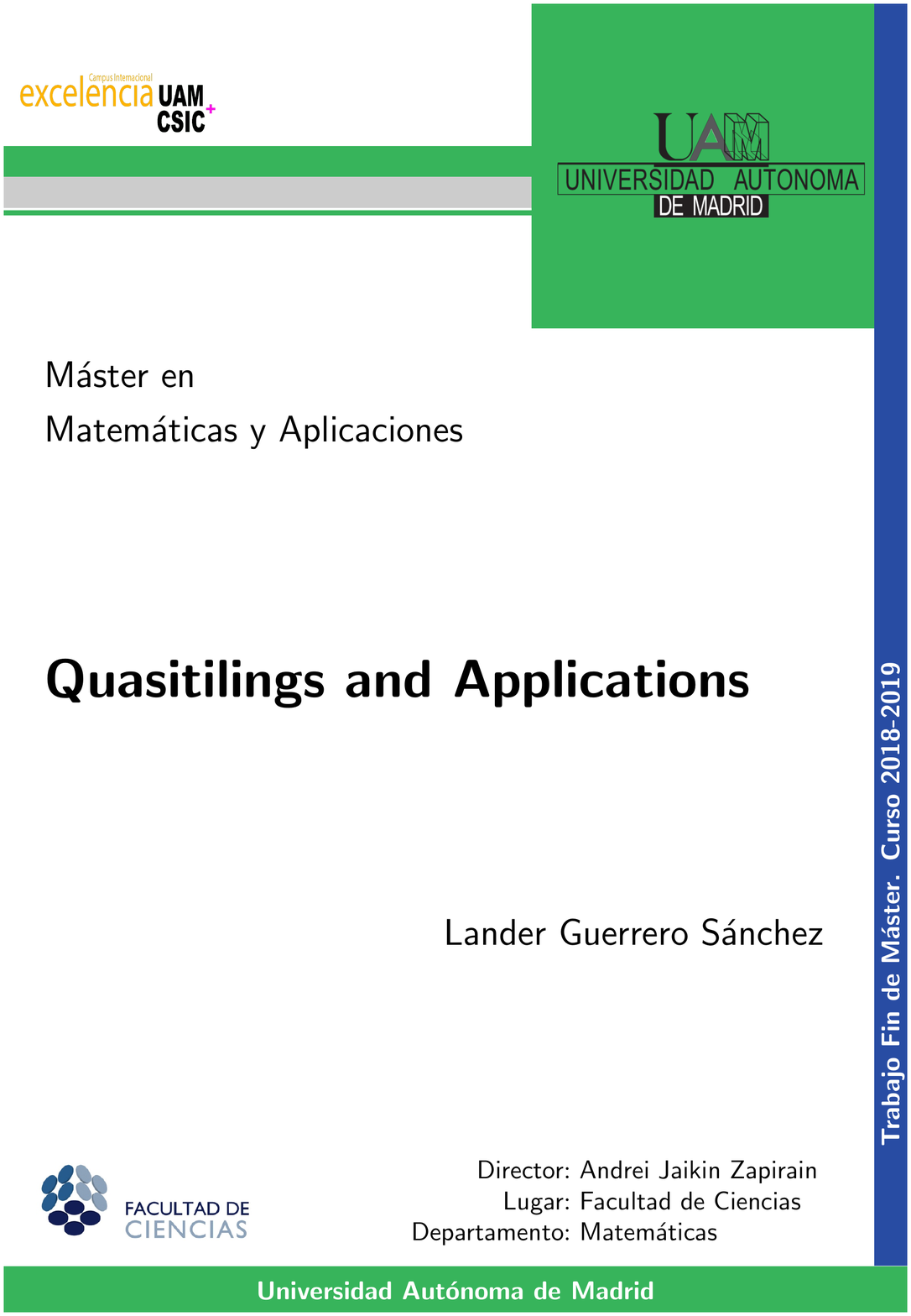}
\end{titlepage}
	
\thispagestyle{empty}
\pagestyle{plain}
\tableofcontents
\clearpage{\pagestyle{empty}\cleardoublepage}
\pagestyle{fancy}
\fancyhf{}
\renewcommand{\chaptermark}[1]{\markboth{#1}{}}
\fancyhead[LO]{\slshape\nouppercase{\leftmark}}
\fancyhead[RO]{\thepage}
\fancyhead[LE]{\thepage}
\fancyhead[RE]{\slshape\nouppercase{\rightmark}}
	
%
%
	
\include{Introduction}\clearpage{\pagestyle{empty}\cleardoublepage}

%
%
	
\mainmatter
\renewcommand{\chaptermark}[1]{\markboth{\chaptername\ \thechapter. #1}{}}
	
%
%

\include{Chapter1}\clearpage{\pagestyle{empty}\cleardoublepage}
\include{Chapter2}\clearpage{\pagestyle{empty}\cleardoublepage}
\include{Chapter3}\clearpage{\pagestyle{empty}\cleardoublepage}
\include{Chapter4}\clearpage{\pagestyle{empty}\cleardoublepage}
\include{Chapter5}\clearpage{\pagestyle{empty}\cleardoublepage}

%
%

\include{Bibliography}\clearpage{\pagestyle{empty}\cleardoublepage}

\end{document}

%% file: Introduction.tex
\chapter{Introduction}

This thesis aims to serve as an introduction to the theory of quasitilings for amenable groups. In order to showcase the power of this theory, we focus on the study of the Sofic L\"uck Approximation Conjecture, which can be proven for amenable groups by making use of quasitilings. The first four chapters of the thesis are an exposition of the aforementioned topics, collected from the literature. After that, we present some new results in the fifth and final chapter.


Amenable groups originated in 1929 from J. von Neumann's work on the Banach-Tarski Paradox in \cite{von-neumann}. This so-called paradox, proved in 1924 by S. Banach and A. Tarski \cite{banach-tarski}, states that a ball in the euclidean three-dimensional space can be decomposed into a finite number of pieces that can then be rearranged to form two new balls of the same size as the original ball, using only translations and rotations. The key to this result lies on the fact that the group of isometries of $\R^3$ contains a copy of the free group of rank two. This led to  von Neumann introducing amenable groups as those with a finitely-additive probability measure that is invariant under the action of the group on itself. These are precisely the groups that cannot cause a paradoxical decomposition akin to the one in the Banach-Tarski Paradox. It was then conjectured that a group is amenable if and only it contains a free subgroup of rank two. This came to be known as the von Neumann Conjecture, and was  disproved in 1980 by A. Y. Ol'shanskii \cite{olshanskii}.

The term amenable was later coined by M. M. Day \cite{day} as a pun on the word mean, after he showed that 
amenable groups are those on which an invariant mean can be defined. Another equivalent definition was found by E. F\o lner \cite{folner}, characterising amenable groups as those with almost-invariant finite subsets. Subsequently, amenable groups have been extensively studied, and a plethora of different characterisations of amenability has been found, making amenable groups ubiquitous across many seemingly distant areas of mathematics.

The theory of quasitilings for amenable groups was first developed by D. S. Ornstein and B. Weiss  \cite{ornstein-weiss}, when they proved that any sufficiently invariant finite subset of an amenable group can be covered almost entirely by almost-disjoint translates of a finite collection of tiles with good invariance properties. The existence of these  quasitilings, obtained by using F\o lner sets, has far-reaching applications in the study of many problems concerning amenable groups. A more general version of this theory, valid not only for finite subsets of the group but also for finite labelled graphs, was introduced by G. Elek \cite{elek}. 

Amenable groups are in a way groups of a finite-like nature, in the sense that they can be approximated by finite F\o lner sets. Residually finite groups, in which elements can be distinguished in finite quotients, are of a similar nature in that they can be approximated by finite groups. As a joint generalisation of both amenable and residually finite groups arise sofic groups, first introduced by M. Gromov \cite{gromov} in 1999 as groups whose Cayley graphs can be approximated by finite graphs. Soon after in 2000, B. Weiss \cite{weiss} gave these groups the name sofic, a term that comes from the Hebrew word for finite. Both amenable and residually finite groups are sofic, and there are currently no known examples of non-sofic groups. In \cite{elek-szabo}, using the theory of quasitilings applied to the sofic approximations of an amenable group, G. Elek and E. Szab\'o were able to characterise amenable groups amongst sofic groups as those whose sofic approximations are all conjugate. 

In this same spirit of using finite approximations to obtain information about infinite objects, we have the Sofic L\"uck Approximation Conjecture, a version of a conjecture that has its origin in a work of W. L\"uck on approximations of $L^2$-Betti numbers of compact manifolds. Given an element of the group algebra of a sofic group over some field, we can naturally define an operator for each element in the sofic approximation of our group. The Sofic L\"uck Approximation Conjecture then asks whether the normalised dimensions of the kernels of these associated operators converge, and whether this convergence is independent of the chosen sofic approximation.

In the case that we are working in a field of characteristic zero, this conjecture has been extensively studied, and was eventually shown to be true for any sofic group by A. Jaikin-Zapirain \cite{jaikin2}. The proof of this fact relies heavily on techniques from functional analysis, in particular the spectral theory of self-adjoint operators, which cannot be readily exported to the case of positive characteristic. As such, the positive characteristic case of the conjecture remains open.

Nonetheless, the conjecture can be shown to hold for amenable groups, independent of the characteristic of the field, by making use of the previously mentioned result by Elek and Szab\'o from \cite{elek-szabo} that says that any two sofic approximations of an amenable group are conjugate.

The proof of the conjecture in characteristic zero relies on the construction of a sequence of measures, each associated to an element of the sofic approximation of the group. Proving the conjecture is then reduced to the problem of showing that these measures converge pointwise at zero, independent of the approximation.

Suppose now that we are working over the field of fractions of some discrete valuation ring, e.g. the ring of $p$-adic integers $\Z_p$ with its field of fractions $\Q_p$. Using the Smith normal form of a matrix over a principal ideal domain, we can define a measure on the space of ideals of our discrete valuation ring for each element of the sofic approximation. For amenable groups, these measures can be shown to converge at each ideal, independent of the sofic approximation.

This construction can be generalised to the case of number fields, whose rings of integers are Dedekind domains. This time, the construction of the associated measures on the space of ideals is done not by using the Smith normal form, but the decomposition of finitely generated modules over Dedekind domain. In this case, we are able to prove the strong convergence of the measures for amenable groups, independent of the sofic approximation.

Chapter~\ref{chapter:AmenableGroups} serves as a standard introduction to the theory of amenable groups and their basic properties. Throughout the chapter, a number of the many different characterisations of amenability are discussed, before eventually proving in the last section the equivalences between them. We also make some room in the middle of the chapter for the proof of the Banach-Tarski Paradox.




In Chapter~\ref{chapter:Quasitilings}, we develop the theory of quasitilings of graphs for amenable groups, starting with by defining what an approximation of a Cayley graph before going on to prove that quasitilings always exist for amenable groups. We then present a different version of this result using linear combinations of graphs.

Chapter~\ref{chapter:SoficGroups} begins with a brief discussion of residually finite groups, after which we introduce the concept of sofic groups via sofic approximations. We also discuss a characterisation of sofic groups making use of ultraproducts of finite symmetric groups.

In Chapter~\ref{chapter:LuckConjecture}, we discuss the Sofic L\"uck Approximation Conjecture and we prove it for amenable groups and over the field $\Q$.

In Chapter~\ref{chapter:LastChapter}, we construct a sequence of measures associated to the operators that appear in the Sofic L\"uck Approximation Conjecture, first over discrete valuation rings, and then over number fields, and study the convergence of these measures for amenable groups.

%% file: Chapter1.tex
\chapter{Amenable Groups}
\label{chapter:AmenableGroups}

In this chapter we will study amenable groups and some of their many characterisations. Our introduction to the concept of amenability will be through the original definition in terms of invariant finitely additive probability measures. We will also discuss its relation to means and the closure properties of the class of amenable groups. Afterwards, we will study the characterisation of amenable groups in terms of the F\o lner condition, and the related concepts of F\o lner nets and sequences. We will then discuss paradoxical decompositions, after which we will make a slight digression to prove the Banach-Tarski Paradox. Then, we will introduce the concepts of ultrafilters and the Stone-\u{C}ech compactification, in order to later prove some characterisations of amenability in terms of fixed points and measures of certain actions. We will conclude  the chapter by finally proving the equivalences between all of the characterisations of amenability that we have discussed. This chapter is mainly based on \cite[\S 4]{ceccherini}, \cite[\S 4.1]{kerr} and \cite{garrido}.

\section{Finitely Additive Measures and Means}
\label{section:AmenableGroups}

In 1924, S. Banach and A. Tarski \cite{banach-tarski} proved that the unit ball in $\R^3$ can be partitioned into five pieces which can then be used to form two disjoint copies of the original ball using only translations and rotations. This result, known as the Banach-Tarski Paradox, can be reformulated as saying that there is no finitely additive measure on $\R^3$ that is invariant under translations and rotations. Amenable groups were originally defined in 1929 by J. von Neumann \cite{von-neumann} whilst studying the Banach-Tarski Paradox. We will now present his original definition in terms of invariant finitely additive probability measures.


\begin{definition}\label{def:ProbMeasure}
	A \emph{finitely additive probability measure} on a group $G$ is a map $$\apl{\mu}{\mathcal{P}(G)}{[0,1]}$$ satisfying the following properties:
	\begin{enumerate}
		\item $\mu(G)=1$.
		
		\item $\mu(X\cupdot Y)=\mu(X)+\mu(Y)$ for all disjoint $X,Y\subseteq G$, where $\cupdot$ denotes disjoint union.
	\end{enumerate}
	Furthermore, we say that $\mu$ is \emph{right-invariant} if it satisfies the following additional property:
	\begin{enumerate}
		\item[(iii)] $\mu(Xg) = \mu(X)$ for all $X\subseteq G$ and $g\in G$. 
	\end{enumerate}
\end{definition}

We can now give our first definition of amenable groups.

\begin{definition}\label{def:AmenableGroup}
	A group $G$ is said to be \emph{amenable} if there exists a right-invariant finitely additive probability measure on $G$.
\end{definition}

\begin{examples}\label{ex:First amenable groups}
	\begin{enumerate}
		\item Let $G$ be a finite group. Then, we can define a right-invariant finitely additive probability measure $\mu$ on $G$ by setting
		\begin{displaymath}
		\mu(X)=\frac{\abs{X}}{\abs{G}}
		\end{displaymath}
		for any $X\subseteq G$, and so $G$ is amenable.


		\item Let $F(a,b)$ be the free group on two generators. Then, $F(a,b)$ is not amenable. Indeed, for each $s\in\{a^{\pm 1}, b^{\pm 1}\}$ denote by $W(s)$ the set of reduced words ending with $s$. Then, we can write 
		\begin{align*}
		F(a,b) & = \{1\}\cupdot W(a)\cupdot W(a^{-1}) \cupdot W(b) \cupdot W(b^{-1}) \\
		& = W(a) \cupdot W(a^{-1})a \\
		& = W(b) \cupdot W(b^{-1})b.
		\end{align*}
		Assume by contradiction that there is a right-invariant finitely additive probability measure $\mu$ on $F(a,b)$. Then, on the one hand we have that 
		\begin{align*}
		1 & = \mu\big(F(a,b)\big) \\ 
		& = \mu\big(W(a)\big) + \mu\big(W(a^{-1})\big) \\
		& = \mu\big(W(b)\big) + \mu\big(W(b^{-1})\big).
		\end{align*}
		On the other hand,
		\begin{align*}
		\mu\big(F(a,b)\big) & \geq  
		\mu\big(W(a)\big) + \mu\big(W(a^{-1})\big) + \mu\big(W(b)\big) + \mu\big(W(b^{-1})\big) \\
		& =2,
		\end{align*}
		so we have a contradiction. Therefore, $F(a,b)$ is not amenable.
	\end{enumerate}
\end{examples}

\begin{remarks}
	\begin{enumerate}
		\item It is not difficult to see that the existence of a right-invariant finitely additive probability measure on $G$ is equivalent to the existence of a left-invariant finitely additive probability measure on $G$, i.e. a finitely additive probability measure $\mu$ on $G$ such that $\mu(gX)=\mu(X)$ for all $X\subseteq G$ and $g\in G$.

		\item We are dealing exclusively with discrete groups. Nonetheless, a more general theory of amenability can be developed for locally compact groups. 
	\end{enumerate}
\end{remarks}

If we denote by $\PM(G)$ the set of finitely additive probability measures on the group $G$, then we can define a right action of $G$ on $\PM(G)$ by setting
\begin{displaymath}
\mu^g(X)=\mu(Xg^{-1})
\end{displaymath}
for $\mu\in \PM(G)$ and $X\subseteq G$. Observe that $\PM(G)^G$, the subset of $G$-invariant elements in $\PM(G)$, is precisely the set of right-invariant finitely additive probability measures on $G$. Thus, $G$ is amenable if and only if $\PM(G)^G\not=\emptyset$.


The definition of amenability that we have given suffers from the fact that, in general, finitely additive measures are not $\sigma$-additive and, as a consequence, we cannot make full use of the theory of Lebesgue integration. For this reason, M. Day gave in \cite{day} a new characterisation of amenability that allows us to use techniques from Functional Analysis.

Recall that for a set $\Omega$ the space $\ell^{\infty}_{\R}(\Omega)$ of bounded functions $x\colon \Omega\longrightarrow \R$ is a Banach space with the supremum norm
\begin{displaymath}
\norm{x}_{\infty}=\sup_{w\in \Omega}\abs{x(w)}.
\end{displaymath}
Given $\lambda\in\R$, we denote by $\lambda$ the constant map in $\ell^{\infty}_{\R}(E)$ taking the constant value $\lambda$ on all $\Omega$. We can order $\ell^{\infty}_{\R}(E)$ by setting $x\leq y$ if and only if $x(w)\leq y(w)$ for all $w\in \Omega$.

If $G$ is a group, then we can consider the action of $G$ on $\ell^{\infty}_{\R}(G)$ given by
\begin{displaymath}
x^g(h)=x(hg^{-1})
\end{displaymath}
for $x\in\ell^{\infty}_{\R}(G)$ and $h\in G$.

Let us now introduce the concept of a mean.

\begin{definition}\label{def:Mean}
	A \emph{mean} on a group $G$ is a linear map $$\apl{m}{\ell^{\infty}_{\R}(G)}{\R}$$ satisfying the following properties:
	\begin{enumerate}
		\item $m(1)=1$.
		 
		\item $m(x)\geq 0$ for all $x\in \ell^{\infty}_{\R}(G)$ such that $x\geq 0$. 
	\end{enumerate}
	Furthermore, we say that $m$ is \emph{right-invariant} if it satisfies the following additional property:
	\begin{enumerate}
		\item[(iii)] $m(x^g) = m(x)$ for all $x\in\ell^{\infty}_{\R}(G)$ and $g\in G$.
	\end{enumerate}
\end{definition}

Denoting by $\M(G)$ the set of means on the group $G$, we have that the action of $G$ on $\ell^{\infty}_{\R}(G)$ restricts to an action of $G$ on $\M(G)$. Then, $\M(G)^G$, the subset of $G$-invariant elements in $\M(G)$, is precisely the set of right-invariant means on $G$. 

Let $m$ be a mean on $G$. Given any $X\subseteq G$, we can consider its characteristic function $\chi_X\in\ell^{\infty}_{\R}(G)$. We can then define a finitely additive probability measure $\mu_m$ on $G$ by setting
\begin{displaymath}
\mu_m(X)=m(\chi_X)
\end{displaymath}
and, if $m$ is right-invariant, then so is $\mu$.

Conversely, given a right-invariant finitely additive probability measure $\mu$ on $G$, we can construct an integral in a manner analogous to the construction of the Lebesgue integral of a $\sigma$-additive measure, although some of the properties of the Lebesgue integral fail in our case due to the lack of $\sigma$-additivity. First, we consider $\R[G]$ the space of functions in $\ell^{\infty}_{\R}(G)$ with finite support. Then, given $x\in \R[G]$ we can define 
\begin{displaymath}
\int_G x\ \dd\mu = \sum_{\lambda\in x(G)} \lambda \mu\big(x^{-1}(\lambda)\big).
\end{displaymath}
It is easy to see that this integral satisfies the following properties:
\begin{enumerate}
	\item It is linear, i.e. 
	\begin{displaymath}
	\int_G (\alpha x+\beta y)\ \dd \mu =\alpha\int_G x\ \dd \mu + \beta \int_G y\ \dd \mu
	\end{displaymath}
	for any $x,y\in\R[G]$ and $\alpha,\beta\in\R$. 
	
	\item It is bounded, i.e.
	\begin{displaymath}
	\abs[\bigg]{\int_G x\ \dd\mu}\leq \norm{x}_{\infty}
	\end{displaymath}
	for any $x\in\R[G]$.
	
	\item It is right-invariant, i.e. 
	\begin{displaymath}
	\int_G x^g\ \dd\mu=\int_G x\ \dd\mu
	\end{displaymath}
	for any $x\in\R[G]$ and $g\in G$.
\end{enumerate} 
As a consequence, this integral defines a right-invariant bounded functional on $\R[G]$. Now, $\R[G]$ forms a dense subspace of $\ell^{\infty}_{\R}(G)$, and so the integral can be extended to a right-invariant bounded functional on the whole $\ell_{\R}^{\infty}(G)$. 
%
Therefore, the map $m_{\mu}\colon \ell^{\infty}_{\R}(G)\longrightarrow \R$ defined by
\begin{displaymath}
m_{\mu}(x)=\int_G x\ \dd \mu
\end{displaymath}
for $x\in\ell^{\infty}_{\R}(G)$ is a right-invariant mean on $G$. 
A more detailed version of this construction can be found in \cite{ceccherini}. 

The previous discussion shows that there is a bijection between the sets $\M(G)^G$ and $\PM(G)^G$, which leads us to the following characterisation of amenability.

\begin{theorem}\label{thm:ProbabilityMeasure<=>Mean}
	Let $G$ be a group. Then, $G$ is amenable if and only if there exists a right-invariant mean on $G$.
\end{theorem}

	

\section{Closure Properties}

We will now study some closure properties of the class of amenable groups. Specifically, we will show that amenability is closed under taking subgroups, extensions and direct limits. But first, let us quickly recall the notion of direct limit of groups.

\begin{definition}
	Let $I$ be a directed set, i.e. a partially ordered set such that for any $i,j\in I$ there is some $k\in I$ with $i,j\leq k$. A \emph{direct system of groups} consists of a collection of groups $(G_i)_{i\in I}$ and homomorphisms $\phi_{ij}\colon A_i\longrightarrow A_j$ for all $i\leq j$ such that the following hold:
	\begin{enumerate}
		\item $\phi_{ii}=\id_{G_i}$ for all $i\in I$.
		
		\item $\phi_{ik}=\phi_{jk}\circ\phi_{ij}$ for all $i,j,k\in I$ with $i\leq j\leq k$.
	\end{enumerate}
	The \emph{direct limit} of the direct system  $(G_i)_{i\in I}$ is then defined as the group
	\begin{displaymath}
	\varinjlim_{i\in I}G_i=\bigg(\bigcup_{i\in I}G_i\bigg)/\sim,
	\end{displaymath}
	with the equivalence relation $\sim$ given by setting $g_i\sim g_j$ for $g_i\in G_i$ and $g_j\in G_j$ if and only if there exists some $k\in I$ with $i,j\leq k$ such that $\phi_{ik}(g_i)=\phi_{jk}(g_j)$.
\end{definition}

\begin{example}\label{ex:Direct Limit fg subgroups}
	Given any group $G$, we can order the family $I$ of finitely generated subgroups of $G$ by inclusion, which is thus turned into a direct system of groups. We can then easily see that
	\begin{displaymath}
	G=\varinjlim_{H\in I}H.
	\end{displaymath}
	Consequently, every group can be written as the direct limit of its finitely generated subgroups.
\end{example}

\begin{proposition}\label{Prop:Amenable closed Sbgp Ext DirLim}
	Let $G$ be a group. Then, the following properties hold:
	\begin{enumerate}
		\item If $G$ is amenable and $H\leq G$, then $H$ is amenable.
		
		\item If $N\unlhd G$, then $G$ is amenable if and only if both $N$ and $G/N$ are amenable.
		
		\item 
		If $(G_i)_{i\in I}$ is a direct system of amenable groups and $$G=\varinjlim_{i\in I}G_i,$$ then $G$ is amenable.
	\end{enumerate}
\end{proposition}

\begin{proof}
	\begin{enumerate}
		\item Let $\mu$ be a right-invariant finitely additive probability measure on $G$, and $T$ be a left transversal of $H$ in $G$, i.e. a set of representatives of the left cosets of $H$ in $G$. Then, we define
		\begin{displaymath}
		\tilde\mu (X)=\mu(TX)
		\end{displaymath}
		for any $X\subseteq H$. We can easily check that $\tilde{\mu}\colon \mathcal{P}(H)\longrightarrow \R$ is a right-invariant finitely additive probability measure on $H$. Indeed, we have that
		\begin{displaymath}
		\tilde{\mu}(H)=\mu(TH)=\mu(G)=1.
		\end{displaymath}
		Furthermore, if $X,Y\subseteq H$ are disjoint, then so are $TX$ and $TY$, and as a consequence
		\begin{align*}
		\tilde{\mu}(X\cupdot Y) & = \mu\big(T(X\cupdot Y)\big) \\
		& = \mu(TX \cupdot TY) \\
		& = \mu(TX)+\mu(TY) \\
		& = \tilde{\mu}(X)+\tilde{\mu}(Y).
		\end{align*}
		Finally, given any $X\subseteq H$ and $h\in H$ we have that
		\begin{displaymath}
		\tilde{\mu}(Xh)=\mu(TXh)=\mu(TX)=\tilde{\mu}(X).
		\end{displaymath} 
		Therefore, $H$ is amenable.
		
		\item Assume first that $G$ is amenable. Then, item (i) implies that $N$ is also amenable. 
		Now, let $\mu$ be a right-invariant finitely additive probability measure on $G$. Then, we define
		\begin{displaymath}
		\tilde\mu (X/N)=\mu(X)
		\end{displaymath}
		for any $X/N\subseteq G/N$. We have that
		\begin{displaymath}
		\tilde{\mu}(G/N)=\mu(G)=1.
		\end{displaymath}
		Furthermore, if $X/N,Y/N\subseteq G/N$ are disjoint, then so are $X$ and $Y$, and as a consequence
		\begin{align*}
			\tilde{\mu}(X/N\cupdot Y/N) & = \tilde\mu\big((X\cupdot Y)/N\big) \\
			& = \mu(X \cupdot Y) \\
			& = \mu(X)+\mu(Y) \\
			& = \tilde{\mu}(X/N)+\tilde{\mu}(Y/N).
		\end{align*}
		Finally, given any $X/N\subseteq G/N$ and $gN\in G/N$ we have that
		\begin{align*}
		\tilde{\mu}\big((X/N)(gN)\big) & =\tilde{\mu}\big((Xg)/N\big) \\
		& =\mu(Xg) \\
		& =\mu(X) \\
		& =\tilde{\mu}(X/N).
		\end{align*} 
		Therefore, $G/N$ is amenable.

		Conversely, assume that both $N$ and $G/N$ are amenable. Let $\mu_{N}$ and $\mu_{G/N}$ be right-invariant, finitely additive probability measures on $N$ and $G/N$, respectively. Then, for any $gN\in G/N$ the map $\mu_{N}^g$ defines a finitely additive probability measure on $gN$. Note that this measure does not depend on the representative of $gN$ chosen, for if $gN=hN$, then $$\mu_{N}^g=\mu_{N}^{gh^{-1}h}=\mu^h$$ because $gh^{-1}\in N$ and $\mu_{N}$ is $N$-invariant. 
		
		Now, given $X\subseteq G$ we set $$\gamma_X(gN)=\mu_{N}^g(X\cap gN)$$ for $gN\in G/N$. Then, it is clear that $\gamma_X\in\ell_{\R}^{\infty}(G/N)$. Furthermore, 
		if $X,Y\subseteq G$ are disjoint, then so are $X\cap gN$ and $Y\cap gN$ for all $gN\in G/N$, and hence,
		\begin{align*}
		\gamma_{X\cupdot Y}(gN) & = \mu_{N}^g\big((X\cupdot Y)\cap gN\big) \\
		& = \mu_{N}^g(X\cap gN)+\mu_{N}^g(Y\cap gN) \\
		& = \gamma_X(gN)+\gamma_Y(gN)
		\end{align*}
		for any $gN\in G/N$. 
		Moreover, given $g\in G$, we have that 
		\begin{align*}
		\gamma_X^g(hN) & =\gamma_X(hNg^{-1}) \\
		& = \mu_{N}^{hg^{-1}}(X\cap hg^{-1}N) \\
		& = \mu_{N}^{h}(Xg\cap hN) \\
		& = \gamma_{Xg}(hN)
		\end{align*}
		for any $hN\in G/N$.
		
		Then, we define
		\begin{displaymath}
		\mu(X)=\int_{G/N} \gamma_X\ \dd\mu_{G/N}
		\end{displaymath}
		for  $X\subseteq G$. It is now clear from the aforementioned properties of $\gamma_X$ that $\mu$ is a right-invariant finitely additive probability measure on $G$. Therefore, $G$ is amenable.

		\item 
		For each $i\in I$, let $\rho_i\colon G_i\longrightarrow G$ be the canonical homomorphism and $H_i=\rho(G_i)$, which is amenable by item (ii), and so it has a right-invariant finitely additive probability measure $\mu_i$
		. 
		
		Consider the set $\PM_i$ of all $\mu\in\PM(G)$ such that $\mu(Xh)=\mu(X)$ for all $X\subseteq G$ and $h\in H_i$. For every $i\in I$, we can set $$\bar\mu_i(X)=\mu_i(X\cap H_i)$$ for $X\subseteq G$, and so $\bar\mu_i\in\PM_i$. The set $[0,1]^{\mathcal{P}(G)}$ is compact by Tychonoff's Theorem, and $\PM_i$ is a closed subset of $[0,1]^{\mathcal{P}(G)}$ because it can be written as the intersection of preimages of closed sets by continuous functions. 

		Furthermore, given any $i,j\in I$ there exists some $k\in I$ such that $H_i,H_j\leq H_k$, and thus  $\PM_k\subseteq \PM_i\cap \PM_j$, which implies that $\PM_i \cap \PM_j \not= \emptyset$. Hence, $\{\PM_i\}_{i\in I}$ is a collection of closed subsets of the compact space $[0,1]^{\mathcal{P}(G)}$ with the finite intersection property, and so their intersection is non-empty, i.e. there exists some $\mu\in\bigcap_{i\in I}\PM_i$. Therefore, $\mu$ is a right-invariant finitely additive probability measure on $G$, and so $G$ is amenable.
	\end{enumerate}
\end{proof}

\begin{remarks}
	\begin{enumerate}
		\item As mentioned in Example~\ref{ex:Direct Limit fg subgroups}, every group can be written as the direct limit of its finitely generated subgroups. In light of Proposition~\ref{Prop:Amenable closed Sbgp Ext DirLim}, this implies that a group is amenable if and only if all of its finitely generated subgroups are amenable.
		
		\item As we saw in Example~\ref{ex:First amenable groups}, the free group of rank $2$ is not amenable. As such, no group with a free non-abelian subgroup can be amenable. It was conjectured for some time that the converse of this result was true as well. This conjecture, which came to be known as the von Neumann Conjecture, was eventually shown to be false by A. Y. Ol'shanskii in \cite{olshanskii}. 
	\end{enumerate}
\end{remarks}

\section{The F\o lner Condition}
\label{section:Folner}

We will now present a characterisation of amenability given by E. F\o lner in \cite{folner}. The so-called F\o lner condition is satisfied when a group has arbitrarily invariant finite subsets. 
This will give us another characterisation of amenable groups as those that satisfy the F\o lner condition.

\begin{definition}\label{def:FolnerCondition}
	A group $G$ is said to satisfy the \emph{F\o lner condition} if for every finite $X\subseteq G$ and every $\varepsilon>0$ there exists a finite non-empty subset $F\subseteq G$ such that 
	\begin{displaymath}
	\frac{\abs{ F\setminus Fg}}{\abs{F}}< \varepsilon
	\end{displaymath}
	for all $g\in X$.
\end{definition}


The F\o lner condition can be restated in terms of nets of almost invariant subsets. Let us now briefly recall the notion of net and some of its basic properties.
\begin{definition}
	Let $X$ be a topological space. A \emph{net} in $X$ is a family $(x_i)_{i\in I}$ of points of $X$ indexed by some directed set $I$.
	
		

	We say that the net $(x_i)_{i\in I}$ \emph{converges} to the point $x\in X$ if, for every neighbourhood $V\subseteq X$ of $x$, there is some $i_0\in I$ such that $x_i\in V$ for all $i\geq i_0$. If the limit is unique, we write
	\begin{displaymath}
	x=\lim_{i\in I}x_i.
	\end{displaymath}
\end{definition}

\begin{proposition}
	Let $X$ be a topological space. Then, the following hold:
	\begin{enumerate}
		\item The space $X$ is Hausdorff if and only if every convergent net has a unique limit point.
		
		\item The space $X$ is compact of and only if every net has a convergent subnet.
	\end{enumerate}
\end{proposition}

The F\o lner condition can then be stated in terms of the existence of a net of finite subsets that grow more and more invariant.

\begin{definition}\label{def:FolnerNet}
	A net $(F_i)_{i\in I}$ of finite non-empty subsets of a group $G$ is said to be a  \emph{F\o lner net} if 
	\begin{displaymath}
	\lim_{i\in I}\frac{\abs{F_i\setminus F_ig}}{\abs{F_i}}=0
	\end{displaymath}
	for every $g\in G$. 
	When $I=\N$, we refer to a sequence $(F_n)_{n\in \N}$ satisfying the above property as a \emph{F\o lner sequence}.
\end{definition}

\begin{examples}\label{ex:Folner Sequences}
	\begin{enumerate}
		\item If $G$ is a finite group, then the constant sequence $(F_n)_{n\in \N}$ with $F_n=G$ for all $n\in\N$ is clearly a F\o lner sequence.
		
		\item Consider the group of integers $\Z$. For each $n\in\N$, consider the finite set $$F_n=[-n,n]\cap\Z.$$ Then, for each $k\in \Z$ we have that
		\begin{displaymath}
		\frac{\abs{F_n\setminus (F_n+k)}}{\abs{F_n}}\leq \frac{\abs{k}}{2n+1}
		\end{displaymath}
		for all $n\in\N$, and so $(F_n)_{n\in \N}$ is a F\o lner sequence in $\Z$.
	\end{enumerate}
\end{examples}

\begin{remark}
	A group $G$ satisfies the F\o lner condition if and only if for every finite $X\subseteq G$ and every $\varepsilon>0$ there exists a finite non-empty subset $F\subseteq G$ such that 
	\begin{displaymath}
	\frac{\abs{F\Delta Fg}}{\abs{F}}< \varepsilon
	\end{displaymath}
	for all $g\in X$. Similarly, the net $(F_i)_{i\in I}$ is F\o lner if and only if 
	\begin{displaymath}
	\lim_{i\in I}\frac{\abs{F_i\Delta F_i g}}{\abs{F_i}}=0
	\end{displaymath}
	for every $g\in G$. We will use  these characterisations when convenient.
\end{remark}

\begin{theorem}\label{thm:Folner Condition <=> Folner Net}
	Let $G$ be a group. Then, $G$ satisfies the F\o lner condition if and only if there is a F\o lner net in $G$.
\end{theorem}

\begin{proof}
	Assume first that there is a F\o lner net $(F_i)_{i\in I}$ in $G$. Then, given $\varepsilon>0$ and a finite subset $X\subseteq G$, there exists some $i\in I$ such that
	\begin{displaymath}
	\frac{\abs{F_i\setminus F_i g}}{\abs{F_i}}< \varepsilon
	\end{displaymath}
	for all $g\in X$. Hence, $G$ satisfies the F\o lner condition.
	
	Conversely, assume that $G$ satisfies the F\o lner condition. Let $I$ be the set of pairs $(X,\varepsilon)$ with $X\subseteq G$ finite and $\varepsilon>0$. We can define a partial order $\preceq$ on $I$ by setting $(X,\varepsilon)\preceq (X',\varepsilon')$ if and only if $X\subseteq X'$ and $\varepsilon\geq \varepsilon'$. Given $(X,\varepsilon),(X',\varepsilon')\in I$, we have that 
	\begin{displaymath}
	(X,\varepsilon),(X',\varepsilon')\preceq \big(X\cup X', \min\{\varepsilon,\varepsilon'\}\big),
	\end{displaymath}
	and so $I$ is a directed set. By the F\o lner condition, for every $i\in I$ there exists some finite non-empty subset $F_i\subseteq G$ such that
	\begin{displaymath}
	\frac{\abs{F_i\setminus F_i g}}{\abs{F_i}}< \varepsilon
	\end{displaymath}
	for all $g\in X$. Hence, $(F_i)_{i\in I}$ is a F\o lner net in $G$.
\end{proof}

We will now show that every group satisfying the F\o lner condition is amenable. Later on, we will be able to prove that the converse also holds, as part of Theorem~\ref{thm:Characterisations Amenablility}.

\begin{theorem}\label{thm:FolnerCondition=>ProbabilityMeasure}
	Let $G$ be a group. If $G$ satisfies the F\o lner condition, then $G$ is amenable.
\end{theorem}

\begin{proof}
	Given any finite subset $X\subseteq G$ and $\varepsilon> 0$, denote by
	$\mathcal{PM}_{X,\varepsilon}$ 
	the set of finitely additive probability measures on $G$ such that 
	\begin{displaymath}
	\abs{\mu(Y)-\mu(Yg)}\leq\varepsilon
	\end{displaymath}
	for all $g\in X$ and  $Y\subseteq G$. 
	We have that $\mathcal{PM}_{X,\varepsilon}$ is a closed subset of $[0,1]^{\mathcal{P}(G)}$, for it can be written as the intersection of zero sets of continuous functions.
	
	Moreover, $[0,1]^{\mathcal{P}(G)}$ is compact as a consequence of Tychonoff's Theorem
	, and so $\mathcal{PM}_{X,\varepsilon}$ is compact. 
	By the F\o lner condition, there exists some finite non-empty subset $F\subseteq G$ such that
	\begin{displaymath}
	\frac{\abs{ F\setminus Fg}}{\abs{F}}< \varepsilon
	\end{displaymath}
	for all $g\in X$, so we can set $$\mu_{X,\varepsilon}(Y)=\frac{\abs{Y\cap F}}{\abs{F}}$$ for  $Y\subseteq G$. Then, $\mu_{X,\varepsilon}\in \mathcal{PM}_{X,\varepsilon}$ and the set $\mathcal{PM}_{X,\varepsilon}$ is non-empty. We also have that 
	\begin{displaymath}
	\mathcal{PM}_{X\cap X',\min\{\varepsilon,\varepsilon'\}}\subseteq\mathcal{PM}_{X,\varepsilon}\cap \mathcal{PM}_{X',\varepsilon'},
	\end{displaymath}
	and the intersection is non-empty. 
	Hence, $\{\mathcal{PM}_{X,\varepsilon}\}$ is a collection of closed non-empty subsets of $[0,1]^{\mathcal{P}(G)}$ with the finite intersection property and, because $[0,1]^{\mathcal{P}(G)}$ is compact, there must exist some $\mu\in\bigcap\mathcal{PM}_{X,\varepsilon}$. This $\mu$ is a right-invariant finitely additive probability measure on $G$, and so $G$ is amenable.
\end{proof}

We can now show, with the help of F\o lner sequences, that a number different classes of groups are amenable.

\begin{examples}
	\begin{enumerate}
		\item The group $\Z$ is amenable, for as we saw in Example~\ref{ex:Folner Sequences}, the sequence $(F_n)_{n\in \N}$ with $$F_n=[-n,n]\cap\Z$$ for each $n\in\N$ is a F\o lner sequence in $\Z$.
		
		\item Abelian groups are amenable. Indeed, every finitely generated abelian group is of the form $G=\Z^r\times H$ with $r\geq 0$ and $H$ finite. Since both $\Z$ and $H$ are amenable, and extensions of amenable groups are amenable, we have that $G$ is amenable. Finally, because amenability is closed under taking direct limits, we reach the conclusion that arbitrary abelian groups are amenable. 
		
		\item Solvable groups are amenable. Recall that a group $G$ is solvable if it has a subnormal series
		\begin{displaymath}
		1=G_0\normaleq G_1\normaleq \dotsb \normaleq  G_n=G
		\end{displaymath}
		such that the quotient $G_k/G_{k-1}$ is abelian for all $k=1,\dotsc, n$. If $n+1$ is the minimum length of any such series, we say that $G$ is solvable of class $n$. By induction on the solvability class $n$ of $G$, assume that every solvable group of class less than $n$ is amenable. Then, $G_{n-1}\normaleq G$ is solvable of class less than $n$, so it is amenable by induction. Furthermore, $G/G_{n-1}$ is also amenable by virtue of being abelian. Therefore, $G$ is an extension of amenable groups, and so it is itself amenable.
	\end{enumerate}
\end{examples}

When our group is countable, and in particular when it is finitely generated, the existence of F\o lner sequences is equivalent to satisfying the F\o lner condition.

\begin{theorem}\label{thm:FolnerSequence<=>FolnerCondition+Countable}
	A group $G$ has a F\o lner sequence if and only if $G$ satisfies the F\o lner condition and is countable.
\end{theorem}

\begin{proof}
	Suppose that $G$ satisfies the F\o lner condition and is countable. Because $G$ is countable, we can write $$G=\bigcup_{n\in \N}X_n$$ with $X_n \subseteq G$ finite and $X_n\subseteq X_{n+1}$ for all $n\in\N$. Now, because $G$ satisfies the F\o lner condition, for each $n\in\N$ there exists a finite subset $F_n\subseteq G$ such that  
	\begin{displaymath}
	\frac{\abs{F_n\setminus F_n g}}{\abs{F_n}} < \frac1n
	\end{displaymath}
	for every $g\in X_n$. From this, we deduce that 
	\begin{displaymath}
	\lim_{n\to\infty} \frac{\abs{F_n\setminus F_n g}}{\abs{F_n}} = 0
	\end{displaymath}
	for every $g\in G$, and so $(F_n)_{n\in \N}$ is a F\o lner sequence in $G$.
	
	Suppose now that $G$ has a F\o lner sequence $(F_n)_{n\in \N}$. 
	Then, $G$ satisfies the F\o lner condition by Theorem~\ref{thm:Folner Condition <=> Folner Net}. Now, for each $n\in\N$ define
	\begin{displaymath}
	X_n=\{xy^{-1}\mid x,y\in F_n\}.
	\end{displaymath}
	Given $g\in G$, there is some $N\in \N$ such that
	\begin{displaymath}
	\frac{\abs{F_n\setminus F_n g}}{\abs{F_n}}<\frac12
	\end{displaymath}
	for all $n\geq N$, implying that $F_n\cap F_n g\not=\emptyset$, and so $g\in X_n$. Therefore, $$G=\bigcup_{n\in\N}X_n$$ and, because every $X_n$ is finite, $G$ is countable.
\end{proof}

Let us now see some alternative characterisations of F\o lner sequences. For that, we will need to introduce some concepts related to invariance of subsets of a group.

\begin{definition}
	Let $G$ be a group, $F,X\subseteq G$ be non-empty finite subsets of $G$ and $\varepsilon>0$. We say that $X$ is \emph{$(F,\varepsilon)$-invariant} if
	\begin{displaymath}
	\abs[\big]{\{g\in X\mid gF\subseteq X\}}> (1-\varepsilon)\abs{X}.
	\end{displaymath}
\end{definition}

\begin{definition}
	Let $G$ be a group and $F,X\subseteq G$ be non-empty finite subsets of $G$. The \emph{$F$-boundary of $X$} is the set
	\begin{displaymath}
	\partial_F X=\big\{g\in G\mid gF\cap X\not=\emptyset \text{ and } gF\cap (G\setminus X)\not=\emptyset\big\}.
	\end{displaymath}
\end{definition}

We can now prove the following characterisations of F\o lner sequences.

\begin{proposition}\label{prop:FolnerSeq Equivalent Defs}
	Let $G$ be a countable group and $(F_n)_{n\in\N}$ be a sequence of non-empty finite subsets of $G$. Then, the following are equivalent:
	\begin{enumerate}
		\item The sequence $(F_n)_{n\in \N}$ is F\o lner, i.e. for every $g\in G$ we have that 
		\begin{displaymath}
		\lim_{n\to\infty}\frac{\abs{F_n\Delta F_n g}}{\abs{F_n}}=0.
		\end{displaymath}
		
		\item For any finite subset $F\subseteq G$ and any $\varepsilon>0$, there exists some $N\in\N$ such that $F_n$ is $(F,\varepsilon)$-invariant for every $n\geq N$.
		
		\item For any finite subset $F\subseteq G$ and any $\varepsilon>0$, there exists some $N\in\N$ such that $\abs{\partial_F F_n}< \varepsilon\abs{F_n}$ for every $n\geq N$.
	\end{enumerate}
\end{proposition}

\begin{proof}
	First, let us see that \textbf{(i) implies (iii)}. Given a finite subset $F\subseteq G$ and $\varepsilon>0$, there exists some $N\in\N$ such that 
	$$\frac{\abs{F_n\Delta F_n g}}{\abs{F_n}}< \frac{\varepsilon}{\abs{F}^2}$$
	for all $g\in FF^{-1}$. Observe that we can write
	\begin{align*}
	\partial_F F_n & =\bigg(\bigcup_{s\in F}F_n s^{-1}\bigg)\setminus\bigg(\bigcap_{s\in F}F_n s^{-1}\bigg) \\
	& = \bigcup_{s,t\in F} (F_n s^{-1}\Delta F_nt^{-1}),
	\end{align*}
	and so
	\begin{align*}
	\abs{\partial_F F_n} & =\abs[\bigg]{\bigcup_{s,t\in F} (F_n s^{-1}\Delta F_n t^{-1})} \\
	& \leq \sum_{s,t\in F}\abs{F_n\Delta F_n t^{-1}s} \\
	& < \varepsilon\abs{F_n}
	\end{align*}
	for every $n\geq N$.
	
	\smallskip
	
	Let us now prove that \textbf{(iii) implies (ii)}. Given a finite subset $F\subseteq G$ and $\varepsilon>0$, if we take the set $F'=F\cup\{1\}$, there exists some $N\in\N$ such that $\abs{\partial_{F'} F_n}< \varepsilon\abs{F_n}$ for every $n\geq N$. Assume by contradiction that $F_n$ is not $(F',\varepsilon)$-invariant, i.e. 
	\begin{displaymath}
	\abs[\big]{\{s\in F_n\mid sF'\subseteq F_n\}}\leq (1-\varepsilon)\abs{F_n}.
	\end{displaymath}
	Because $1\in F'$, we can write
	\begin{displaymath}
	\{s\in F_n\mid sF'\subseteq F_n\}=\bigcap_{s\in F'}(F_n\cap F_n s^{-1})=\bigcap_{s\in F'}F_n s^{-1},
	\end{displaymath}
	and so
	\begin{align*}
	\abs{\partial_{F'}F_n} & =\abs[\bigg]{\bigcup_{s\in F'} F_n s^{-1}}-\abs[\bigg]{\bigcap_{s\in F'} F_n s^{-1}} \\
	& \geq  \abs{F_n}-(1-\varepsilon)\abs{F_n} \\
	& = \varepsilon\abs{F_n}.
	\end{align*}
	Therefore, $F_n$ must be $(F',\varepsilon)$-invariant for every $n\geq N$ and, because $F\subseteq F'$ and
	\begin{displaymath}
	\{s\in F_n\mid sF'\subseteq F_n\}\subseteq \{s\in F_n\mid sF\subseteq F_n\},
	\end{displaymath}
	this implies that $F_n$ is $(F,\varepsilon)$-invariant for every $n\geq N$.
	
	\smallskip
	
	Finally, let us show that \textbf{(ii) implies (i)}. Given $g\in G$ and $\varepsilon>0$, there exists some $N\in\N$ such that $F_n$ is $(\{g^{-1}\},\frac{\varepsilon}{2})$-invariant for every $n\geq N$. 
	Now, we have that
	\begin{align*}
		2\abs{F_n\cap F_n g} & = \big(\abs{F_n}-\abs{F_n\setminus F_n g}\big)+\big(\abs{F_n g}-\abs{F_n g\setminus F_n}\big) \\
		& = 2\abs{F_n}-\abs{F_n\Delta F_n g}.
	\end{align*}
	Thus,
	\begin{align*}
		\Big(1-\frac{\varepsilon}{2}\Big)\abs{F_n} 
		& < \abs[\big]{\{s\in F_n\mid s g^{-1}\in F_n\}} \\
		& = \abs{F_n\cap F_n g} \\
		& = \abs{F_n}-\frac12\abs{F_n\Delta gF_n}, 
	\end{align*}
	from where we obtain that
	\begin{displaymath}
	\frac{\abs{F_n\Delta F_n g}}{\abs{F_n}}< \varepsilon
	\end{displaymath}
	for every $n\geq N$. 
	Therefore, $(F_n)_{n\in\N}$ is a F\o lner sequence.
\end{proof}

A particular type of F\o lner sequence is what we will call \emph{F\o lner exhaustion}, i.e. a F\o lner sequence $(F_n)_{n\in\N}$ in the group $G$ such that 
\begin{displaymath}
1\in F_1\subseteq \dotsb \subseteq F_n \subseteq \dotsb
\end{displaymath}
and $$G=\bigcup_{n\in\N}F_n.$$ We will now see that the existence of F\o lner sequences is equivalent to the existence of F\o lner exhaustions.

\begin{proposition}
	Let $G$ be a countable group. Then, $G$ has a F\o lner sequence if and only if it has a F\o lner exhaustion.
\end{proposition}

\begin{proof}
	Every F\o lner exhaustion is by definition a F\o lner sequence. Thus, we only need to show that whenever we have a F\o lner sequence we can obtain a F\o lner exhaustion. 
	
	Let $(F_n)_{n\in\N}$ be a F\o lner sequence in $G$. First, we will see that we can obtain from $(F_n)_{n\in\N}$ a nested F\o lner sequence, i.e. a F\o lner sequence $(F_{k}')_{k\in \N}$ such that 
	\begin{displaymath}
	1\in F_1'\subseteq \dotsb \subseteq F_k' \subseteq \dotsb.
	\end{displaymath} 
	Without loss of generality, assume that $1\in F_1$, and take $F_1'=F_1$. Suppose by induction that we have constructed finite subsets $F_1'\subseteq \dotsb \subseteq F_{k-1}'$ of $G$. Because $(F_n)_{n\in\N}$ is a F\o lner sequence, by Proposition~\ref{prop:FolnerSeq Equivalent Defs} there is some ${n_k}\in\N$ such that $F_{n_k}$ is $(F_{k-1}',1)$-invariant, i.e. 
	\begin{displaymath}
	\abs[\big]{\{g\in F_{n_k} \mid g F_{k-1}'\subseteq F_{n_k}\}}>0,
	\end{displaymath}
	and so there exists some $g_k\in F_{n_k}$ such that $g_kF_{k-1}' \subseteq F_{n_k}$. If we define $F_{k}'=g_k^{-1}F_{n_k}$, then $F_{k-1}'\subseteq F_{k}'$. Furthermore, given any $g\in G$ we have that
	\begin{displaymath}
	\abs{F_k'\setminus F_k' g}=\abs[\big]{g_k(F_{n_k}\setminus F_{n_k}g)}=\abs{F_{n_k}\setminus F_{n_k}g},
	\end{displaymath}
	and so the sequence $(F_{k}')_{k\in \N}$ that we have constructed is a nested F\o lner sequence.

	Assume now that $(F_n)\ninN$ is a nested F\o lner sequence in $G$. Because $G$ is countable, we can write $$G=\bigcup_{n\in \N}X_n$$ with $X_n \subseteq G$ finite and $X_k\subseteq X_{k+1}$ for all $k\in\N$. Define now $F_n^r=F_n X_r$ for each $n,r\in \N$. Observe that $F_n\subseteq F_n^r$, and so $\abs{F_n}\leq \abs{F_n^r}$ for any $n,r\in\N$. Then, for any finite subset $F\subseteq G$ we have that 
	\begin{displaymath}
	\abs{\partial_F F_n^r}\leq \sum_{g\in X_r}\abs[\big]{\partial_F(F_n g)}\leq \abs{X_r}\abs{\partial_F F_n},
	\end{displaymath}
	and so 
	\begin{displaymath}
	\frac{\abs{\partial_F F_n^r}}{\abs{F_n^r}}\leq \abs{X_r}\frac{\abs{\partial_F F_n}}{\abs{F_n}}
	\end{displaymath}
	for every $n,r\in\N$. Now, for each $r\in\N$ take $n_r\in\N$ such that 
	\begin{displaymath}
	\frac{\abs{\partial_F F_{n_r}}}{\abs{F_{n_r}}}<\frac{1}{r\abs{X_r}}
	\end{displaymath}
	and $n_r\geq n_{r-1}$ for $r>1$.
	Thus, if we set $F_r'=F_{n_r}^r$ for each $r\in \N$, we have that
	\begin{align*}
	\abs{\partial_F F_r'} & \leq\abs{X_r}\abs{\partial_F F_{n_r}} \\
	& <\frac{\abs{F_{n_r}}}{r} \\
	& \leq \frac{\abs{F_{r}'}}{r}.
	\end{align*}
	Therefore, 
	\begin{displaymath}
	\lim_{r\to\infty}\frac{\abs{\partial_F F_{r}'}}{\abs{F_{r}'}}=0
	\end{displaymath}
	for any finite subset $F\subseteq G$, meaning that $(F_r')_{r\in \N}$ is a F\o lner sequence in $G$. Furthermore, it is a nested sequence because $(F_n)_{n\in\N}$ is nested. Finally,  we have that $X_r\subseteq F_r'$ for every $r\in \N$, and so
	\begin{displaymath}
	G=\bigcup_{r\in \N}F_r',
	\end{displaymath}
	which implies that $(F_r')_{r\in \N}$ is a F\o lner  exhaustion.
\end{proof}

\section{Paradoxical Decompositions
}
\label{section:ParadoxicalDecompositions}

The characterisation of amenability that we will study in this section is also intimately related to the Banach-Tarski Paradox. Essentially, we will characterise amenable groups as those for which a Banach-Tarski-like paradox cannot happen, i.e. the pieces of any finite decomposition of an amenable group cannot be rearranged  in such a way that we obtain two copies of the group. 

\begin{definition}
	%
	Let $G$ be a group acting on a set $\Omega$. Then, the action of $G$ on $\Omega$ is said to be \emph{paradoxical}, and $\Omega$ is said to be \emph{$G$-paradoxical}, if there exist pairwise disjoint subsets $X_1,\dotsc,X_n$ and $Y_1,\dotsc, Y_m$ of $\Omega$, and elements $g_1,\dotsc,g_n$ and $h_1,\dotsc,h_m$ in $G$ such that
	\begin{displaymath}
	\Omega= \bigg(\bigcupdot_{i=1}^n X_i\bigg) \cupdot \bigg( \bigcupdot_{j=1}^m Y_j\bigg) = \bigcupdot_{i=1}^n X_ig_i = \bigcupdot_{i=1}^m Y_jh_j.
	\end{displaymath}
	In that case, we also say that $\Omega$ has a \emph{$G$-paradoxical decomposition}. 
	The group $G$ is said to be \emph{paradoxical} if the action of $G$ on itself by right multiplication is paradoxical.
\end{definition}

 Using the terminology we have just introduced, amenable groups can be characterised as those that are non-paradoxical, as we will show later.

We will now see that the requirements in the definition of paradoxical decompositions can be relaxed.

\begin{proposition}\label{prop:Paradoxical Equivalent Defs}
	Let $G$ be a group acting on a set $\Omega$. Then, the following are equivalent:
	\begin{enumerate}
		\item  There exist pairwise disjoint subsets $X_1,\dotsc,X_n$ and $Y_1,\dotsc, Y_m$ of $\Omega$, and elements $g_1,\dotsc,g_n$ and $h_1,\dotsc,h_m$ in $G$ such that
		\begin{displaymath}
		\Omega= \bigg(\bigcupdot_{i=1}^n X_i\bigg) \cupdot \bigg( \bigcupdot_{j=1}^m Y_j\bigg) = \bigcupdot_{i=1}^n X_ig_i = \bigcupdot_{i=1}^m Y_jh_j.
		\end{displaymath}
		
		\item  There exist pairwise disjoint subsets $X_1,\dotsc,X_n$ and $Y_1,\dotsc, Y_m$ of $\Omega$, and elements $g_1,\dotsc,g_n$ and $h_1,\dotsc,h_m$ in $G$ such that
		\begin{displaymath}
		\Omega= \bigcupdot_{i=1}^n X_ig_i = \bigcupdot_{i=1}^m Y_jh_j.
		\end{displaymath}
		
		\item  There exist pairwise disjoint subsets $X_1,\dotsc,X_n$ and $Y_1,\dotsc, Y_m$ of $\Omega$, and elements $g_1,\dotsc,g_n$ and $h_1,\dotsc,h_m$ in $G$ such that
		\begin{displaymath}
		\Omega= \bigcup_{i=1}^n X_ig_i = \bigcup_{i=1}^m Y_jh_j.
		\end{displaymath}
	\end{enumerate}
\end{proposition}
\begin{proof}
	The fact that \textbf{(i) implies (iii)} is trivial. 
	
	\smallskip
	
	Let us show that \textbf{(iii) implies (ii)}. Assume that there exist pairwise subsets $X_1,\dotsc,X_n$ and $Y_1,\dotsc, Y_m$ of $\Omega$, and elements $g_1,\dotsc,g_n$ and $h_1,\dotsc,h_m$ in $G$ such that
	\begin{displaymath}
	\Omega= \bigcup_{i=1}^n X_ig_i = \bigcup_{i=1}^m Y_jh_j.
	\end{displaymath}
	Without loss of generality, we may assume that $g_1=h_1=1$. Take $X_1'=X_1$ and define inductively $$X_k'=X_k\setminus \bigg(\bigcup_{i=1}^{k-1}X_i'g_i\bigg)g_k^{-1}$$
	for $k=2,\dotsc,n$. 
	Similarly, take $Y'_1=Y_1$ and define inductively $$Y_k'=Y_k\setminus \bigg(\bigcup_{j=1}^{k-1}Y_j'h_j\bigg)h_k^{-1}$$
	for $k=2,\dotsc,m$. We can check that the sets $X'_1,\dotsc,X'_n$ and $Y'_1,\dotsc, Y'_m$ are pairwise disjoint, and 
	\begin{displaymath}
	\Omega= \bigcupdot_{i=1}^n X'_ig_i = \bigcupdot_{i=1}^m Y'_jh_j.
	\end{displaymath}
	
	\smallskip
	
	Finally, let us see that \textbf{(ii) implies (i)}. Assume that there exist pairwise disjoint subsets $X_1,\dotsc,X_n$ and $Y_1,\dotsc, Y_m$ of $\Omega$, and elements $g_1,\dotsc,g_n$ and $h_1,\dotsc,h_m$ in $G$ such that
	\begin{displaymath}
	\Omega= \bigcupdot_{i=1}^n X_ig_i = \bigcupdot_{i=1}^m Y_jh_j.
	\end{displaymath}
	Without loss of generality, we may assume that $h_1=1$. Write $$X=\bigcup_{i=1}^n X_i, \quad Y=\bigcup_{j=1}^m Y_j.$$ Observe that $X\cap Y=\emptyset$. Now, given any $\alpha\in \Omega$ there exist a unique $j\in\{1,\dotsc,m\}$ and some $f(\alpha)\in Y_j$ such that $\alpha=f(\alpha)h_j$. This defines a map $f\colon \Omega\longrightarrow Y$. Let $$Z=X\cup\bigg(\bigcup_{k\in\N}f^k(X)\bigg),\quad Z_0=(G\setminus X)\setminus f(Z).$$ Then, we have that $$X\cap f(Z)=\emptyset,\quad X\cup f(Z)=Z.$$ Moreover, if we put $Z_j=Y_j\cap Z h_j^{-1}$, since $h_1=1$, we obtain that
	\begin{align*}
	\Omega & =X\cupdot  \Bigg((Z_0\cupdot Z_1) \cupdot \bigg(\bigcupdot_{j=2}^m Z_j \bigg) \Bigg) \\
	& = \bigcupdot_{i=1}^n X_ig_i \\
	& = (Z_0\cupdot Z_1)h_1\cupdot \bigg(\bigcupdot_{j=2}^mZ_j h_j\bigg).
	\end{align*}
\end{proof}

\begin{remark}
	As we can see in the proof of Proposition~\ref{prop:Paradoxical Equivalent Defs}, the number of pieces is preserved when we go from one type of decomposition to another. This allows us to define the \emph{Tarski number} of a $G$-set $\Omega$ as the smallest number of pieces of any $G$-paradoxical decomposition of $\Omega$.
\end{remark}

\begin{example}
	Consider the free $F(a,b)$ on two generators. As we saw in Examples~\ref{ex:First amenable groups} (ii), we can write 
	\begin{align*}
		F(a,b) & = \{1\}\cupdot W(a)\cupdot W(a^{-1}) \cupdot W(b) \cupdot W(b^{-1}) \\
		& = W(a) \cupdot W(a^{-1})a \\
		& = W(b) \cupdot W(b^{-1})b,
	\end{align*}
	where $W(s)$ is the set of reduced words ending with $s\in\{a^{\pm 1}, b^{\pm 1}\}$. Therefore, $F(a,b)$ is paradoxical. Furthermore, it is clear that any paradoxical decomposition must have at least 4 pieces, and so the Tarski number of $F(a,b)$ is 4. It can actually be shown that a group has Tarski number 4 if and only if it contains a subgroup isomorphic to $F(a,b)$, see \cite[Theorem 5.8.38]{sapir}.
\end{example}

We will now show that whether a group is paradoxical is entirely dependent on whether it has paradoxical free actions.

\begin{theorem}\label{thm:Paradoxical Action Group  Equivalent}
	Let $G$ be a group. Then, the following are equivalent:
	\begin{enumerate}
		\item The group $G$ is paradoxical.

		\item Every free action of $G$ is paradoxical.
		
		\item There exists a paradoxical free action of $G$. 
	\end{enumerate}
\end{theorem}
\begin{proof}
	First, let us show that \textbf{(i) implies (ii)}. Assume that there exist pairwise disjoint subsets $X_1,\dotsc,X_n$ and $Y_1,\dotsc, Y_m$ of $G$, and elements $g_1,\dotsc,g_n$ and $h_1,\dotsc,h_m$ in $G$ such that
	\begin{displaymath}
	G= \bigcup_{i=1}^n X_ig_i = \bigcup_{i=1}^m Y_jh_j.
	\end{displaymath}
	Let $\Omega$ be a set on which $G$ acts freely. Using the Axiom of Choice, we can select a set $T\subseteq \Omega$ of representatives of the orbits of $\Omega$ under the action of $G$. Then, we can write 
	\begin{displaymath}
	\Omega = \bigcupdot_{g\in G} Tg,
	\end{displaymath}
	for if $\alpha g=\beta h$ for some $\alpha,\beta\in T$ and $g,h\in G$, then $\alpha=\beta$ by the definition of $T$, and the action being free implies that $g=h$. Now, define $$\tilde X_i=\bigcupdot_{g\in X_i}Tg, \quad \tilde Y_j=\bigcupdot_{g\in Y_j} Tg$$ for $i=1,\dotsc,n$ and $j=1,\dotsc,m$. Then, the $\tilde X_1,\dotsc,\tilde X_n$ and $\tilde Y_1,\dotsc, \tilde Y_m$ are pairwise disjoint, and
	\begin{displaymath}
	\Omega= \bigcup_{i=1}^n \tilde X_ig_i = \bigcup_{i=1}^m \tilde Y_jh_j.
	\end{displaymath}
	Therefore, the action of $G$ on $\Omega$ is paradoxical. 
	
	\smallskip 
	
	It is trivial that \textbf{(ii) implies (iii)}, for the action of $G$ on itself by right multiplication is free. 
	
	\smallskip
	
	Finally, let us prove that \textbf{(iii) implies (i)}. Assume that there is a free paradoxical action of $G$ on some set $\Omega$. Then, if we fix an element $\alpha\in\Omega$, the action of $G$ on the orbit $\alpha G$ must also be paradoxical and, because $G$ acts freely on $\Omega$, the action of $G$ on $\alpha G$ is equivalent to the action of $G$ on itself by right multiplication. Thus, $G$ itself is paradoxical.
\end{proof}

\section{The Banach-Tarski Paradox}

We will now prove the Banach-Tarski Paradox, which states that the closed unit ball in the euclidean space $\R^3$ can be decomposed into a finite number of pieces that can then be rearranged using only isometries of $\R^3$. The proof of the paradox relies on the paradoxicality of the free group of rank $2$. The group of rotations of $\R^3$ contains a free subgroup of rank $2$, which produces a paradoxical decomposition of the unit sphere. This decomposition of the unit sphere can then be extended to a paradoxical decomposition of the whole unit ball.

Recall that $\SO(3)$ is the group of rotations about the origin in $\R^3$ under composition, and is identified with the group of orthogonal $3\times 3$ real matrices with determinant $1$ under matrix multiplication. We will also need to consider $\E(3)$, the group of isometries of the euclidean space $\R^3$. 

Throughout the rest of this section, 
we will denote the unit sphere centred at the origin in $\R^3$ by $\SSS^2$, and the closed unit ball centred at the origin in $\R^3$ by $\B^3$.

The key fact in the proof of the Banach-Tarski is the following result. 

\begin{proposition}
	The group $\SO(3)$ contains a subgroup $H$ which is isomorphic to the free group $F(a,b)$.
\end{proposition}


\begin{proof}
	Consider the matrices $A,B\in\SO(3)$ given by \begin{displaymath}
	A=\frac{1}{7}\begin{pmatrix}
	6 & 2 & -3 \\
	2 & 3 & 6 \\
	3 & -6 & 2
	\end{pmatrix}, \quad
	B=\frac{1}{7}\begin{pmatrix}
	2 & 6 & -3 \\
	-6 & 3 & 2 \\
	3 & 2 & 6
	\end{pmatrix},
	\end{displaymath}
	and the group $H=\scalar{A,B}\leq \SO(3)$. 
	Let $w\in F(a,b)$ be a non-trivial reduced word. We will now show that $w(A,B)\not= I$, thus proving that $H\isom F(a,b)$. For the sake of simplicity, we will write $w=w(A,B)$.

	We may assume without loss of generality that $w$ begins with $A$, otherwise conjugate $w$ by a sufficiently high power of $A$ and invert if necessary. Then, we can write $w=AA^{k_1}B^{\pm k_2}\dotsm A^{\pm k_t}$ with $k_i\geq 0$ for all $i=1,\dotsc,t$. 
	
	Write $\bar A_{\pm},\bar B_{\pm}$ for the reductions modulo $7$ of the matrices $7A^{\pm 1},7B^{\pm 1}$, respectively. Then, if we put $\bar w= \bar A_+\bar A_+^{k_1}\bar B_{\pm}^{ k_2}\dotsm \bar A_{\pm}^{ k_t}$, it is enough to show that $(1,0,0)\bar w\not=(1,0,0)$. Define
	\begin{align*}
		V_{\bar A_+}&=\{(3,1,2),(5,4,1),(6,2,4)\}, \\
		V_{\bar A_{-}}&=\{(3, 2, 6), (5, 1, 3), (6, 4, 5)\}, \\
		V_{\bar B_+}&=\{(3, 5, 1), (5, 6, 4), (6, 3, 2)\}, \\
		V_{\bar B_{-}}&=\{(1, 5, 4), (2, 3, 1), (4, 6, 2)\}.
	\end{align*}
	Firstly, we have that $$\bar A_+(1,0,0)=(6,2,4)\in V_{\bar A_+}.$$ Doing matrix computations, we can see that the following hold:
	\begin{enumerate}
		\item If $$v\in V_{\bar A_+}\cup V_{\bar B_+}\cup V_{\bar B_-},$$ then $\bar A_+v\in V_{\bar A_+}$.
		
		\item If $$v\in V_{\bar A_-}\cup V_{\bar B_+}\cup V_{\bar B_-},$$ then $\bar A_-v\in V_{\bar A_-}$.
		
		\item If $$v\in V_{\bar B_+}\cup V_{\bar A_+}\cup V_{\bar A_-},$$ then $\bar B_+v\in V_{\bar B_+}$.
		
		\item If $$v\in V_{\bar B_-}\cup V_{\bar A_+}\cup V_{\bar A_-},$$ then $\bar B_-v\in V_{\bar B_-}$.
	\end{enumerate}
	Now, $\bar A_{+}(1,0,0)\in V_{\bar A_+}$, so  $\bar A_+^{k_1}(1,0,0)\in V_{\bar A_+}$. Then, multiplying by $\bar B_{\pm}^{ k_2}$ we arrive at $V_{\bar B_+}\cup V_{\bar B_-}$, and the next multiplication takes us to $V_{\bar A_+}\cup V_{\bar A_-}$. As we move right through $\bar w$, at each step we are either in $V_{\bar A_+}\cup V_{\bar A_-}$ or in $V_{\bar B_+}\cup V_{\bar B_-}$, which means that $$\bar w(1,0,0)\in V_{\bar A_+}\cup V_{\bar A_-}\cup V_{\bar B_+}\cup V_{\bar B_-},$$ and so $\bar w(1,0,0)\not=0$.
\end{proof}

In order to prove the Banach-Tarski Paradox, we will use the following result, known as the Hausdorff Paradox.

\begin{theorem}[Hausdorff]
	There exists a countable subset $X\subseteq \SSS^2$ such that $\SSS^2\setminus X$ is $\SO(3)$-paradoxical.
\end{theorem}
\begin{proof}
	Every non-trivial rotation in $\SO(3)$ fixes two antipodal points in $\SSS^2$. Consider the set $X\subseteq \SSS^2$ of all points fixed by some rotation in $H\leq\SO(3)$, which is countable because $H\cong F(a,b)$ is finitely generated. Then, the paradoxical group $H$ acts freely on $\SSS^2\setminus X$, and so $\SSS^2\setminus X$ is paradoxical by Theorem~\ref{thm:Paradoxical Action Group  Equivalent}.
\end{proof}

\begin{definition}
	Let $G$ be a group acting on a set $\Omega$. We say that two subsets $X,Y\subseteq \Omega$ are \emph{$G$-equidecomposable}, and write $X\sim Y$, if there exist subsets $X_1,\dotsc,X_n\subseteq X$ and $Y_1,\dotsc,Y_n\subseteq Y$ with  $$X=\bigcupdot_{i=1}^n X_i, \qquad Y=\bigcupdot_{i=1}^n Y_i,$$ and elements $g_1,\dotsc,g_n\in G$ such that $Y_i=X_ig_i$ for all $i=1,\dotsc,n$.
\end{definition}

\begin{remarks}
	\begin{enumerate}
		\item It is easy to see that being $G$-equidecomposable is an equivalence relation on the family of subsets of $\Omega$.
		
		\item The condition of $\Omega$ being $G$-paradoxical can be reformulated by saying that there exist subsets $X,Y\subseteq \Omega$ such that $X\sim \Omega \sim Y$.
		
		\item Clearly, if $X$ is $G$-paradoxical and $X\sim Y$, then $Y$ is $G$-paradoxical as well. 
	\end{enumerate}
\end{remarks}

\begin{proposition}
	Given a countable subset $X\subseteq \SSS^2$, then we have that $\SSS^2\setminus D$ is $\SO(3)$-equidecomposable to $\SSS^2$.
\end{proposition}
\begin{proof}
	Because $X$ is countable, there is some line $L\subseteq \R^3$ going through the origin such that $L\cap X=\emptyset$. 
	Consider now the set $\Gamma$ of all angles $\theta \in [0,2\pi)$ such that, if we denote by $\rho_{\theta}$ the rotation about $L$ of angle $\theta$, we have that  $x\rho_{n\theta}\in X$ for some $n\in\N$ and some $x\in X$. Then, $\Gamma$ is countable, and so there is some angle $\theta\in [0,2\pi)$ such that $X\rho_{n\theta}\cap X=\emptyset$ for any $n\in \N$. If we consider the set $$\bar{X}=\bigcup_{n=0}^{\infty}X\rho_{n\theta},$$
	 we have that 
	\begin{align*}
	\SSS^2 & = \bar X\cup (\SSS^2 \setminus \bar X) \\
	& \sim \bar X\rho \cup (\SSS^2\setminus \bar X) \\
	& = (\bar X \setminus X) \cup (\SSS^2 \setminus\bar X) \\
	& = \SSS^2\setminus X.
	\end{align*}
\end{proof}

\begin{corollary}[Banach-Tarski]
	The sphere $\SSS^2$ is $\SO(3)$-paradoxical.
\end{corollary}

Connecting every point on $\SSS^2$ with a half-open segment to the origin, the paradoxical decomposition of $\SSS^2$  yields a paradoxical decomposition of the unit ball without the origin. 

\begin{corollary}\label{cor:B-T Paradox punct ball paradoxical}
	The punctured ball $\B^3\setminus\{0\}$ is $\SO(3)$-paradoxical.
\end{corollary}


There is just one final step left in order to prove the Banach-Tarski Paradox.

\begin{proposition}\label{prop:B-T Paradox punct ball equidec}
	The punctured ball $\B^3\setminus\{0\}$ is $\E(3)$-equidecomposable to $\B^3$.
\end{proposition}
\begin{proof}
	Let $\rho\in\E(3)$ be a rotation of infinite order about an axis crossing $\B^3$ but without going through the origin, and such that $0\cdot \rho^n\in \B^3$ for all $n\in\N$. Then, if we take $X=\{0\}$ and $$\bar{X}=\{0\cdot \rho^n\mid n\geq 0\},$$ we have that 
	\begin{align*}
	\B^3 & =\bar X\cup(\B^3\setminus \bar X) \\
	& \sim \bar X\rho \cup (\B^3\setminus\bar X) \\
	& = \B^3\setminus \{0\}.
	\end{align*}
\end{proof}

Finally, combining the previous results we obtain the Banach-Tarski Paradox. 

\begin{theorem}[Banach-Tarski]
	The ball $\B^3$ is $\E(3)$-paradoxical. 
\end{theorem}
\begin{proof}
	By Proposition~\ref{prop:B-T Paradox punct ball equidec}, the ball $\B^3$ is $\E(3)$-equidecomposable to the punctured ball $\B^3\setminus\{0\}$, which is in turn $\E(3)$-paradoxical due to Corollary~\ref{cor:B-T Paradox punct ball paradoxical}. Therefore, we can conclude that $\B^3$ is $\E(3)$-paradoxical.
\end{proof}



\section{Ultrafilters, the Stone-\u{C}ech Compactification and Fixed Point Properties}


The concept of amenability can be further characterised by the fixed points of certain kinds of actions of our group on some spaces. One such characterisation says that a group is amenable if and only if every affine continuous action of the group on a non-empty convex compact subset of a Hausdorff topological vector space has a fixed point.

Amenable groups can also be characterised as those whose every continuous action on a non-empty compact  Hausdorff topological space fixes some Borel probability measure.

In order to be able to prove these characterisations, we will introduce the concepts of filters and ultrafilters and the Stone-\u{C}ech compactification of a discrete topological space. We will also make use of ultrafilters in the following chapters. For a more through exposition of the topics of ultrafilters and the Stone-\u{C}ech compactification, see \cite{hindman}.





\begin{definition}\label{def:Filter}
	A \emph{filter} on a set $\Omega$ is a collection $\F$ of subsets of $\Omega$ satisfying the following properties:
	\begin{enumerate}
		\item $\Omega\in\F$ and $\emptyset\not\in\F$.
		
		\item If $X\in\F$ and $X\subseteq Y$, then $Y\in\F$.
		
		\item If $X,Y\in \F$, then $X\cap Y\in\F$.
	\end{enumerate} 
	An \emph{ultrafilter} on $\Omega$ is a maximal filter, i.e. a filter that is not properly contained in any other filter on $X$.
\end{definition}

\begin{examples}
	\begin{enumerate}
		\item If $\Omega$ is a topological space, then given any point $x\in \Omega$ the set $\mathcal{N}_x$ of all neighbourhoods of $x$ is a filter on $\Omega$. 
		
		\item Given an element $x\in \Omega$, we can define the ultrafilter
		\begin{displaymath}
		\F_x=\{X\subseteq \Omega\mid x\in X\},
		\end{displaymath}
		which is called the \emph{principal ultrafilter based on $x$}.
	\end{enumerate}
\end{examples}


We can talk about convergence along filters on topological spaces. Given a filter $\F$ on a topological space $\Omega$ and a point $x\in \Omega$, we say that $\F$ \emph{converges} to $x$ if $\mathcal{N}_x\subseteq \F$. We then have the following properties. 

\begin{proposition}
	Let $\Omega$ be a topological space. Then, the following hold:
	\begin{enumerate}
		\item The space $\Omega$ is Hausdorff if and only if every convergent filter on $\Omega$ has a unique limit.
		
		\item The space $\Omega$ is compact if and only if every ultrafilter on $\Omega$ is convergent.
	\end{enumerate}
\end{proposition}


Filters also allow us to generalise the notion of limit of a function. Given a set $\Omega$, a topological space $\Upsilon$ and a filter $\F$ on $\Omega$, we say that a map $f\colon \Omega\longrightarrow \Upsilon$ \emph{converges} to the point $y\in \Upsilon$ along $\F$ if $f^{-1}(V)\in\F$ for every $V\in \mathcal{N}_y$. If the limit is unique, we write
\begin{displaymath}
y=\lim_{x\to \F}f(x).
\end{displaymath}

\begin{example}
	Let $(x_n)_{n\in \N}$ be a sequence in the topological space $\Omega$. Then, $(x_n)_{n\in \N}$ converges to the point $x\in \Omega$ in the usual sense if and only if it converges along the filter $$
	\{X\subseteq\N\mid \N\setminus X\text{ is finite}\}$$ 
	on $\N$. 

	Further suppose that $(x_n)_{n\in \N}$ is bounded. Then, we have that $(x_n)_{n\in \N}$ is convergent in the usual sense with 
	\begin{displaymath}
	\lim_{n\to \infty}x_n=x
	\end{displaymath}
	if and only if
	\begin{displaymath}
	\lim_{n\to \omega}x_n=x
	\end{displaymath}
	for every non-principal ultrafilter $\omega$ on $\N$.
\end{example}

\begin{proposition}
	Let $\Omega$ be a set, $\Upsilon$ a compact topological space and $\omega$ an ultrafilter on $\Omega$. Then, a map $f\colon \Omega\longrightarrow \Upsilon$ has a limit which is unique.
\end{proposition}

The concept of ultrafilter now allows us to define the Stone-\u{C}ech compactification of a discrete topological space.

\begin{definition}\label{def:SCCompactification}
	Let $\Omega$ be a discrete topological space. The set of all ultrafilters on $\Omega$ is called the \emph{Stone-\u{C}ech compactification} of $\Omega$, and is denoted by $\beta \Omega$.
\end{definition}

Given $X\subseteq \Omega$ non-empty, we can consider
\begin{displaymath}
\beta X=\{\F\in\beta \Omega\mid X\in \F\}\subseteq \beta \Omega.
\end{displaymath}
This set can be naturally identified with the Stone-\u{C}ech compactification of $X$, which justifies our abuse of notation.

\begin{proposition} Let $\Omega$ be a discrete topological space and $X,Y\subseteq \Omega$. Then, the following properties hold:
	\begin{enumerate}
		\item Given $X,Y\subseteq \Omega$, we have that $$\beta(X\cap Y)=\beta{X}\cap\beta{Y}.$$
		
		\item Given $X,Y\subseteq \Omega$, we have that $$\beta(X\cup Y)=\beta{X}\cup\beta{Y}.$$
		
		\item Given $X\subseteq \Omega$, we have that $$\beta(\Omega\setminus X)=\beta \Omega\setminus \beta{X}.$$
	\end{enumerate}
\end{proposition}

\begin{proof}
	Let us first prove (i). Given $X,Y\subseteq \Omega$ and $\omega\in\beta \Omega$, we have that $\omega\in\beta(X\cap Y)$ if and only if $X\cap Y\in \omega$, which is in turn equivalent to $X,Y\in\omega$. But this is precisely the condition that $\omega\in\beta X\cap \beta Y$.
	
	\smallskip
	
	Let us now show that (ii) holds. Given $X,Y\subseteq \Omega$ and $\omega\in\beta X$, we have that $\omega\in\beta(X\cup Y)$ if and only if $X\cup Y\in \omega$. 
	Assume by contradiction that $X,Y\not\in \omega$. Then, we must have that $\Omega\setminus X,\Omega\setminus Y\in\omega$, leading us to deduce that 
	\begin{displaymath}
	\Omega\setminus(X\cup Y)=(\Omega\setminus X)\cap(\Omega\setminus Y)\in\omega,
	\end{displaymath}
	which implies that $X\cup Y\not\in\omega$. Therefore, $X\cup Y\in \omega$ is equivalent to having $X\in\omega$ or $Y\in\omega$. But this is precisely the condition that $\omega\in \beta{X}\cup\beta{Y}$.
	
	\smallskip 
	
	Finally, let us  prove (iii). Let $X\subseteq \Omega$ and $\omega\in\beta \Omega$. Because $\omega$ is an ultrafilter, it is easy to see that either $X\in\omega$ or $\Omega\setminus X\in\omega$, and the two possibilities are mutually exclusive. This implies that $\omega\in \beta X$ if and only if $\omega\not\in\beta(\Omega\setminus X)$. 
	%
\end{proof}

The above result shows that the family
\begin{displaymath}
\{\beta X\mid X\subseteq \Omega\}
\end{displaymath}
forms the basis for a topology on $\beta \Omega$. The Stone-\u{C}ech compactification of a discrete space $\Omega$ is thus the largest compact Hausdorff space into which $\Omega$ can be embedded as a dense subset, as can be gleaned from its universal property.

\begin{theorem}\label{thm:SC Compactification Universal Property}
	Let $\Omega$ be a discrete topological space. Then, $\beta \Omega$ is a compact Hausdorff topological space containing $\Omega$ as a dense subset. Furthermore, if $\Upsilon$ is a compact Hausdorff space, any continuous map $f\colon \Omega\longrightarrow \Upsilon$ admits a unique continuous extension $\beta{f}\colon \beta \Omega\longrightarrow Y$.
\end{theorem}

\begin{proof}
	We can identify $\Omega$ with the subspace of $\beta \Omega$ formed by the principal ultrafilters, i.e. 
	\begin{displaymath}
	\Omega=\{\omega_x\mid x\in \Omega\}.
	\end{displaymath}
	Then, given any non-empty subset $X\subseteq \Omega$ and a point $x\in X$ we have that $\omega_x\in\beta X$, and so $\beta X\cap \Omega\not=\emptyset$. Hence, $\Omega$ is dense in $\beta \Omega$.
	
	Let us now show that $\beta \Omega$ is a Hausdorff space. Given $\omega_1,\omega_2\in\beta \Omega$ with $\omega_1\not=\omega_2$, there must be some subset $X\subseteq \Omega$ with $X\in\omega_1$  and $X\not\in \omega_2$. But then, $\Omega\setminus X\in\omega_2$. Hence, $\beta X,\beta(\Omega\setminus X)\subseteq \beta \Omega$ are open, disjoint subsets with $\omega_1\in\beta X$ and $\omega_2\in\beta(\Omega\setminus X)$. Therefore, $\beta \Omega$ is a Hausdorff space.
	
	Now, we need to prove that $\beta \Omega$ is compact. Let $\{\beta X_i\}_{i\in I}$ be a covering of $\beta \Omega$ by basic open sets. Suppose by contradiction that $$\bigcup_{i\in J}X_i\not= \beta \Omega$$ for every finite subset $J\subseteq I$. Then, this implies that $$\bigcap_{i\in J}(\Omega\setminus X_i)\not=\emptyset$$
	for every finite subset $J\subseteq I$, i.e. $\{\Omega\setminus X_i\}_{i\in I}$ has the finite intersection property. Hence, using Zorn's Lemma we can find an ultrafilter $\omega \in \beta \Omega$ such that $\Omega\setminus X_i\subseteq \omega$ for all $i\in I$.  Then, we have that $$\beta \Omega\setminus\bigg(\bigcup_{i\in I}\beta X_i\bigg)=\bigcap_{i\in I}\beta(\Omega\setminus X_i)\not=\emptyset,$$
	contradicting that $\{\beta X_i\}_{i\in I}$ is a covering of $\beta \Omega$. Therefore, we can extract from $\{\beta X_i\}_{i\in I}$ a finite subcovering, and so $\beta \Omega$ is a compact space.

	Finally, let $\Upsilon$ be a compact Hausdorff space and $f\colon \Omega\longrightarrow \Upsilon$ be a continuous map. Then, because $\Upsilon$ is both compact and Hausdorff, the map $f$ has a unique limit along every ultrafilter $\omega\in\beta \Omega$, and so we can define 
	\begin{displaymath}
	\beta f(\omega)=\lim_{x\to \omega} f(x)
	\end{displaymath}
	for  $\omega \in \beta \Omega$. We can then easily check that $\beta f\colon \beta\Omega\longrightarrow \Upsilon$ defined in this manner is the unique continuous extension of $f$ to $\beta \Omega$. 
\end{proof}

We can use the Stone-\u{C}ech compactification to prove the characterisation of amenability in terms of continuous actions fixing Borel measures. The key fact will be that the action of a group $G$ on itself can be extended to an action on $\beta G$ by using the universal property of the Stone-\u{C}ech compactification.

\begin{proposition}\label{prop:Action SC Compactification}
	Let $G$ be a group. Then, the action of $G$ on itself by right multiplication can be extended uniquely to an action of $G$ on $\beta G$ by homeomorphisms.
\end{proposition}
\begin{proof}
	Given $g\in G$, consider the right translation $\tau_g\colon G\longrightarrow G$ given by $\tau_g(h)=hg$. Then, the universal property of $\beta G$ given in Theorem~\ref{thm:SC Compactification Universal Property} implies that there is a unique continuous extension $\beta\tau_g\colon \beta G\longrightarrow \beta G$ of $\tau_g$ to $\beta G$ for each $g\in G$. Now, because $\tau_1=\id_G$ and the extension is unique, we have that $$\beta\tau_1=\id_{\beta G}.$$ Furthermore, given $g,h\in G$, using that $\tau_g\circ\tau_h=\tau_{gh}$ and that the extension is unique, we obtain that
	\begin{displaymath}
	\beta\tau_g\circ\beta\tau_h= \beta\tau_{gh}.
	\end{displaymath} 
	In particular, we have that $$\beta\tau_g\circ\beta\tau_{g^{-1}}=\beta\tau_{g^{-1}}\circ \beta\tau_{g}=\id_{\beta G}$$ for any $g\in G$, and so $\beta\tau_g$ is a homeomorphism of $\beta G$ for every $g\in G$.
	Therefore, the action of $G$ by right multiplication extends uniquely to an action of $G$ on $\beta G$ by homeomorphisms. 
\end{proof}

\section{Characterisations of Amenability}

We are now ready to come full circle and prove that all the different characterisations of amenability that we have discussed up to this point are actually equivalent.

\begin{theorem}\label{thm:Characterisations Amenablility}
	Let $G$ be a group. Then, the following are equivalent:
	\begin{enumerate}
		\item There is a right-invariant finitely additive probability measure on $G$.
		
		\item There is a right-invariant mean on $G$.
		
		\item The group $G$ satisfies the F\o lner condition.
		
		\item There is a F\o lner net in $G$.
		
		\item The group $G$ is non-paradoxical.
		
		\item Every affine continuous action of $G$ on a non-empty convex compact subset of a Hausdorff topological vector space has a fixed point.
		
		\item Every continuous action of $G$ on a non-empty compact Hausdorff topological space has an invariant Borel probability measure.
	\end{enumerate}
\end{theorem}

\begin{proof}
	We will prove the implications in the following diagram:
	\[\begin{tikzpicture}[
	> = implies, 
	shorten > = 1pt, 
	auto,
	node distance = 4cm, 
	scale=1, 
	transform shape, align=center, 
	state/.style={circle, draw, minimum size=0.8cm}]

	\node[state] (i) at (0,0) {i};
	\node[state] (ii) at (1.25,1.25) {ii};
	\node[state] (iii) at (3,0) {iii};
	\node[state] (iv) at (3.5,1.5) {iv};
	\node[state] (v) at (1.5,-1.25) {v};
	\node[state] (vi) at (1.5,3.) {vi};
	\node[state] (vii) at (-0.35,1.75) {vii};

	\draw[double distance=0.75mm,double,<->,bend left=8] (ii) to (i);
	\draw[double distance=0.75mm,double,->,bend left=15] (iii) to (i);
	\draw[double distance=0.75mm,double,->,bend right=8] (v) to (iii);
	\draw[double distance=0.75mm,double,->,bend right=8] (i) to (v);
	\draw[double distance=0.75mm,double,<->,bend right=10] (iii) to (iv);
	\draw[double distance=0.75mm,double,->,bend right=20] (iv) to (vi);
	\draw[double distance=0.75mm,double,->,bend right=15] (vi) to (vii);
	\draw[double distance=0.75mm,double,->,bend right=10] (vii) to (i);
	
	\end{tikzpicture}\]

	First, the fact that \textbf{(ii) is equivalent to (i)} is precisely Theorem~\ref{thm:ProbabilityMeasure<=>Mean}. 
	
	\smallskip
	
	Furthermore, the fact that \textbf{(iii) implies (i)} is a consequence of Theorem~\ref{thm:FolnerCondition=>ProbabilityMeasure}.

	\smallskip
	
	We also know that \textbf{(iii) is equivalent to (iv)} by Theorem~\ref{thm:Folner Condition <=> Folner Net}.
	
	\smallskip
	
	Let us see that \textbf{(v) implies (iii)}. We will actually show that $G$ not satisfying the F\o lner condition implies the existence of a paradoxical decomposition of $G$. Suppose that $G$ does not satisfy the F\o lner condition. Then, there exist a finite subset $X_0\subseteq G$ and $\varepsilon>0$ such that, for every finite non-empty subset $F\subseteq G$, there is some $g\in X_0$ satisfying that
	\begin{displaymath}
	\frac{\abs{F\setminus Fg}}{\abs{F}}>\varepsilon.
	\end{displaymath}
	Without loss of generality, we may assume that $1\in X_0$. Thus, for any finite non-empty subset $F\subseteq G$ we have that 
	\begin{align*}
	\abs{F}-\abs{FX_0}& =\abs{F\setminus FX_0}\\
	& \geq \abs{F\setminus Fg} \\
	& > \varepsilon\abs{F},
	\end{align*}
	and so we have 
	a finite subset $X_0\subseteq G$ and some $\lambda>1$ such that 
	\begin{displaymath}
	\abs{FX_0}\geq \lambda\abs{F}.
	\end{displaymath}
	Taking $n\in\N$ large enough that $\lambda^n\geq 2$ and writing $X=X_0^n$, we obtain a finite subset $X\subseteq G$ such that 
	\begin{displaymath}
	\abs{FX}\geq 2\abs{F}
	\end{displaymath}
	for every finite subset $F\subseteq G$.

	Let $\Omega$ be the collection of families  $$\{X_{(g,i)}\}_{(g,i)\in G\times \{1,2\}}$$ of finite subsets of $G$ satisfying the following conditions:
	\begin{itemize}
		\item For any finite subset $\Phi
		\subseteq G\times\{1,2\}$, 
		we have that
		\begin{displaymath}
		\abs[\bigg]{\bigcup_{(g,i)\in\Phi}X_{(g,i)}}\geq \abs{\Phi}.
		\end{displaymath}
		
		\item For every $(g,i)\in G\times\{1,2\}$, we have that
		\begin{displaymath}
		X_{(g,i)}\subseteq gX.
		\end{displaymath}
	\end{itemize}
	Note that $\Omega$ is non-empty, for $\{gX\}_{(g,i)\in G\times \{1,2\}}\in \Omega$. Indeed, any finite subset $\Phi\subseteq G\times\{0,1\}$ can be written as $$\Phi=\big(F_1\times\{1\}\big)\cup\big(F_2\times\{2\}\big)$$
	with $F_1,F_2\subseteq G$ finite, and so
	\begin{align*}
	\abs[\bigg]{\bigcup_{(g,i)\in\Phi}gX} & =\abs[\big]{(F_1\cup F_2)X} \\
	& \geq 2\abs{F_1\cup F_2} \\
	& \geq \abs{\Phi}.
	\end{align*}
	We can order $\Omega$ by component-wise inclusion. Then, every chain
	\begin{displaymath}
	\big\{X_{(g,i)}^{1}\big\}_{(g,i)\in G\times \{1,2\}}\geq \big\{X_{(g,i)}^{2}\big\}_{(g,i)\in G\times \{1,2\}} \geq \dotsb 
	\end{displaymath}
	has a lower bound, namely
	\begin{displaymath}
	\bigg\{\bigcap_{r\in\N} X_{(g,i)}^{r}\bigg\}_{(g,i)\in G\times \{1,2\}}.
	\end{displaymath}
	By Zorn's Lemma, $\Omega$ has a minimal element $\{M_{(g,i)}\}_{(g,i)\in G\times \{1,2\}}$. Let us see that $\abs{M_{(g,i)}}=1$ for all $(g,i)\in G\times \{1,2\}$. The construction of $\Omega$ implies that the $M_{(g,i)}$ are all non-empty. Assume by contradiction that $\abs{M_{(g_0,i_0)}}>1$ for some $(g_0,i_0)\in G\times \{1,2\}$, and take $g_1,g_2\in M_{(g_0,i_0)}$ distinct. For $l=1,2$, construct the family $\{M_{(g,i)}^{l}\}_{(g,i)\in G\times \{1,2\}}$ by replacing in $\{M_{(g,i)}\}_{(g,i)\in G\times \{1,2\}}$ the set $M_{(g_0,i_0)}$ with $M_{(g_0,i_0)}\setminus\{g_l\}$. By the minimality of $M_{(g,i)}$, neither of the families $\{M_{(g,i)}^{l}\}_{(g,i)\in G\times \{1,2\}}$ are in $\Omega$. Thus, there exist finite sets $\Phi_l\subseteq G\times\{1,2\}$ not containing $(g_0,i_0)$ such that 
	\begin{displaymath}
	\abs[\bigg]{M_{(g_0,i_0)}^{l}\cup\bigcup_{(g,i)\in\Phi_l}M_{(g,i)}^{l}}<\abs{\Phi_l}+1.
	\end{displaymath}
	Write
	\begin{displaymath}
	M^{l}=M_{(g_0,i_0)}^{l}\cup\bigcup_{(g,i)\in\Phi_l}M_{(g,i)}^{l}.
	\end{displaymath}
	Then,
	\begin{align*}
		\abs{\Phi_1}+\abs{\Phi_2} & \geq \abs{M^{1}}+\abs{M^{2}} \\
		& = \abs{M^{1}\cup M^{2}}+\abs{M^{1}\cap M^{2}} \\
		& = \abs[\Bigg]{M_{(g_0,i_0)}\cup \bigg(\bigcup_{(g,i)\in\Phi_1\cap\Phi_2}M_{(g,i)}\bigg)} \\
		& 
		\qquad
		+\abs[\Bigg]{\big(M_{(g_0,i_0)}\setminus\{g_1,g_2\}\big)\cup \bigg(\bigcup_{(g,i)\in\Phi_1\cap\Phi_2}M_{(g,i)}\bigg)} \\
		& \geq 1+\abs{\Phi_1\cup \Phi_2}+\abs{\Phi_1\cap \Phi_2} \\
		& =1+\abs{\Phi_1}+\abs{\Phi_1},
	\end{align*}
	a contradiction. This shows that $\abs{M_{(g,i)}}=1$ for all $(g,i)\in G\times \{1,2\}$. Also, the singletons $M_{(g,i)}$ must be pairwise disjoint by the properties of $\Omega$.
	
	
	Now, we define for each $x\in X$ the sets
	\begin{displaymath}
	Y_x=\{g\in G \mid gx\in M_{(g,1)}\}, \quad Z_x=\{g\in G \mid gx\in M_{(g,2)}\}.
	\end{displaymath}
	Write $M_{(g,i)}=\{h_{(g,i)}\}$. 
	Given $g\in G$, by the properties of $\Omega$ we have that $M_{(g,i)}\subseteq Xg$ for $i=1,2$, so there exists $x_i\in X$ such that $gx_i=h_{(g,i)}$, meaning that $g\in Y_{x_1}$ and $g\in Z_{x_2}$. Furthermore, if $g\in Y_x\cap Y_{x'}$ then $gx=gx'$, implying that $x=x'$ and the $Y_x$ are pairwise disjoint. The same is clearly true for the $Z_x$. Note also that all the $Y_x$ and the $Z_x$ are distinct due to the elements $h_{(g,i)}$ being distinct. Therefore, we can write
	\begin{displaymath}
	G=\bigcupdot_{x\in X}Y_x=\bigcupdot_{x\in X}Z_x.
	\end{displaymath}
	Finally, we have that 
	\begin{displaymath}
	Y_x x\cap Z_{x'}x'=Y_x x\cap Y_{x'}x'=Y_x x\cap Z_{x}x=\emptyset,
	\end{displaymath}
	for all distinct $x,x'\in X$, 
	and so $G$ is paradoxical by Proposition~\ref{prop:Paradoxical Equivalent Defs}.
	
	\smallskip

	Let us now show that \textbf{(i) implies (v)}. We will prove that if $G$ is paradoxical, then there cannot be any right-invariant finitely additive probability measure on $G$. Suppose that we have pairwise disjoint subsets $X_1,\dotsc,X_n$ and $Y_1,\dotsc, Y_m$ of $G$, and elements $g_1,\dotsc,g_n$ and $h_1,\dotsc,h_m$ in $G$ such that
	\begin{displaymath}
	G= \bigg(\bigcupdot_{i=1}^n X_i\bigg) \cupdot \bigg( \bigcupdot_{j=1}^m Y_j\bigg) = \bigcupdot_{i=1}^n X_ig_i = \bigcupdot_{i=1}^m Y_jh_j.
	\end{displaymath}
	Assume now by contradiction that there is a right-invariant finitely additive probability measure $\mu$ on $G$. On the one hand, we have that 
	\begin{displaymath}
	\sum_{i=1}^{n}\mu(X_i)=\mu\bigg(\bigcupdot_{i=1}^n X_ig_i\bigg) = \mu(G)=1,
	\end{displaymath}
	and analogously,
	\begin{displaymath}
	\sum_{j=1}^{m}\mu(Y_i)=1.
	\end{displaymath}
	On the other hand,
	\begin{align*}
	\mu(G) & =\mu\Bigg(\bigg(\bigcupdot_{i=1}^n X_i\bigg) \cupdot \bigg( \bigcupdot_{j=1}^m Y_j\bigg)\Bigg) \\
	& = \sum_{i=1}^{n}\mu(X_i) + \sum_{j=1}^{m}\mu(Y_i) \\
	& = 2,
	\end{align*}
	which contradicts the fact that $\mu(G)=1$. Therefore, no such a $\mu$ can exist on $G$.
	
	\smallskip
	
	Let us prove that \textbf{(iv) implies (vi)}. Let $X$ be  Hausdorff topological vector space and $C\subseteq X$ a non-empty convex compact subset. Assume that $G$ acts on $X$, and consequently on $C$, via an affine continuous action. Let $(F_i)_{i\in I}$ be a F\o lner net in $G$. Choose a point $x\in C$ and, for each $i\in I$, define
	\begin{displaymath}
	c_i=\frac{1}{\abs{F_i}}\sum_{h\in F_i}x h.
	\end{displaymath}
	Note that $c_i\in C$ because $C$ is convex. Without loss of generality, we may assume that the net $(c_i)_{i\in I}$ converges -- otherwise, since $C$ is compact, we can take a convergent subnet. Write then 
	\begin{displaymath}
	c=\lim_{i\in I}c_i.
	\end{displaymath}
	For any $g\in G$, we have that
	\begin{displaymath}
	c_i g= \frac{1}{\abs{F_i}}\sum_{h\in F_ig}xh,
	\end{displaymath}
	and so
	\begin{align*}
		c_ig-c_i & = \frac{1}{\abs{F_i}} \Bigg(\sum_{h\in F_ig\setminus F_i}xh - \sum_{h\in F_i\setminus F_ig}xh\Bigg) \\
		& = \frac{\abs{F_i\setminus F_ig}}{\abs{F_i}} \Bigg(\frac{1}{\abs{F_i\setminus F_ig}} \sum_{h\in F_ig\setminus F_i}xh - \frac{1}{\abs{F_i\setminus F_ig}} \sum_{h\in F_i\setminus F_ig}xh\Bigg)
	\end{align*}
	for every $i\in I$. Because $(F_i)_{i\in I}$ is a F\o lner net and $C$ is compact, we can conclude that
	\begin{displaymath}
	cg-c=\lim_{i\in I}(c_ig - c_i)=0.
	\end{displaymath}
	Therefore, $c$ is a fixed point for the action of $G$ on $C$.

	\smallskip
	
	Let us now see that \textbf{(vi) implies (vii)}. If $\Omega$ is a compact Hausdorff topological space, then the space $\mathfrak{M}(\Omega)$ of complex regular Borel measures on $\Omega$ can be identified by the Riesz Representation Theorem with the dual space of  $\mathcal{C}(\Omega)$, the space of continuous functions from $\Omega$ to $\C$. Denote by $\mathfrak{P}(\Omega)$ the set of Borel probability measures on $\Omega$. Then, we have that $\mathfrak{P}(\Omega)\subseteq \mathfrak{M}(\Omega)$ and $\mathfrak{P}(\Omega)$ is clearly convex. Furthermore, $\mathfrak{P}(\Omega)$ can be written as the zero set of a continuous map on $\mathfrak{M}(\Omega)$. Moreover, $\mathfrak{P}(\Omega)$ is contained in the unit ball of $\mathfrak{M}(\Omega)$, which is compact by the Banach-Alaoglu Theorem. Finally, the action of $G$ on $\Omega$ naturally induces an action on $\mathfrak{P}(\Omega)$, given by
	\begin{displaymath}
	\mu^g(X)=\mu(Xg^{-1})
	\end{displaymath}
	for $X\subseteq \Omega$ and $g\in G$. Therefore, by hypothesis we must have $\mathfrak{P}(\Omega)^G\not=\emptyset$.

	\smallskip
	
	Finally, let us see that \textbf{(vii) implies (i)}. By Proposition~\ref{prop:Action SC Compactification}, the action of $G$ on itself extends to a continuous action on its Stone-\u{C}ech compactification $\beta G$, which is a non-empty compact Hausdorff topological space. Then, by hypothesis there must be some $G$-invariant Borel probability measure $\mu$ on $\beta G$. Define then $\bar{\mu}\colon\mathcal{P}(G)\longrightarrow[0,1]$ by $$\bar{\mu}(X)=\mu(\beta{X})$$ 
	for $X\subseteq G$. It is clear that $$\bar{\mu}(G)=\mu(\beta G)=1.$$
	%
	Moreover, given any disjoint $X,Y\subseteq G$, we have that $\beta X,\beta Y\subseteq \beta G$ are disjoint as well, and so
	\begin{align*}
		\bar{\mu}(X\cupdot Y) & = \mu\big(\beta(X\cupdot Y)\big) \\
		& = \mu(\beta{X}\cupdot\beta{Y}) \\
		& =\mu(\beta{X})+\mu(\beta{Y}) \\ 
		& = \bar{\mu}(X)+\bar{\mu}(Y).
	\end{align*}
	%
	%
	Finally, given $X\subseteq G$ and $g\in G$, we have that $\beta(Xg)=(\beta X)g$ because $G$ acts by homeomorphisms, and so
	\begin{displaymath}
	\bar{\mu}(Xg)=\mu\big(\beta(Xg)\big)=\mu(\beta{X}g)=\mu(\beta{X})=\bar{\mu}(X).
	\end{displaymath}
	Therefore, $\bar \mu$ is a right-invariant finitely additive probability measure on $G$, and so $G$ is amenable.
\end{proof}

%% file: Chapter2.tex
\chapter{Quasitilings}
\label{chapter:Quasitilings}

In this chapter we will develop the theory of quasitilings for finitely generated amenable groups, and prove a result originally by D. S. Ornstein and B. Weiss stating that quasitilings always exist. Quasitilings prove to be a key tool in the proof of a number of results for finitely generated amenable groups. This chapter is primarily based on \cite{elek}, \cite{elek-szabo2} and \cite[\S 4.5]{kerr}.

\section{Cayley Graphs and Graph Approximations}

Before introducing the Cayley graph of a finitely generated group, let us fix some notation. A graph $X$ will consist of a set of vertices $V(X)$ and a set of edges $E(X)$. We will frequently identify $X$ with its set of vertices.

	Let $S$ be a finite set. An \emph{$S$-labelled graph} is a graph $X$ such that every directed edge $(x,y)\in E(X)$ is labelled by some $s\in S^{\pm 1}$, in such a way that $(y,x)$ is labelled by $s^{-1}\in S^{\pm 1}$ and for each $x\in X$ and $s\in S^{\pm 1}$ there is at most one edge from $x$ labelled by $s$. 

\begin{definition}
	Let $G$ be a group generated by a finite set $S$. 
	Then, the \emph{Cayley graph} of $G$ with respect to $S$, denoted by $\Cay(G,S)$, is the $S$-labelled graph with vertex set $G$ and directed edges $(g,gs)$ labelled by $s\in S^{\pm 1}$, with $g\in G$. 
\end{definition}

\begin{remark}
	Given a finite set $S$, consider the free group $F(S)$ on $S$. Then, an $S$-labelled graph is the same as an $F(S)$-set. Indeed, an action of $F(S)$ on any set automatically turns it into an $S$-labelled set, whereas the labels of an $S$-labelled graph give us an action of $F(S)$.
\end{remark}

We can consider F\o lner sequences in the Cayley graph of a finitely generated group, which turn out to be the same as  F\o lner sequences in the group itself.

\begin{definition}
	Let $G$ be a group generated by a finite set $S\subseteq G$. A sequence $(F_n)_{n\in\N}$ of finite subgraphs of $\Cay(G,S)$ is called a \emph{F\o lner sequence} in $\Cay(G,S)$ if for all $\varepsilon>0$ there exists some $N\in\N$ such that
	\begin{displaymath}
	\frac{\abs{\partial_SF_n}}{\abs{F_n}}< \varepsilon
	\end{displaymath}
	for all $n\geq N$.
\end{definition}

\begin{remark}
	It can be easily seen that this definition does not depend on the generating set, and so a F\o lner sequence in some Cayley graph of a group is a F\o lner sequnce in any Cayley graph of the group, and is in fact a F\o lner sequence in the group itself.
\end{remark}


We will now study approximations of graphs. 
Whenever we have a graph $X$, we can define a metric on $X$ by setting the distance between any two vertices of $X$ to be the shortest length of a path between them. Given $r>0$ and $x\in X$, we denote by $B_r(x)$ the ball of radius $r$ centred at $x$.

\begin{definition}
	Let $G$ be a group generated by a finite set $S\subseteq G$, $X$ be a finite $S$-labelled graph and $r\in\N$. We say that $X$ is an \emph{$r$-approximation} of $\Cay(G,S)$ if there exists some subgraph $X'\subseteq X$ such that
	\begin{displaymath}
	\abs{X'}> \Big(1-\frac1r\Big)\abs{X}
	\end{displaymath}
	and $B_r(x)$ is isomorphic to $B_r(1)$ as an $S$-labelled graph for every $x\in A$.
\end{definition}

This definition allows us to give another characterisation of F\o lner sequences.

\begin{proposition}
	Let $G$ be a group generated by a finite set $S\subseteq G$. A sequence $(F_n)_{n\in\N}$ of finite subgraphs of $\Cay(G,S)$ is a F\o lner sequence if and only if for every $r\in\N$ there exists some $N\in\N$ such that $F_n$ is an $r$-approximation of $\Cay(G,S)$ for all $n\geq N$.
\end{proposition}

Now, we proceed to define the notion of $r$-isomorphism of labelled graphs.

\begin{definition}
	Let $S$ be a finite set and $r\in\N$. Two $S$-labeled graphs $X_1$ and $X_2$ are said to be \emph{$r$-isomorphic} if there are subgraphs $X_i'\subseteq X_i$ such that
	\begin{displaymath}
	\abs[\big]{E(X_i')}\geq \Big(1-\frac1r\Big)\abs[\big]{E(X_i)}
	\end{displaymath}
	for $i=1,2$ and $X_1'$ is isomorphic to $X_2'$ as an $S$-labelled graph.
\end{definition}

\begin{lemma}\label{lem:Labelled Graphs r-Isomorphic}
	Let $X$, $Y$ and $Z$ be $S$-labelled graphs and $r\in\N$. If $X$ is $2r$-isomorphic to $Y$ and $Y$ is $2r$-isomorphic to $Z$, then $X$ is $r$-isomorphic to $Z$. 
	%
\end{lemma}

\begin{proof}
	Because $X$ and $Y$ are $2r$-isomorphic, there exist $X'\subseteq X$ and $Y'\subseteq Y$, and an isomorphism $\varphi\colon X'\longrightarrow Y'$, such that $$\abs[\big]{E(X')}\geq \Big(1-\frac{1}{2r}\Big)\abs[\big]{E(X)}, \quad \abs[\big]{E(Y')}\geq \Big(1-\frac{1}{2r}\Big)\abs[\big]{E(Y)}.$$ Similarly, there are $Y''\subseteq B$ and $Z'\subseteq Z$, and an isomorphism $\psi\colon Y''\longrightarrow Z'$, such that $$\abs{E(Y'')}\geq \Big(1-\frac{1}{2r}\Big)\abs[\big]{E(Y)}, \quad \abs[\big]{E(Z')}\geq \Big(1-\frac{1}{2r}\Big)\abs[\big]{E(Z)}.$$ 
	We can obtain $Y'\cap Y''$ from $Y'$ by erasing at most $\frac{1}{2r}\abs{E(Y)}$ edges from $Y'$. If we write $$X''=\varphi^{-1}(Y'\cap Y''),$$ we have that
	\begin{align*}
	\abs[\big]{E(X'')} & =\abs[\big]{E(Y'\cap Y'')} \\
	&\geq  \Big(1-\frac{1}{2r}\Big) \abs[\big]{E(Y')} \\
	& = \Big(1-\frac{1}{2r}\Big) \abs[\big]{E(X')}.
	\end{align*}
	Thus, we can obtain $X''$ from $X'$ by erasing at most $\frac{1}{2r}\abs{E(X')}$ edges from $X'$, and so  we can obtain $X''$ from $X$ by erasing at most $\frac{1}{r}\abs{E(X)}$ edges from $X$. Analogously, we can obtain $$Z''=\psi(Y'\cap Y'')$$ from $Z$ by erasing at most $\frac{1}{r}\abs{E(Z)}$ edges from $Z$. Hence, $\psi\circ\varphi\colon X''\longrightarrow Z''$ is an isomorphism, and $$\abs[\big]{E(X'')}\geq \Big(1-\frac{1}{r}\Big)\abs[\big]{E(X)},\quad \abs[\big]{E(Z'')}\geq \Big(1-\frac{1}{r}\Big)\abs[\big]{E(Z)},$$ meaning that $X$ is $r$-isomorphic to $Z$.
\end{proof}

\section{Quasitilings}

We will now prove a version of the Ornstein-Weiss Quasitiling Theorem for graphs presented in \cite{elek}. 
Before talking about quasitilings, we will need a number of auxiliary concepts about coverings of finite sets. 

\begin{definition}
	Let $F$ be a finite set, $(X_i)_{i\in I}$ a family of subsets of $F$, and $\lambda,\varepsilon\geq 0$.
	\begin{enumerate}
		\item We say that $(X_i)_{i\in I}$ is a \emph{$\lambda$-even covering of $F$ with multiplicity $M$} if
		\begin{displaymath}
		\sum_{i\in I}\chi_{X_i}\leq M,
		\end{displaymath}
		where $\chi_{X_i}$ is the characteristic function of $X_i$, 
		and
		\begin{displaymath}
		\sum_{i\in I}\abs{X_i}\geq \lambda M\abs{F}.
		\end{displaymath}
		
		\item We say that $(X_i)_{i\in I}$  \emph{$\lambda$-covers $F$} if
		\begin{displaymath}
		\abs[\Big]{\bigcup_{i\in I} X_i}\geq \lambda \abs{F}.
		\end{displaymath}
		
		\item  We say that $(X_i)_{i\in I}$ is \emph{$\varepsilon$-disjoint} if for each $i\in I$ there exists $Y_i\subseteq X_i$ such that 
		\begin{displaymath}
		\abs{Y_i}\geq (1-\varepsilon) \abs{X_i}
		\end{displaymath}
		and $(Y_i)_{i\in I}$ is a family of pairwise disjoint sets.
	\end{enumerate}
\end{definition}

\begin{lemma}\label{lem:Quasitilings1}
	Let $F$ be a finite set, $0<\lambda<1$ and $(X_i)_{i\in I}$ a $\lambda$-even covering of $F$. Then, for every 
	subset $Y\subseteq F$ there exists some $i\in I$ such that 
	\begin{displaymath}
	\frac{\abs{X_i\cap Y}}{\abs{X_i}}\leq \frac{\abs{Y}}{\lambda\abs{F}}.
	\end{displaymath}
\end{lemma}

\begin{proof}
	Suppose by contradiction that there is some $Y\subseteq F$ such that
	\begin{displaymath}
	\frac{\abs{X_i\cap Y}}{\abs{X_i}}> \frac{\abs{Y}}{\lambda\abs{F}}
	\end{displaymath}
	for all $i\in I$. Then, if the $\lambda$-even covering $(X_i)_{i\in I}$ has multiplicity $M$, we have that
	\begin{align*}
		\sum_{i\in I} \abs{X_i\cap Y} & > \frac{\abs{Y}}{\lambda\abs{F}}\sum_{i\in I}\abs{X_i} \\
		& \geq \abs{Y} M \\
		& \geq \sum_{y\in Y}\chi_{Y}(y)\sum_{i\in I}\chi_{X_i}(y) \\
		& = \sum_{i\in I}\sum_{y\in Y}\chi_{X_i\cap Y}(y) \\
		& = \sum_{i\in I}\abs{X_i\cap Y},
	\end{align*}
	a contradiction.
\end{proof}

\begin{lemma}\label{lem:Quasitilings2}
	Let $F$ be a finite set, $0\leq \varepsilon\leq \frac12$ and $0<\lambda\leq 1$. If $(X_i)_{i\in I}$ is a $\lambda$-even covering of $F$ by non-empty sets, then we can extract an $\varepsilon$-disjoint subcollection of $(X_i)_{i\in I}$ that $\varepsilon\lambda$-covers $F$.
\end{lemma}

\begin{proof}
	Let $\Omega$ be the collection of families $\{(X_i,Y_i)\}_{i\in I'}$ with $I'\subseteq I$ and $Y_i\subseteq X_i$ satisfying that $$\abs{Y_i}\geq(1-\varepsilon)\abs{X_i}$$ for every $i\in I'$ and the sets $Y_i$ are pairwise disjoint. We can order $\Omega$ by setting $$\big\{(X_i,Y_i)\big\}_{i\in I'} \preceq \big\{(X_i,Z_i)\big\}_{i\in I''}$$ if $I'\subseteq I''$ and $Y_i=Z_i$ for all $i\in I'$. It is clear that $\Omega$ is non-empty, for given any $i_0\in I$ we have that $\{(X_{i_0},X_{i_0})\}\in\Omega$. 
	Thus, 
	$\Omega$ has a maximal element, say $\{(X_i,Y_i)\}_{i\in J}$. Assume by contradiction that $(X_i)_{i\in J}$ does not $\varepsilon\lambda$-cover $F$, i.e. $$\abs[\Big]{\bigcup_{i\in J}X_i}<\varepsilon\lambda\abs{F}.$$ Then, Lemma~\ref{lem:Quasitilings1} implies that there exists some $i_0\in I$ such that 
	\begin{displaymath}
	\frac{\abs{X_{i_0}\cap \bigcup_{i\in J}X_i}}{\abs{X_{i_0}}} \leq \frac{\abs{\bigcup_{i\in J}X_i}}{\lambda\abs{F}} <\varepsilon.
	\end{displaymath}
	Thus, we can add the pair $(X_{i_0},X_{i_0}\setminus \bigcup_{i\in J} X_i)$ to the collection $\{(X_i,X_i)\}_{i\in J}$, contradicting its maximality.
\end{proof}

We are now ready to study quasitilings. 
We will introduce the version of quasitilings developed in \cite{elek} for graphs.  

Let $G$ be a group generated by a finite set $S\subseteq G$. Given a finite $S$-labelled graph $X$ and $r>0$, we denote by $Q_r(X)$ the set of vertices $x\in X$ such that the ball $B_r(x)\subseteq X$ is isomorphic to the ball $B_r(1)\subseteq \Cay(G,S)$ as an $S$-labelled graph. 

For each point $x\in Q_r(X)$, we have an isomorphism of $S$-labelled graphs $\phi_x\colon B_r(1)\longrightarrow B_r(s)$. Given $\varepsilon>0$, we will say that a collection $(T_1,\dotsc,T_m)$ of finite subsets of $B_{r/2}(1)\subseteq \Cay(G,S)$ is an \emph{$\varepsilon$-quasitiling} of $X$ if there exist $C_1,\dotsc, C_m\subseteq Q_r(X)$ such that the family
\begin{displaymath}
\bigcup_{k=1}^m\big\{\phi_x(T_k)\mid x\in C_k\big\}
\end{displaymath}
is $\varepsilon$-disjoint and $(1-\varepsilon)$-covers $X$.

Whenever our group is amenable, every finite graph which is a sufficiently good approximation of the Cayley graph of the group can be quasitiled by elements of a F\o lner sequence.

\begin{theorem}[{\cite[Theorem~2]{elek}}]\label{thm:Quasitilings Graphs}
	Let $G$ be a finitely generated amenable group with $S\subseteq G$ a finite generating set and  $(F_n)_{n\in \N}$ a F\o lner exhaustion of $G$. Given $\varepsilon>0$, there exist some $r>0$ and a finite subcollection $(T_1,\dotsc, T_m)$ of $(F_n)_{n\in \N}$ with $T_i\subseteq B_{r/2}(1)$ for $i=1,\dotsc, m$ such that every finite $S$-labelled graph $X$ satisfying that 
	\begin{displaymath}
	\frac{\abs[\big]{Q_r(X)}}{\abs{X}}>1-\frac{\varepsilon}{4}
	\end{displaymath}
	is $\varepsilon$-quasitiled by $(T_1,\dotsc, T_m)$.
\end{theorem}

\begin{proof}
	Take $m\in\N$ such that $(1-\frac\varepsilon2)^m<\varepsilon$. 
	Choose some $n_1\in \N$ and write $T_1=F_{n_1}$. Then, take $r_1\geq 1$ such that $T_{1}\subseteq B_{r_1/2}(1)$. Now, because $(F_n)_{n\in \N}$ is a F\o lner exhaustion, we can take $n_2\in\N$ such that $B_{r_1/2}\subseteq F_{n_2}$ and 
	\begin{displaymath}
	\frac{\abs[\big]{B_{r_1}(F_{n_2})\setminus F_{n_2}}}{\abs{F_{n_2}}}=\frac{\abs[\big]{F_{n_2}\cdot B_{r_1}(1)\setminus F_{n_2}}}{\abs{F_{n_2}}}<\frac{\varepsilon}{8}.
	\end{displaymath}
	Write $T_2=F_{n_2}$ and choose $r_2\geq 8 r_i$. Continuing in this manner, we can  extract from $(F_n)_{n\in \N}$ a subcollection $(T_1,\dotsc, T_m)$ such that
	\begin{displaymath}
	T_1\subseteq B_{r_1/2}(1)\subseteq T_2\subseteq\dotsb\subseteq T_m\subseteq B_{r_m/2},
	\end{displaymath}
	and 
	\begin{displaymath}
	\frac{\abs[\big]{B_{r_i}(T_{i+1})\setminus T_{i+1}}}{\abs{T_{i+1}}}< \frac{\varepsilon}{8},
	\end{displaymath}
	with $r_1\geq1$ and $r_i\leq \frac{r_{i+1}}{8}$ for $i=1,\dotsc,m-1$.
	
	Let $X$ be a finite $S$-labelled graph satisfying that 
	\begin{displaymath}
	\frac{\abs[\big]{Q_{r_m}(X)}}{\abs{X}}>1-\frac{\varepsilon}{4}.
	\end{displaymath}
	For each $x\in Q_{r_m}(X)$, we can consider the isomorphism of $S$-labelled graphs $\phi_x\colon B_{r_m}(1)\longrightarrow B_{r_m}(x)$. Note that $
	(\phi_x(T_m))_{x\in Q_{r_m}(X)}$ is a $\frac12$-even covering of $X$. Indeed, any $y\in \phi_x(T_m)$ is also in $Q_{r_m/2}$, and so $x\in \phi_y(T_m^{-1})$. Thus, every vertex of $X$ is covered by at most $\abs{T_m}$ tiles, i.e.
	\begin{displaymath}
	\sum_{x\in Q_{r_m}(X)} \chi_{\phi_x(T_m)} \leq \abs{T_m}.
	\end{displaymath} 
	Furthermore, 
	\begin{align*}
	\sum_{x\in Q_{r_m}(X)} \abs{\phi_x(T_m)} & = 
	\abs{Q_{r_m}(X)}\abs{T_m} \\
	& > \Big(1-\frac\varepsilon 4\Big)\abs{X}\abs{T_m} \\
	& >\frac12 \abs{T_m}\abs{X},
	\end{align*}
	thus showing that $
	(\phi_x(T_m))_{x\in Q_{r_m}(X)}$ is a $\frac12$-even covering of $X$ with multiplicity $\abs{T_m}$. Therefore, Lemma~\ref{lem:Quasitilings2} allows us to extract an $\varepsilon$-disjoint subcollection $(\phi_x(T_m))_{x\in C_m}$ with $C_m\subseteq Q_{r_m}(X)$ that $\frac{\varepsilon}{2}$-covers $X$.
	
	Assume now that for $1\leq k<m$ we have constructed sets $C_{k+1},\dotsc, C_m$ with $C_i\subseteq Q_{r_i}(X)$ such that  
	the collection $$\bigcup_{i=k+1}^m\big\{\phi_x(T_i)\mid x\in C_i\big\}$$ is $\varepsilon$-disjoint and $\lambda$-covers $X$, where
	\begin{displaymath}
	\lambda=\min\bigg\{1-\varepsilon,1-\Big(1-\frac{\varepsilon}{2}\Big)^{m-k+2}\bigg\}.
	\end{displaymath}
	Let $$X_K=X\setminus\bigg( \bigcup_{i=k+1}^m\bigcup_{x\in C_i} \phi_x(T_i)\bigg).$$ 
	If $\abs{X_k}<\varepsilon \abs{X}$,  we can simply take $C_1,\dotsc,C_k$ to be the empty set and we are finished. Assume then that $\abs{X_k}\geq \varepsilon \abs{X}$. 
	Because
	$$\bigcup_{i=k+1}^m\big\{\phi_x(T_i)\mid x\in C_i\big\}$$ is $\frac12$-disjoint and
	\begin{displaymath}
	\abs[\bigg]{\bigcup_{i=k+1}^m\bigcup_{x\in C_i} \phi_x(T_i)}\leq\abs{X}\leq \frac1\varepsilon \abs{X_k},
	\end{displaymath}
	we have that
	\begin{align*}
		\abs[\bigg]{\bigcup_{i=k+1}^m\bigcup_{x\in C_i} \Big(\phi_x\big(B_{r_k}(T_i)\big)\setminus\phi_x(T_i)\Big)} & \leq \frac{\varepsilon}{8}\sum_{i=k+1}^m\abs{T_i}\abs{C_i} \\
		& \leq \frac{\varepsilon}{4}\abs[\bigg]{\bigcup_{i=k+1}^m\bigcup_{x\in C_i} \phi_x(T_i)} \\
		& \leq\frac14\abs{X_k}.
	\end{align*}
	Observe that 
	\begin{displaymath}
	\frac{\abs[\big]{Q_{r_k}(X)}}{\abs{X}}>1-\frac{\varepsilon}{4}
	\end{displaymath}
	because $r_k\leq r_m$. 
	Consequently,
	\begin{align*}
		\abs[\big]{Q_{r_k}(X_k)} & = \abs[\bigg]{Q_{r_k}(X)\setminus \bigcup_{i=k+1}^m\bigcup_{x\in C_i} \Big(\phi_x(T_i) \cup  \phi_x\big(B_{r_k}(T_i)\big)\setminus\phi_x(T_i)\Big)} \\
		& = \abs[\big]{Q_{r_k}(X)} -  \abs[\bigg]{\bigcup_{i=k+1}^m\bigcup_{x\in C_i} \phi_a(T_i)} \\
		& \qquad -  \abs[\bigg]{\bigcup_{i=k+1}^m\bigcup_{x\in C_i}\Big(\phi_x\big(B_{r_k}(T_i)\big)\setminus\phi_x(T_i)\Big)} \\
		& > \Big(1-\frac{\varepsilon}{4}\Big) \abs{X} - \big(\abs{X}-\abs{X_k}\big) - \frac14\abs{X_k} \\
		& > \frac12\abs{X_k},
	\end{align*}
	and so 
	\begin{displaymath}
	\frac{\abs[\big]{Q_{r_k}(X_k)}}{\abs{X_k}}>1-\frac{1}{2}.
	\end{displaymath} 
	Therefore, the collection $
	(\phi_x(T_k))_{x\in Q_{r_k}(X_k)}$ is a $\frac12$-even covering of $X_k$, and so Lemma~\ref{lem:Quasitilings2} implies that we can extract an $\varepsilon$-disjoint subcollection $(\phi_x(T_k))_{x\in C_k}$ with $C_k\subseteq Q_{r_k}(X_k)$ that $\frac{\varepsilon}{2}$-covers $X_k$. Hence,
	$$\bigcup_{i=k}^m\big\{\phi_x(T_i)\mid x\in C_i\big\}$$
	is $\varepsilon$-disjoint and $(1-(1-\frac{\varepsilon}{2})^{m-k+1})$-covers $X$. Therefore, we can recursively construct sets $C_1,\dotsc, C_m$ with $C_i\subseteq Q_{r_i}(X)$ such that  $$\bigcup_{i=1}^m\big\{\phi_x(T_i)\mid x\in C_i\big\}$$
	is $\varepsilon$-disjoint and $(1-\varepsilon)$-covers $X$, i.e. the finite $S$-labelled graph $X$ can be $\varepsilon$-quasitiled by $(T_1,\dotsc, T_m)$.
\end{proof}

We will now present a classical application of quasitilings in the form of the Subadditive Function Theorem.

\begin{definition}
	Let $G$ be a group and $f$ 
	be a function from the set of finite subsets of $G$ to 
	$\R$
	. We say that $f(X)$ \emph{converges to $\lambda$ as $X$ becomes more and more invariant} if for every $\varepsilon>0$ there exist a finite subset  $F\subseteq G$ and $\delta>0$ such that $\abs{f(X)-\lambda}<\varepsilon$ for every non-empty $(F,\delta)$-invariant finite subset $X\subseteq G$.
\end{definition}

\begin{theorem}
	Let $G$ be a finitely generated amenable group with $S\subseteq G$ a finite generating set, and $\varphi$ be a function from the set of finite subsets of $G$ to $[0,\infty)$ satisfying the following conditions:
	\begin{enumerate}
		\item $\varphi(Xg)=\varphi(X)$ for every finite  subset $X\subseteq G$ and $g\in G$.
		
		\item $\varphi(X\cup Y)\leq \varphi(X)+\varphi(Y)$ for all finite $X,Y\subseteq G$.
	\end{enumerate}
	Then, $\frac{\varphi(X)}{\abs{X}}$ converges to a limit as $X$ becomes more and more invariant.
\end{theorem}

\begin{proof}
	Let $(F_n)_{n\in\N}$ be a nested F\o lner sequence in $G$. 
	First, observe that the subadditivity of $\varphi$ implies that $$\varphi(X)\leq \abs{X}\varphi\big(\{1\}\big)$$ for any finite subset $X\subseteq G$. Thus, the sequence $(\frac{\varphi(F_n)}{\abs{F_n}})_{n\in\N}$ is bounded in $\R$, and so we can define
	\begin{displaymath}
	\lambda=\liminf_{n\to\infty} \frac{\varphi(F_n)}{\abs{F_n}}.
	\end{displaymath}
	
	Let $\varepsilon>0$. Given any $0<\eta <\frac12$, we can extract from $(F_n)_{n\in\N}$ a finite subcollection $(T_1,\dotsc, T_m)$ that $\eta$-quasitile every $(T_m,\frac{\eta}{4})$-invariant finite subset of $G$, and such that 
	\begin{displaymath}
	\frac{\varphi(T_i)}{\abs{T_i}}\leq \lambda +\frac{\varepsilon}{2}
	\end{displaymath}
	for every $i=1,\dotsc,m$. Let $X\subseteq G$ be a $(T_m,\frac{\eta}{4})$-invariant finite subset of $G$. Then, there exist
	$C_1,\dotsc, C_m\subseteq G$ such that $$\bigcup_{i=1}^nT_i C_i\subseteq X$$ and the family $$\bigcup_{k=1}^n\{T_ic\mid c\in C_i\}$$ is $\varepsilon$-disjoint and $(1-\varepsilon)$-covers $X$. Then, for each $T_ic$ there exists a subset $\tilde T_i c$ with
	\begin{displaymath}
	\abs{\tilde T_i c}\geq (1-\varepsilon)\abs{T_i c}=(1-\varepsilon)\abs{T_i}
	\end{displaymath}
	and such that the $\tilde T_ic$ are pairwise disjoint. Then, we have that
	\begin{displaymath}
	\abs{X}\geq \sum_{i=1}^n\sum_{c\in C_i}\abs{\tilde T_i c}\geq (1-\varepsilon)\sum_{i=1}^n\sum_{c\in C_i}\abs{\tilde T_i},
	\end{displaymath}
	and so
	\begin{align*}
		\varphi(X) & \leq \varphi\bigg(\bigcup_{i=1}^n\bigcup_{c\in C_i}T_ic\bigg) + \varphi\bigg(X\setminus \bigcup_{i=1}^n\bigcup_{c\in C_i}T_ic\bigg) \\
		& \leq \sum_{i=1}^n\sum_{c\in C_i}\varphi(T_i) + \varepsilon\abs{X}\varphi\big(\{1\}\big) \\
		& \leq \Big(\lambda +\frac{\varepsilon}{2}\Big) \sum_{i=1}^n\sum_{c\in C_i}\abs{T_i} + \varepsilon\abs{X}\varphi\big(\{1\}\big) \\
		& \leq \Big(\lambda +\frac{\varepsilon}{2}\Big)\frac{\abs{X}}{1-\eta} + \varepsilon\abs{X}\varphi\big(\{1\}\big).
	\end{align*}
	Hence,
	\begin{displaymath}
	\frac{\varphi(X)}{\abs{X}} \leq \frac{1}{1-\eta} \Big(\lambda +\frac{\varepsilon}{2}\Big) + \eta \varphi\big(\{1\}\big),
	\end{displaymath}
	and we can take $\eta$ small enough that 
	\begin{displaymath}
	\frac{\varphi(X)}{\abs{X}} < \lambda +\varepsilon
	\end{displaymath}
	for every $(T_m,\frac{\eta}{4})$-invariant subset $X\subseteq G$. Therefore, $\frac{\varphi(X)}{\abs{X}}$ converges to $\lambda$ as $X$ becomes more and more invariant.
\end{proof}

\section{Approximations by Linear Combinations}

In this section, we will introduce the concept of linear combination of a sequence of graphs. This will allow us to reformulate Theorem~\ref{thm:Quasitilings Graphs} and then give a stronger version of this same result that we will need in the following chapter.

\begin{definition}
	Let $\mathcal{T}=(T_1,\dotsc,T_m)$ be a finite sequence of $S$-labelled graphs and $\alpha =(\alpha_1,\dotsc,\alpha_m)\in\N_0^m$. We define the \emph{linear combination of $\mathcal{T}$ with coefficient vector $\alpha$}, denoted by $\alpha\mathcal{T}$, as the disjoint union of $\alpha_i$ copies of $T_i$ for each $i=1,\dotsc,m$.
\end{definition}

\begin{lemma}\label{lem:Combination Graphs}
	Let $\mathcal{T}=(T_1,\dotsc,T_m)$ be a finite sequence of $S$-labelled graphs, each of them with at least one edge. Then, given $r\in\N$ there exists some $M\in\N$ such that, for all $\alpha,\beta\in \N_0^m\setminus\{0\}$ satisfying that $\norm{\beta}\geq M\norm{\alpha}$ and $$\norm[\bigg]{\frac{\alpha}{\norm{\alpha}}-\frac{\beta}{\norm{\beta}}}\leq \frac{1}{M},$$
	we have that $\beta\mathcal{T}$ is $r$-isomorphic to $t\alpha\cdot\mathcal{T}$ for some $t\in\N$.
\end{lemma}
\begin{proof}
	Let $$E=\max_{1\leq i\leq m}\abs[\big]{E(T_i)}$$ and $r>0$. Take $M\geq 2mEr$ and $\alpha,\beta \in \N^m\setminus\{0\}$ such that that $\norm{\beta}\geq M\norm{\alpha}$ and $$\norm[\bigg]{\frac{\alpha}{\norm{\alpha}}-\frac{\beta}{\norm{\beta}}}\leq \frac{1}{M}.$$
	Now, consider $t\in\N$ the largest integer such that $\norm{t\alpha}\leq\norm{\beta}$. Note that we must have $M\leq t$ and $\norm{\alpha}\geq 1$. Consequently,
	\begin{align*}
		\norm{\beta-t\alpha} & \leq \norm{\alpha}+t\norm[\bigg]{\frac{\alpha}{\norm{\alpha}}-\frac{\beta}{\norm{\beta}}} \\
		& \leq \norm{\alpha} +\frac{t}{M} \\
		& \leq  \frac{2t}{M}\norm{\alpha},
	\end{align*}
	and so each coordinate of $\beta$ differs from the corresponding coordinate of $t\alpha$ by at most $\frac{2t}{M}\norm{\alpha}$. Therefore, $\beta\mathcal{T}$ can be obtained from $t\alpha\mathcal{T}$ by adding or deleting at most $\frac{2t}{M}\norm{\alpha}$ copies of each $T_i$, and the number of edges either added or deleted is at most $\frac{2mEt}{M}\norm{\alpha}$.
	
	Furthermore, the graph $t\alpha\mathcal{T}$ has at least $t\norm{\alpha}$ edges, and $\beta\mathcal{T}$ has at least $\norm{\beta}\geq t\norm{\alpha}$ edges. Thus, because $M\geq 2mEr$, we obtain that $\beta\mathcal{T}$ is $r$-isomorphic to $t\alpha\mathcal{T}$.
\end{proof}

The following is a reformulation of Theorem~\ref{thm:Quasitilings Graphs} in terms of graph approximations and $r$-isomorphisms of labelled graphs.

\begin{theorem}[{\cite[Proposition~2.7]{elek-szabo}}]\label{thm:Quasitilings Graphs Reformulation}
	Let $G$ be a finitely generated amenable group with $S\subseteq G$ a finite generating set and $(F_n)_{n\in \N}$ a F\o lner exhaustion in $G$. Then, given $r\in\N$ there exists some $R\in\N$ and a finite subcollection $(T_1,\dotsc, T_m)$ of $(F_n)_{n\in \N}$ such that every $R$-approximation of $\Cay(G,S)$ is $r$-isomorphic to some linear combination of $(T_1,\dotsc, T_m)$.
\end{theorem}

This result can be refined in order to obtain a stronger version of the theorem that gives us quasitings of a particular type, which will prove to play a key role in the proof of Theorem~\ref{thm:Sofic Embeddings Conjugate}. 

\begin{theorem}[{\cite[Proposition~2.8]{elek-szabo}}]\label{thm:Quastitilings Graph 2}
	Let $G$ be a finitely generated amenable group with $S\subseteq G$ a finite generating set and $(F_n)_{n\in \N}$ a F\o lner exhaustion in $G$. Then, given $r\in\N$ there exists some $R\in\N$, a finite subsequence $\mathcal{T}=(T_1,\dotsc, T_m)$ of $(F_n)_{n\in \N}$ and some $\alpha\in\N_0^m\setminus\{0\}$ such that every $T_i$ is a $2r$-approximation of $\Cay(G,S)$, and every $R$-approximation of $\Cay(G,S)$ is $2r$-isomorphic to $t\alpha\cdot\mathcal{T}$ for some $t\in\N$.
\end{theorem}

\begin{proof}
	Let $r\in\N$. Assume without loss of generality that every set in $(F_n)_{n\in \N}$ is a $2r$-approximation of $\Cay(G,S)$. Using Theorem~\ref{thm:Quasitilings Graphs Reformulation} we can obtain a subsequence $\mathcal{T}=(T_1,\dotsc, T_m)$ of $(F_n)\ninN$ and some $R_0\in\N$ such that every $R_0$-approximation of $\Cay(G,S)$ is $8r$-isomorphic to some linear combination of $\mathcal{T}$. Now, take $(F_n')\ninN$ a subsequence of $(F_n)\ninN$ such that every $F_n'$ is an $R_0$-approximation of $\Cay(G,S)$. 
	
	For each $n\in\N$, take $\beta_n\in\N_0^m\setminus\{0\}$ such that $F_n$ is $8r$-isomorphic to $\beta_n\cdot\mathcal{T}$. Then, $(\frac{\beta_n}{\norm{\beta_n}})_{n\in\N}$ is a sequence of unit vectors in $\R^m$, so it must have an accumulation point $v\in\R^m$. Applying Lemma~\ref{lem:Combination Graphs}, we get that there exists some $M\in\N$ such that for all $\alpha,\beta\in \N_0^m\setminus\{0\}$ satisfying that $\norm{\beta}\geq M\norm{\alpha}$ and $$\norm[\bigg]{\frac{\alpha}{\norm{\alpha}}-\frac{\beta}{\norm{\beta}}}\leq \frac{1}{M}$$
	we have that $\beta\mathcal{T}$ is $8r$-isomorphic to some integer multiple of $\alpha\mathcal{T}$. Take some $N\in\N$ such that 
	$$\norm[\bigg]{\frac{\beta_{N}}{\norm{\beta_{N}}}-v}< \frac{1}{2M},$$
	and extract from $(F_n')\ninN$ a subsequence $(F_{n_k}')_{k\in \N}$ such that 
	$$\norm[\bigg]{\frac{\beta_{n_k}}{\norm{\beta_{n_k}}}-v}< \frac{1}{2M}$$
	and $\norm{\beta_{n_k}}\geq M$ for all $k\in\N$. The sequence $(F_{n_k}')_{k\in\N}$ is still a F\o lner exhaustion, and it satisfies that
	$$\norm[\bigg]{\frac{\beta_N}{\norm{\beta_N}}-\frac{\beta_{n_k}}{\norm{\beta_{n_k}}}}< \frac{1}{M}$$ for all $k\in\N$. 
	Thus, we have that $\beta_{n_k}\cdot \mathcal{T}$ is $8r$-isomorphic to $t_0\beta_N\cdot\mathcal{T}$ for some $t_0\in\N$.
	
	Applying Theorem~\ref{thm:Quasitilings Graphs} with the sequence $(F_{n_k}')_{k\in\N}$, we obtain a finite subfamily $\mathcal{Q}=(Q_1,\dotsc, Q_l)$ of $(F_{n_k}')_{k\in\N}$ and some $R_1\in\N$ such that every $R_1$-approximation of $\Cay(G,S)$ is $8r$-isomorphic to a linear combination of $\mathcal{Q}$. Let $X$ be an $R_1$-approximation of $\Cay(G,S)$. Then, every $Q_i$ is $8r$-isomorphic to the corresponding $\beta_{n_k}\cdot\mathcal{T}$, which is in turn $8r$-isomorphic to $t_0\beta_N\cdot\mathcal{T}$, and so Lemma~\ref{lem:Labelled Graphs r-Isomorphic} yields the result that $X$ is $2r$-isomorphic to $t\beta_N\cdot\mathcal{T}$ for some $t\in\Z$.
\end{proof}

%% file: Chapter3.tex
\chapter{Sofic Groups}
\label{chapter:SoficGroups}

In this chapter we will review the concept of sofic groups. Before fully focusing on them, we will discuss the family of residually finite groups. After that, we will introduce sofic groups as groups whose Cayley graphs can be approximated by finite graphs. In the last section, we will present a characterisation of soficity in terms of ultraproducts of finite symmetric groups. The main references for this chapter are \cite[\S 4, \S7]{ceccherini}, \cite{elek-szabo}, \cite{elek-szabo2} and \cite{szoke}.

\section{Residually Finite Groups}

We will now briefly introduce the concept of residually finite groups, which are groups in which elements can be distinguished by taking finite quotients. They serve as a generalisation of finite groups. 


\begin{definition}
	A group $G$ is said to be residually finite if for any $g\in G$ with $g\not=1$ there is some normal subgroup $N\normaleq G$ of finite index such that $g\not\in N$.
\end{definition}

\begin{proposition}\label{prop: Res Finite Charact}
	Let $G$ be a group. Then, $G$ is residually finite if and only if for every finite subset $F\subseteq G$ there exist a finite group $H$ and a homomorphism $\varphi\colon G\longrightarrow H$ which is injective in $F$.
\end{proposition}

\begin{proof}
	Assume that $G$ is residually finite and let $F\subseteq G$ be finite. Then, for every $g,h\in F$ with $g\not=h$ there exists some normal subgroup $N_{gh}\normaleq G$ of finite index such that $gh^{-1}\not\in N_{gh}$. Hence, the natural projection 
	\begin{displaymath}
	\varphi\colon G\longrightarrow \prod_{\substack{g,h\in F \\ g\not=h}}G/N_{gh}
	\end{displaymath}
	is a homomorphism from $G$ to a finite group which is injective in $F$. 
	
	Let us now prove the converse. Assume that for any finite subset of $G$ there is a homomorphism as required. Given $g\in G$ with $g\not =1$, there exist a finite group $H$ and a homomorphism $\varphi\colon G\longrightarrow H$ such that $\varphi(g)\not= 1$. Then, $\ker\varphi\normaleq G$ has finite index and $g\not\in \ker \varphi$.
\end{proof}

\begin{examples}
	\begin{enumerate}
		\item Finite groups are residually finite.

		\item Infinite simple groups, as well as groups with no non-trivial finite quotients, are never residually finite.

		\item The direct product of residually finite groups is once again residually finite. Indeed, let $$G=\prod_{i\in I}G_i$$ with $G_i$ a residually finite group for each $i\in I$. Then, given any $g=(g_i)_{i\in I}\in G$ with $g\not=1$, there is some $i_0\in I$ such that $g_{i_0}\not= 1$. Because $G_{i_0}$ is residually finite, we can find some normal subgroup $N_{i_0}\normaleq G_{i_0}$ of finite index in $G_{i_0}$ such that $g_{i_0}\not\in N_{i_0}$. Set now $N_i=G_i$ for each $i\not=i_0$, and $$N=\prod_{i\in I}N_i.$$
		It is then clear that $N\normaleq G$ has finite index and $g\not\in N$.
		
		\item Every finitely generated abelian group is residually finite. Due to the fact that finitely generated groups are direct products of finite groups and copies of $\Z$, it follows from the previous examples that we only have to show that  the group $\Z$ is residually finite. Indeed, given any $n\in\Z$ with $n\not= 1$, take any $m\in\Z$ such that $m\nmid n$. Then, we have that $m\Z\normaleq \Z$ is of finite index and $n\not\in m\Z$.

		\item The group $\Q$ is not residually finite. Given any finite group $H$ with $\abs{H}=n$, every homomorphism $\varphi\colon \Q\longrightarrow H$ is trivial, for we must have that  $$\varphi(x)=\varphi\bigg(\frac{x}{n}\bigg)^n=1$$ for every $x\in\Q$. This same argument implies that the additive group of a field of characteristic zero is never residually finite.
		
		\item Free groups are residually finite. Let $F(S)$ be the free group on the set $S$. Take $w\in F(S)$ with $w\not=1$, and write $w=s_1\dotsm s_n$ in reduced form with $s_i\in S^{\pm 1}$. We can now define a map $f\colon S\longrightarrow S_{n+1}$. Let $s\in S$.  If $s\not\in \{a_1^{\pm 1},\dotsc, a_n^{\pm 1}\}$ we set $f(s)=\id$. Otherwise, consider the sets $$X_s=\{i\mid s_i=s\},\quad Y_s=\{i\mid s_i^{-1}=s\},$$ and set $f(s)$ to be some $\sigma\in S_{n+1}$ such that $i\cdot \sigma=i+1$ for $i\in X_s$ and $(i+1)\cdot \sigma=i$ for $i\in Y_s$. Note that this is well-defined due to $w$ being in reduced form. Making use of the universal property of $F(S)$, we can extend the map $f\colon S\longrightarrow S_{n+1}$ to a homomorphism $\bar f\colon F(S)\longrightarrow S_{n+1}$. Furthermore, we can see that $\bar f(w)$ sends $1$ to $n+1$, and so $\bar f(w)\not=\id$. Therefore, $\ker \bar f\normaleq F(S)$ has finite index and $w\not\in\ker \bar f$.
	\end{enumerate}
\end{examples}

\section{Sofic Groups}
\label{section:Sofic Groups}

Sofic groups serve as a generalisation of both amenable and residually finite groups. They were first introduced in \cite{gromov} by M. Gromov in 1999. Their name, given to them in 2000 by Weiss in \cite{weiss}, comes from the Hebrew word for finite.

Let $G$ be a finitely generated group with $S\subseteq G$ a finite generating set of $G$, and consider the free group $F(S)$. Then, we can write $G=F(S)/N$ for some normal subgroup $N\unlhd F(S)$. Now, let $(X_n)\ninN$ be a sequence of finite $F(S)$-sets, and write 
\begin{displaymath}
P_r(X_n)=\big\{x\in X_n\mid \text{if } w\in B_r(1_{F(S)}), \text{ then } xw=x \text{ if and only if } w\in N\big\}
\end{displaymath}
for each $n,r\in\N$. 

\begin{definition}
	Let $G$ be a finitely generated group with $S\subseteq G$ a finite generating set. We say that a sequence $(X_n)_{n\in\N}$ of finite $F(S)$-sets is a \emph{sofic approximation} of $G$ if 
	\begin{displaymath}
	\lim_{n\to \infty} \frac{\abs[\big]{P_r(X_n)}}{\abs{X_n}} =1
	\end{displaymath}
	for all $r\in \N$. 
\end{definition}

This concept of sofic approximation leads us to the definition of sofic groups.
	
\begin{definition}
	A finitely generated group with a sofic approximation is called a \emph{sofic group}. In general, a group is said to be sofic if all of its finitely generated subgroups are sofic.
\end{definition}


This definition of soficity can be interpreted in a more geometric manner. Let $G$ be a finitely generated group with a finite generating set $S\subseteq G$ and  $(X_n)_{n\in \N}$ a sequence of finite $F(S)$-sets. 
The action of $F(S)$ on $X_n$ turns it into an $S$-labelled graph for each $n\in\N$. Then, for every $n,r\in\N$ we have that 
\begin{displaymath}
P_r(X_n) \subseteq Q_r(X_n) \subseteq P_{2r}(X_n),
\end{displaymath}
where $Q_r(X_n)$ is the set of all $x\in X_n$ such that $B_r(x)\subseteq X_n$ is isomorphic to $B_r(1)\subseteq \Cay(G,S)$ as an $S$-labelled graph. It  follows from this that $(X_n)_{n\in\N}$ is a sofic approximation of $G$ if and only if 
\begin{displaymath}
\lim_{n\to \infty} \frac{\abs[\big]{Q_r(X_n)}}{\abs{X_n}} = 1
\end{displaymath}
for all $r\in \N$. This last condition means that for every $r\in\N$ there exists some $N\in\N$ such that $X_n$ is an $r$-approximation of $\Cay(G,S)$ for all $n\geq N$.

We will now introduce another characterisation of soficity that looks quite similar to the characterisation of residually finite groups in Proposition~\ref{prop: Res Finite Charact}.

Given $n\in\N$, we can consider the symmetric group $S_n$, on which we can define the \emph{normalised Hamming distance} $d_n\colon S_n\times S_n\longrightarrow [0,1]$ by
\begin{displaymath}
d_n(\sigma,\tau)=\frac{\abs[\big]{\{i\mid \sigma(i)\not=\tau(i)\}}}{n}
\end{displaymath}
for $\sigma,\tau\in S_n$. It is not difficult to see that $d_n$ defines a bi-invariant metric on the group $S_n$, i.e. a metric such that 
\begin{displaymath}
d_n(\gamma\sigma,\gamma\tau)=d_n(\sigma\gamma,\tau\gamma)=d_n(\sigma,\tau)
\end{displaymath}
for every $\sigma,\tau,\gamma\in S_n$.

\begin{proposition}\label{prop:Sofic Group Charac}
	Let $G$ be a finitely generated group. Then, $G$ is sofic if and only if for every finite subset $F\subseteq G$ and $\varepsilon>0$ there exist some $n\in\N$ and a map $\varphi\colon G\longrightarrow S_n$ satisfying the following conditions:
	\begin{enumerate}
		\item For every $g,h\in F$, we have that 
		\begin{displaymath}
		d_n\big(\varphi(gh),\varphi(g)\phi(h)\big)<\varepsilon.
		\end{displaymath}

		\item We have that
		\begin{displaymath}
		d_n\big(\varphi(1),\id_n\big)<\varepsilon.
		\end{displaymath}

		\item For every $g\in F\setminus\{1\}$, we have that 
		\begin{displaymath}
		d_n\big(\varphi(g),\id_n\big)\geq 1-\varepsilon.
		\end{displaymath}
	\end{enumerate}
\end{proposition}
\begin{proof}
	Assume first that $G$ is sofic, with $S\subseteq G$ a finite generating set of $G$. Let $\varepsilon>0$ and $F\subseteq G$ a finite subset. We choose $r\in \N $ such that $F^2\subseteq B_r(1_G)$. Let $X$ be an $S$-labelled graph with a subgraph $X'\subseteq X$ such that $$\abs{X'}\geq (1-\varepsilon)\abs{X}$$ and for every $x\in X'$ there is an  isomorphism $\psi_x\colon B_r(1_G)\longrightarrow B_r(x)$ of $S$-labelled graphs.  
	We can then define a map $\varphi\colon G\longrightarrow S(F)$ by setting $x\cdot\varphi (g)= \psi_x(g)$ if $g\in B_r(1_G)$, and choosing $\varphi(g)$ arbitrarily otherwise. We can easily check that conditions (i), (ii) and (iii) are satisfied.


	Let us now prove the converse. Given any $r\in\N$, take $F=B_{2r+2}(1_G)$ and $\varepsilon>0$. Let $\varphi\colon G\longrightarrow S_n$ be a map satisfying conditions (i), (ii) and (iii) for our chosen $F$ and $\varepsilon$, and write $X=\{1,\dotsc,n\}$. For each $x\in X$, define $\psi_x\colon B_{r+1}(1_G)\longrightarrow X$ by setting $\psi_x(g)=x\cdot\varphi (g)$ for $g\in B_{r+1}(1_G)$. Let $X'\subseteq X$ be the set of points $x\in X$ satisfying the following conditions: 
	\begin{enumerate}
		\item[(a)] $\psi_x(gs)=\psi_{\psi_x(g)}(s)$ for all $g\in B_r(1_G)$ and $s\in S$.
		
		\item[(b)] $\psi_x(g)\not=\psi_x(h)$ for any $g,h\in B_{r+1}(1_G)$.
	\end{enumerate} 
	Now, the directed edges of $X$ are $(x,\psi_x(s))$ labelled by $s$, with $x\in X$ and $s\in S$. It can be seen that $B_r(x)=\psi_x(B_r(1_G))$ and all the edges coming out of it are in $B_{r+1}(x)=\psi_x(B_{r+1}(1_G))$, for each $x\in X'$. Conditions (a) and (b) imply that $\psi_x$ preserves the edges coming out of $\psi_x(g)$ and that $\psi_x$ is injective for $x\in X'$ and $g\in B_{r}(1_G)$, respectively. Therefore, we have that $\psi_x\colon B_r(1_G)\longrightarrow B_r(x)$ is an isomorphism of $S$-labelled graphs for every $x\in X'$.
	
	Moreover, condition (a) gives us  $\abs{B_r(1_G)}\abs{S}$ equations to check. But condition (i) says that each of these can only fail on a subset of $X$ of size at most $\varepsilon\abs{X}$.  Furthermore, condition (c) gives us $\abs{B_{r+1}}^2$ inequalities. Given $g,h\in B_{r+1}(1_G)$ and $x\in X$, applying conditions (i), (ii) and (ii) we get that 
	\begin{align*}
	\psi_x(g)\cdot \varphi(g^{-1}) & = x\cdot \varphi(g^{-1}g) \\
	& = x \\
	& \not= x\cdot \varphi(g^{-1}h) \\
	& =\psi_x(h)\cdot \varphi(g^{-1}),
	\end{align*} 
	and so $\psi_x(g)\not=\psi_x(h)$ for every $x\in X$ save for those in a set of size at most $4\varepsilon\abs{X}$. Therefore, if we take
	\begin{displaymath}
	\varepsilon<\frac{1}{r\big(4\abs{B_{r+1}(1_G)}^2+\abs{B_r(1_G)}\abs{S}\big)},
	\end{displaymath}
	we have that 
	\begin{displaymath}
	\abs{X'}\geq (1-\varepsilon)\abs{X},
	\end{displaymath}
	and so $G$ is sofic.
\end{proof}

We are now ready to show that amenable groups and residually finite groups are sofic. These classes of groups constitute our main two examples of sofic groups.

\begin{examples}
	\begin{enumerate}
		\item Amenable groups are sofic. Let $G$ be a finitely generated amenable group with $S\subseteq G$ a finite generating set and $(F_n)\ninN$ a F\o lner sequence in $G$. Given $r\in\N$, we have that
		\begin{displaymath}
		Q_r(F_n)= \big\{x\in F_n\mid B_r(1_G)x\subseteq F_n\big\}
		\end{displaymath}
		for all $n\in\N$. By Proposition~\ref{prop:FolnerSeq Equivalent Defs}, for any $\varepsilon>0$ there exists some $N\in\N$ such that 
		\begin{displaymath}
		\abs[\big]{Q_r(F_n)}\geq (1-\varepsilon)\abs{F_n}
		\end{displaymath}
		for all $n\geq N$, implying that 
		\begin{displaymath}
		\lim_{n\to \infty} \frac{\abs{Q_r(F_n)}}{\abs{F_n}} = 1
		\end{displaymath}
		for every $r\in\N$. 
		Therefore, $(F_n)_{n\in\N}$ is a sofic approximation of $G$. 
		
		\item Residually finite groups are sofic. This is a direct consequence of Proposition~\ref{prop: Res Finite Charact} and Proposition~\ref{prop:Sofic Group Charac}, along with Cayley's Theorem stating that every finite group is a subgroup of some finite symmetric group.
	\end{enumerate}
\end{examples}

\begin{remarks}
	\begin{enumerate}
		\item The class of sofic groups is closed under taking subgroups, direct products, inverse limits, direct limits, free products and amenable extensions, as shown in \cite{elek-szabo2}.
		
		\item One of the big open problems on the topic of sofic groups is the question of whether every group is sofic. There are no currently known non-sofic groups.
	\end{enumerate}
\end{remarks}

\section{Ultraproducts}

There is another characterisation of soficity via ultraproducts of finite symmetric groups. Our main interest in this characterisation is that it allows us to state Theorem~\ref{thm:Sofic Embeddings Conjugate}, which gives us a characterisation from \cite{elek-szabo} of amenable groups as those sofic groups whose sofic approximations are conjugate.

If 
we fix a non-principal ultrafilter $\omega$ on $\N$ and take a sequence $(k_n)_{n\in\N}$ in $\N$, we can define  a map $$d_{\omega}\colon \prod_{n\in\N}S_{k_n}\times \prod_{n\in\N}S_{k_n}\longrightarrow [0,1]$$ given by
\begin{displaymath}
d_{\omega}(\sigma,\tau)=\lim_{n\to \omega} d_{k_n}(\sigma_{k_n},\tau_{k_n})
\end{displaymath}
for elements $$\sigma=(\sigma_{k_n})_{n\in\N},\tau=(\tau_{k_n})_{n\in\N}\in\prod_{n\in\N}S_{k_n}.$$ We can see that 
\begin{displaymath}
N_{\omega}=\bigg\{\sigma\in \prod_{n\in\N}S_{k_n} \mid d_{\omega}(\id,\sigma)=0\bigg\}
\end{displaymath}
is a normal subgroup of $\prod_{n\in\N}S_{k_n}$, and so we can consider the quotient $$\Sigma_{\omega}=\bigg(\prod_{n\in\N}S_{k_n}\bigg)/N_{\omega}.$$
We say that $\Sigma_{\omega}$ is the \emph{ultraproduct} of $(S_{k_n})_{n\in\N}$ with respect to $\omega$.

We can then show that a finitely generated group is sofic if and only if it can be embedded in a faithful manner into such an ultraproduct of finite symmetric groups. We will only give an outline of the proof of this characterisation. 

Let $G$ be a finitely generated group with a finite generating set $S\subseteq G$ and a sofic approximation $(X_n)_{n\in \N}$.  Given a non-principal ultrafilter $\omega$ on $\N$, consider the ultraproduct $$\Sigma_{\omega}=\bigg(\prod_{n\in\N}S(X_n)\bigg)/N_{\omega}.$$ Write $G=F(S)/N$ with $N\normaleq F(S)$. The actions of $F(S)$ on the  $X_n$ induce a homomorphism $\theta\colon F(S) \longrightarrow \Sigma_{\omega}$ by considering each element $w\in F(S)$ as an element of the symmetric group $S(X_n)$.
We have that $\ker \theta =N$, and so this in turn induces an embedding $\varphi\colon G\longrightarrow \Sigma_{\omega}$. Furthermore, this embedding is \emph{faithful}, i.e. $d_{\omega}(\id,\varphi(g))=1$ for every $g\in G\setminus\{1\}$.

On the other hand, if we are given a faithful embedding $\varphi\colon G\longrightarrow \Sigma_{\omega}$ for some non-principal ultrafilter $\omega$ on $\N$, then we can consider a representative $$\tilde\varphi\colon G\longrightarrow \prod_{n\in\N}S_{k_n}$$ of $\varphi$. Write $\tilde{\varphi}=(\tilde{\varphi}_n)_{n\in\N}$ with $\tilde\varphi_n\colon G\longrightarrow S_{k_n}$. For each $n\in \N$, take the set $X_n=\{1,\dotsc,k_n\}$ and turn it into an $S$-labelled graph with directed edges $(i,i\cdot\tilde\varphi_n(s))$ for $i\in X_n$ labelled by $s\in S^{\pm 1}$. From the definition of the ultraproduct, we have that
\begin{displaymath}
\lim_{n\to \omega}d_{n_k}\big(\tilde\varphi_n(gh),\tilde\varphi_n(g) \tilde\varphi_n(h)\big) = 0
\end{displaymath}
for all $g,h\in G$ and
\begin{displaymath}
\lim_{n\to \omega}d_{k_n}\big(\id_{k_n},\tilde\varphi_n(g)\big) = 1
\end{displaymath}
for all $g\in G\setminus \{1\}$. This implies that
\begin{displaymath}
\lim_{n\to \omega} \frac{\abs[\big]{Q_r(X_n)}}{\abs{X_n}} = 1,
\end{displaymath}
and so 
\begin{displaymath}
\{n\in\N\mid X_n \text{ is an $r$-approximation}\}\in \omega
\end{displaymath}
for every $r\in\N$. In particular, we can extract a sofic approximation of $G$ from $(X_n)_{n\in\N}$.

The previous discussion can be summarised in the form of the following result.

\begin{proposition}
	Let $G$ be a finitely generated sofic group. Then, the following conditions are equivalent.
	\begin{enumerate}
		\item The group $G$ is sofic.
		
		\item For every non-principal ultrafilter $\omega$ on $\N$, there is an ultraproduct $\Sigma_{\omega}$ of finite symmetric groups $(S_{k_n})_{n\in\N}$ and a faithful embedding of $G$ into $\Sigma_{\omega}$.
		
		\item For some non-principal ultrafilter $\omega$ on $\N$, there is an ultraproduct $\Sigma_{\omega}$ of finite symmetric groups $(S_{k_n})_{n\in\N}$ and a faithful embedding of $G$ into $\Sigma_{\omega}$.
	\end{enumerate}
\end{proposition}

The next result is extracted from \cite{elek-szabo}, and states that all sofic approximations of an amenable group are, asymptotically speaking, conjugate. The proof of this theorem  relies heavily on the theory of quasitilings developed during Chapter~\ref{chapter:Quasitilings}.

\begin{theorem}[{\cite[Theorem~2]{elek-szabo}}] \label{thm:Sofic Embeddings Conjugate}
	Let $G$ be a finitely generated amenable group with $S\subseteq G$ a finite generating set, $\omega$  a non-principal ultrafilter on $\N$ and $\Sigma_{\omega}$  an ultraproduct of finite symmetric groups $(S_{k_n})_{n\in\N}$. Then, any two faithful embeddings $\varphi,\psi\colon G\longrightarrow \Sigma_{\omega}$ are conjugate, i.e. there exists some $\sigma\in \Sigma_{\omega}$ such that $$\varphi(g)=\sigma^{-1} (\psi(g))\sigma$$ for all $g\in G$.
\end{theorem}

\begin{proof}
	For each $n\in \N$, take $$X_n=Y_n=\{1,\dotsc,k_n\}$$ to be $S$-labelled graphs, with $X_n$ associated to $\varphi$ and $Y_n$ to $\psi$ in the same manner as before. For each $n\in\N$, let $h(n)\in\N$ be the largest integer such that both $X_n$ and $Y_n$ are $h(n)$-approximations of $\Cay(G,S)$. Then, 
	\begin{displaymath}
	\{n\in\N\mid h(n)\geq r\}\in \omega
	\end{displaymath}
	for every $r\in\N$. Given a F\o lner exhaustion $(F_n)_{n\in\N}$ of $G$ and $r\in\N$, using Theorem~\ref{thm:Quastitilings Graph 2} we can obtain a finite subsequence $\mathcal{T}=(T_1,\dotsc,T_m)$ of $(F_n)_{n\in\N}$ and some $\alpha\in\N_0^m\setminus\{0\}$ such that, if $h(n)$ is large enough, then $X_n$ is $4r$-isomorphic to $t_1\alpha\cdot\mathcal{T}$ and $Y_k$ is $4r$-isomorphic to $t_2\alpha\cdot\mathcal{T}$ with $t_1,t_2\in\N$. 
	We have that 
	\begin{displaymath}
	\lim_{n\to \omega}\frac{\abs[\big]{E(X_n)}}{\abs[\big]{E(Y_n)}}=1,
	\end{displaymath}
	and so
	for $n$ large enough we have that 
	$t_1\alpha\cdot\mathcal{T}$ is $4r$-isomorphic to $t_2\alpha\cdot\mathcal{T}$. Applying Lemma~\ref{lem:Labelled Graphs r-Isomorphic} twice we obtain that $X_n$ and $Y_n$ are $r$-isomorphic for $n$ large enough.
	
	Now, for each $n\in\N$ let $l(n)\in\N$ be the largest integer for which $X_n$ and $Y_n$ are $l(n)$-isomorphic, and note that 
	\begin{displaymath}
	\{n\in \N\mid l(n)\geq r\}\in \omega
	\end{displaymath}
	for all $r\in\N$. Let $\sigma_n\in S_{k_n}$ be a bijection with $X_n''\subseteq X_n$  and $Y_n''\subseteq Y_n$ such that $$\abs[\big]{E(X_n'')}\geq\bigg(1-\frac{1}{l(n)}\bigg)\abs[\big]{E(X_n)},\quad \abs[\big]{E(Y_n'')}\geq\bigg(1-\frac{1}{l(n)}\bigg)\abs[\big]{E(Y_n)},$$ 
	and $\sigma_n\colon X_n''\longrightarrow Y_n''$ is an isomorphism of $S$-labelled graphs. Then, if we write $\sigma=(\sigma_n)_{n\in\N}$, for every $g\in G$ and $i\in \{1,\dotsc,k_n\}$ we have that
	\begin{displaymath}
	\lim_{n\to \omega}d_{k_n}\big(i\cdot \sigma_n \tilde\varphi_{n}(g),i\cdot\tilde\psi_n(g)\sigma_n\big) = \lim_{n\to \omega} \frac{1}{l(n)} = 0,
	\end{displaymath}
	and so we can conclude that $\psi\sigma =\sigma \varphi$.
\end{proof}

\begin{remark}
	As shown in \cite{elek-szabo}, the converse of this result is also true, i.e. if $G$ is a finitely generated sofic group such that any two faithful embeddings into an ultraproduct $\Sigma_{\omega}$ are conjugate, then $G$ is amenable. Hence, it characterises finitely generated amenable groups as those finitely generated sofic groups for which all their sofic approximations are conjugate. 
\end{remark}

%% file: Chapter4.tex
\chapter{The Sofic L\"{u}ck Approximation Conjecture}
\label{chapter:LuckConjecture}

In this chapter we will introduce the Sofic L\"uck Approximation Conjecture, a version of a conjecture that has its origins in the study of $L^2$-invariants. We will begin by explaining the general statement of the conjecture. Afterwards, we will show how the conjecture can be proved for amenable groups by making use of the techniques developed in the previous chapters. We will then conclude the chapter with a brief discussion of the proof of the conjecture for general groups over the field $\Q$, in order to motivate the results that will be developed in the next chapter. This chapter is primarily based on  \cite[\S2,\S7,\S10]{jaikin} and \cite[\S4]{hindman}.

\section{Statement of the Conjecture}

Let $G$ be a finitely generated sofic group with a finite generating subset $S\subseteq G$ and a sofic approximation $(X_n)_{n\in\N}$, and let $K$ be a field. For each $n\in\N$, the free group $F(S)$ acts on $X_n$, and so given any matrix $A\in\Mat_{k\times l}(K[F(S)])$ we can consider the linear map of $K$-vector spaces $\phi_{X_n}^A\colon K[X_n]^k \longrightarrow K[X_n]^l$ defined by $$\phi_{X_n}^A(x_1,\dotsc,x_k)= (x_1,\dotsc,x_k) A.$$
Now, we can define the \emph{rank of $A$ with respect to $X_n$} as
\begin{displaymath}
\rk_{X_n}(A)=k-\frac{\dim_K \ker \phi_{X_n}^A}{\abs{X_n}}=\frac{\dim_K \im \phi_{X_n}^A}{\abs{X_n}}.
\end{displaymath}
We will be interested in the study of the convergence of these ranks. More specifically, we want to know whether these ranks converge, and whether convergence depends in any way on the sofic approximation that we have chosen. We will now state the Sofic L\"uck Approximation Conjecture, which tries to provide an answer to these questions.

\begin{conjecture}[L\"uck]\label{conj:Sofic Luck}
	Let $G$ be a finitely generated group with $S\subseteq G$ a finite generating set and $(X_n)_{n\in\N}$ a sofic approximation of $G$. Let $K$ be a field and $A\in\Mat_{k\times l}(K[F(S)])$. Then, the following hold:
	\begin{enumerate}
		\item The limit $\lim_{n\to \infty}\rk_{X_n}(A)$ exists.
		
		\item The limit $\lim_{n\to \infty}\rk_{X_n}(A)$ is independent of the sofic approximation $(X_n)_{n\in\N}$.
	\end{enumerate}
\end{conjecture}

It is not currently known whether the Sofic L\"uck Approximation Conjecture holds in general, i.e. for any sofic group and over any field. Nevertheless, in some particular instances, such as when the group is amenable, or when the field has characteristic $0$, the conjecture has been proven to be true.

\section{The Conjecture for Amenable Groups}

The Sofic L\"uck Approximation Conjecture can easily be shown to be true over any field in the case where our group is amenable by using Theorem~\ref{thm:Sofic Embeddings Conjugate}.

We will study the case where $G$ is a finitely generated amenable group with a finite generating subset $S\subseteq G$ and a sofic approximation $(X_n)_{n\in\N}$. 
For the sake of simplicity, we will only consider elements of the group algebra and not matrices. Hence, if $K$ is a field and $a\in K[F(S)]$, for each $n\in \N$ we have a $K$-linear map $\phi_{X_n}^a\colon K[X_n] \longrightarrow K[X_n]$, and 
\begin{displaymath}
\rk_{X_n}(a)=1-\frac{\dim_K \ker \phi_{X_n}^a}{\abs{X_n}}=\frac{\dim_K \im \phi_{X_n}^a}{\abs{X_n}}.
\end{displaymath}

Observe then that the convergence of the ranks $\rk_{X_n}(a)$ is equivalent to the convergence of the normalised dimensions 
\begin{displaymath}
\frac{\dim_K \ker \phi_{X_n}^a}{\abs{X_n}}.
\end{displaymath}
As a consequence, in order to prove the Sofic L\"uck Approximation Conjecture, we will study the normalised dimensions of these kernels.

\begin{theorem}\label{thm:preLuck holds Amenable}
	Let $G$ be a finitely generated amenable group with $S\subseteq G$ a finite generating subset and $(X_n)_{n\in\N}$ a sofic approximation of $G$. Given any field $K$, an element $a\in K[F(S)]$ and a non-principal ultrafilter $\omega$ on $\N$, the limit
	\begin{displaymath}
	\lim_{n\to\omega}\frac{\dim_K \ker \phi_{X_n}^a}{\abs{X_n}}
	\end{displaymath}
	exists and is independent of the sofic approximation $(X_n)_{n\in\N}$.
\end{theorem}
\begin{proof}
	Any subsequence of a sofic approximation of $G$ is once again a sofic approximation of $G$. Thus, if the statement of the theorem fails, we can find two sofic approximations $(X_n^{1})_{n\in\N}$ and $(X_n^{2})_{n\in\N}$ of $G$ such that $$\lim_{n\to \omega}\frac{\dim_K \ker \phi_{X_n^i}^a}{\abs{X_n^i}}$$ exists for $i=1,2$ but
	\begin{equation}\label{eqn:preLuck amenable 1}
	\lim_{n\to \omega}\frac{\dim_K \ker \phi_{X_n^1}^a}{\abs{X_n^1}} \not= \lim_{n\to \omega}\frac{\dim_K \ker \phi_{X_n^2}^a}{\abs{X_n^2}}.
	\end{equation}
	Assume by contradiction that two such sofic approximations of $G$ exist. 
	Consider then
	\begin{displaymath}
	Y_n^{1}=Y_n^{2} = X_n^{1}\times X_n^{2},
	\end{displaymath}
	with $F(S)$ acting on $Y_n^{i}$ by acting on the $i$-th coordinate for $i=1,2$. Then, $(Y_n^1)_{n\in\N}$ and $(Y_n^2)_{n\in\N}$ are both sofic approximations of $G$. Furthermore, we have that $\abs{Y_n^1}=\abs{Y_n^2}$  and 
	\begin{equation}\label{eqn:preLuck amenable 2}
	\frac{\dim_K \ker \phi_{Y_n^i}^a}{\abs{Y_n^i}}=\frac{\dim_K \ker \phi_{X_n^i}^a}{\abs{X_n^i}}
	\end{equation}
	for $i=1,2$ and for all $n\in\N$.
	
	By Theorem~\ref{thm:Sofic Embeddings Conjugate}, for each $n\in\N$ there is some bijection $\sigma_n\colon Y_n^1\longrightarrow Y_n^2$ such that, if we denote by
	\begin{displaymath}
	{Y_n^1}'=\big\{x\in Y_n^1\mid (\sigma_n^{-1}\circ\phi_{Y_n^2}^a\circ\sigma_n)(x)=\phi_{Y_n^1}^a(x)\big\}, \qquad {Y_n^2}'=\sigma_n({Y_n^1}'),
	\end{displaymath}
	then
	\begin{displaymath}
	\lim_{n\to \omega}\frac{\abs{{Y_n^1}'}}{\abs{Y_n^1}}=\lim_{n\to \omega}\frac{\abs{{Y_n^2}'}}{\abs{Y_n^2}}=1.
	\end{displaymath}
	Observe that, given $x\in K[{Y_n^1}']$, we have that $x\in\ker \phi_{Y_n^1}$ if and only if $\sigma_n(x)\in\ker\phi_{Y_n^2}$. 
	Consider then the restrictions of $\phi_{Y_n^1}^a$ to $K[{Y_n^1}']$ and of $\phi_{Y_n^2}^a$ to $K[{Y_n^2}']$, which we will denote by $\phi_{{Y_n^1}'}^a$ and $\phi_{{Y_n^2}'}^a$, respectively. Thus, we have that
	\begin{displaymath}
	\ker \phi_{{Y_n^1}'}^a=\ker \phi_{Y_n^1}^a \cap K[{Y_n^1}'], 
	\qquad \ker \phi_{{Y_n^2}'}^a =\ker \phi_{Y_n^2}^a \cap K[{Y_n^2}'],
	\end{displaymath}
	and
	\begin{displaymath}
	\ker \phi_{{Y_n^1}'}^a \cong \ker \phi_{{Y_n^2}'}^a
	\end{displaymath}
	for every $n\in\N$. Furthermore, the Second Isomorphism Theorem implies that $$\ker \phi_{Y_n^1}^a/\ker \phi_{{Y_n^1}'}^a \lesssim K[Y_n^1]/K[{Y_n^1}'].$$ 
	As a consequence, 
	we obtain that 
	\begin{align*}
		\dim_K\ker \phi_{Y_n^1}^a&=\dim_K\ker \phi_{{Y_n^1}'}^a+\dim_K(\ker \phi_{Y_n^1}^a/\ker \phi_{{Y_n^1}'}^a) \\
		& \leq \dim_K\ker \phi_{{Y_n^1}'}^a+ \dim_K\big(K[Y_n^1]/K[{Y_n^1}']\big) \\
		& = \dim_K\ker \phi_{{Y_n^1}'}^a + \abs{Y_n^1}-\abs{{Y_n^1}'}
	\end{align*}
	for all $n\in\N$. 
	Analogously, we obtain the inequality 
	\begin{displaymath}
	\dim_K\ker \phi_{Y_n^2}^a\leq \dim_K\ker \phi_{{Y_n^2}'}^a + \abs{Y_n^2}-\abs{{Y_n^2}'}
	\end{displaymath}
	for all $n\in\N$. 
	Therefore, we obtain that 
	\begin{align*}
		\lim_{n\to \omega}\frac{\dim_K\ker \phi_{Y_n^1}^a}{\abs{Y_n^1}} & = \lim_{n\to \omega}\frac{\dim_K\ker \phi_{{Y_n^1}'}^a}{\abs{Y_n^1}} \\
		& = \lim_{n\to \omega}\frac{\dim_K\ker \phi_{{Y_n^2}'}^a}{\abs{Y_n^2}} \\
		& = \lim_{n\to \omega}\frac{\dim_K\ker \phi_{Y_n^2}^a}{\abs{Y_n^2}},
	\end{align*}
	which along with (\ref{eqn:preLuck amenable 2}) contradicts (\ref{eqn:preLuck amenable 1}). Hence, we can conclude that 
	\begin{displaymath}
	\lim_{n\to\omega}\frac{\dim_K \ker \phi_{X_n}^a}{\abs{X_n}}
	\end{displaymath}
	exists and is independent of the sofic approximation $(X_n)_{n\in\N}$.
%
%
%
\end{proof}

This result now automatically gives us a proof of the conjecture in the case that our group is amenable. 

\begin{theorem}\label{thm:Luck holds Amenable}
	The Sofic L\"uck Approximation Conjecture holds for finitely generated amenable groups.
\end{theorem}

\begin{proof}
	Let $G$ be a finitely generated amenable group with a finite generating set $S\subseteq G$
	. Let $K$ be a field and $a\in K[F(S)]$.
	
	Any subsequence of a sofic approximation of $G$ is once again a sofic approximation of $G$. Thus, if either condition (i) or (ii) in Conjecture~\ref{conj:Sofic Luck} fails, we can find two sofic approximations $(X_n)_{n\in\N}$ and $(Y_n)_{n\in\N}$ of $G$ such that both  limits $$\lim_{n\to \infty}\rk_{X_n}(a), \quad \lim_{n\to \infty}\rk_{Y_n}(a)$$ exist but
	\begin{displaymath}
	\lim_{n\to \infty}\rk_{X_n}(a) \not= \lim_{n\to \infty}\rk_{Y_n}(a).
	\end{displaymath}
	Using the same argument as in the proof of Theorem~\ref{thm:preLuck holds Amenable}, we can assume without loss of generality that $\abs{X_n}=\abs{Y_n}$ for all $n\in\N$.
	
	Now, for any non-principal ultrafilter $\omega$ on $\N$, we have that
	\begin{displaymath}
	\lim_{n\to\omega}\frac{\dim_K \ker \phi_{X_n}^a}{\abs{X_n}}=\lim_{n\to\omega}\frac{\dim_K \ker \phi_{Y_n}^a}{\abs{Y_n}}
	\end{displaymath}
	as a consequence of Theorem~\ref{thm:preLuck holds Amenable}. But then,
	\begin{align*}
	\lim_{n\to \infty} \rk_{X_n}(a) & = \lim_{n\to \omega} \rk_{X_n}(a) \\
	& = 1 - \lim_{n\to \omega} \frac{\dim_K \ker \phi_{X_n}^a}{\abs{X_n}} \\
	& = 1 - \lim_{n\to \omega} \frac{\dim_K \ker \phi_{Y_n}^a}{\abs{Y_n}} \\
	& = \lim_{n\to \omega} \rk_{Y_n}(a) \\
	& = \lim_{n\to \infty} \rk_{Y_n}(a),
	\end{align*}
	which is a contradiction. Therefore, we conclude that the Sofic L\"uck Approximation Conjecture holds in this case.
\end{proof}

\section{The Conjecture over the Field $\Q$}

We will now discuss the proof of the Sofic L\"{u}ck Approximation Conjecture over the field $\Q$ for general groups. We will refrain from giving all the details, as our main goal is to motivate the techniques that we will develop in the next chapter. For a detailed proof, see \cite{jaikin}.

Let $G$ be a finitely generated sofic group with $S\subseteq G$ a finite generating set, $G=F(S)/N$ with $N\normaleq F(S)$ 
and $(X_n)_{n\in\N}$ a sofic approximation of $G$. For the sake of simplicity, we will once more restrict our attention to the case where $a\in \Q[F(S)]$. This element defines a linear map of $\Q$-vector spaces $\phi_{X_n}^a\colon \Q[X_n] \longrightarrow \Q[X_n]$ by $$\phi_{X_n}^a(x)= x a.$$
Then, the rank of $a$ with respect to $X_n$ is
\begin{displaymath}
\rk_{X_n}(a)=1-\frac{\dim_K \ker \phi_{X_n}^a}{\abs{X_n}}=\frac{\dim_K \im \phi_{X_n}^a}{\abs{X_n}}.
\end{displaymath}

The group $G$ acts on the Hilbert space $\ell^2(G)$ by both left and right multiplication. The element $a\in\Q[F(S)]$ thus defines a bounded operator $\phi_G^a\colon \ell^2(G)\longrightarrow\ell^2(G)$ by
\begin{displaymath}
\phi_G^a(v)=va
\end{displaymath}
for $v\in\ell^2(G)$. Given a left-invariant closed subspace $V\leq \ell^2(G)$, we can consider $\proj_V\colon \ell^2(G) \longrightarrow V$, the orthogonal projection onto $V$. Then, we can define 
\begin{displaymath}
\dim_G V = \scalar[\big]{\proj_V(1_G),1_G}.
\end{displaymath}

In order to prove the Sofic L\"uck Approximation Conjecture, we will work towards proving that
\begin{displaymath}
\lim_{n\to \infty}\frac{\dim_{\Q}\ker \phi_{X_n}^a}{\abs{X_n}} = \dim_G \ker \phi_G^a
\end{displaymath}
independent of the sofic approximation $(X_n)_{n\in\N}$. In order to do this, we will construct a sequence of measures whose values at zero coincide with the normalised dimensions of the kernels of the associated operators. Then, we will show that these measures converge weak-$*$, and finally that their values at zero converge.

Write the element $a\in\Q[F(S)]$ as $$a=\sum_{w\in F(S)} a_w w$$ with $a_w\in \Q$ for  $w\in F(S)$. The adjoint of $a$ is then the element  $$a^*=\sum_{w\in F(S)} 
a_w w^{-1}.$$
Given any $x\in \Q[X_n]$, we have that $xa=0$ if and only if $x(aa^*)=0$, and so 
$$\dim_{\Q} \ker \phi_{X_n}^a= \dim_{\Q} \ker \phi_{X_n}^{aa^*}.$$
Moreover, the adjoint of the operator $\phi_{X_n}^a$ is $(\phi_{X_n}^a)^*=\phi_{X_n}^{a^*}$. Analogously, we have that 
$$\dim_G \ker \phi_{G}^a= \dim_G \ker \phi_{G}^{aa^*}$$
and $(\phi_{G}^a)^*=\phi_{G}^{a^*}$.

Therefore, we may assume that $a=bb^*$ for some $b\in \Q[F(S)]$, and so $a=a^*$. Thus, $\phi_{X_n}^a$ and $\phi_{G}^a$ are positive self-adjoint operators, and so their spectra are compact and contained in $[0,\infty)$. In fact, we have that $$\spec \phi_{X_n}^a \subseteq \big[0,\norm{\phi_{X_n}^a}\big],\quad \spec \phi_{G}^a \subseteq \big[0,\norm{\phi_{G}^a}\big].$$

For each $n\in \N$, we can define a probability measure 
\begin{displaymath}
\mu_{X_n}^a = \frac{1}{\abs{X_n}} \sum_{\lambda\in\spec \phi_{X_n}^a} \delta_{\lambda}
\end{displaymath}
on $\spec \phi_{X_n}^a$, where $\delta_{\lambda}$ denotes the Dirac measure concentrated at the point $\lambda
$. Moreover, we have that 
\begin{equation}\label{eqn:Mu(0) Xn}
\mu_{X_n}^a\big(\{0\}\big)=\frac{\dim_{\Q}\ker\phi_{X_n}^a}{\abs{X_n}}.
\end{equation}

We can also associate a probability measure to the operator $\phi_G^a$, using the concept of spectral measures of self-adjoint operators. Let $H$ be a Hilbert space, $A\in \mathcal{B}(H)$ be a bounded self-adjoint operator, and $v\in H$. Then, there exists a unique positive Radon measure $\mu_{A,v}$ on $\spec A$, called the \emph{spectral measure associated to $A$ and $v$}, satisfying that 
\begin{displaymath}
\int_{\spec A}f\ \dd \mu_{A,v} = \scalar[\big]{f(A)v,v}
\end{displaymath}
for every continuous function $f$ on $\spec A$. Furthermore, we have that 
\begin{displaymath}
\mu_{A,v}(\spec A)=\norm{v}^2<\infty.
\end{displaymath}
For more information on spectral measures associated to bounded self-adjoint operators, see \cite[\S 3]{kowalski}.

Because $\phi_{G}^a$ is a self-adjoint operator, we can thus define a probability measure $\mu_{G}^a$ on $\spec \phi_{G}^a$ associated to $\phi_{G}^a$ by  
\begin{displaymath}
\mu_{G}^a=\mu_{\phi_{G}^a,1_G}.
\end{displaymath}
Moreover, we can show that 
\begin{equation}\label{eqn:Mu(0) G} 
\mu_{G}^a\big(\{0\}\big)=\scalar[\big]{\chi_0(\phi_{G}^a)1_G,1_G}=\dim_{G}\ker\phi_{G}^a.
\end{equation}

Because of \ref{eqn:Mu(0) Xn} and \ref{eqn:Mu(0) G}, our goal will now be to prove that 
\begin{displaymath}
\lim_{n\to \infty}\mu_{X_n}^a\big(\{0\}\big)=\mu_G^a\big(\{0\}\big).
\end{displaymath}
Before doing this, we will need some compact space on which all of the measures are defined, which will then allow us to show that there is weak-$*$ convergence. 

It is possible to find a uniform bound for the norms of the operators $\phi_{X_n}^a$. If we denote by
\begin{displaymath}
S(a)=\abs{\{w\in F(S)\mid a_w\not=0\}}
\end{displaymath}
the size of the support of $a\in \Q[F(S)]$ and set
\begin{displaymath}
	\abs{a}=\sum_{w\in F(S)}\abs{a_w},
\end{displaymath}
we can prove the following result.

\begin{lemma}
	Given $a\in \Q[F(S)]$ with $a=a^*$, we have that 
	\begin{displaymath}
	\norm{\phi_{X_n}^a}\leq S(a)\abs{a}
	\end{displaymath}
	for all $n\in\N$.
\end{lemma}

As a consequence of this bound, we can deduce that there is some constant $c>0$ such that $$\norm{\phi_{G}^a}\leq c,\quad \norm{\phi_{X_n}^a}\leq c$$ for all $n\in\N$. Due to the fact that 
\begin{displaymath}
\spec \phi_{X_n}^a \subseteq \big[0,\norm{\phi_{X_n}^a}\big],\quad \spec \phi_{G}^a \subseteq \big[0,\norm{\phi_{G}^a}\big],
\end{displaymath}
this leads us to the conclusion that $\mu_G^a$ and $\mu_{X_n}^a$ are probability measures on the same interval $[0,c]$ for all $n\in\N$.

The set of complex regular Radon measures on the compact Hausdorff space $[0,c]$ is identified by the Riesz Representation Theorem with the dual space $\mathcal{C}([0,c])^*$. As such, we say that a sequence of measures $\mu_n\in \mathcal{C}([0,c])^*$ converges weak-$*$ to $\mu\in \mathcal{C}([0,c])^*$ if
\begin{displaymath}
\lim_{n\to\infty}\int_0^c f\ \dd\mu_n= \int_0^c f\ \dd\mu 
\end{displaymath} 
for all $f\in \mathcal{C}([0,c])$.

Now, we will see that $(\mu_{X_n}^a)_{n\in \N}$ converges weak-$*$ to $\mu_{G}$.

\begin{lemma}\label{lem:Integral Luck}
	For every $n,l\in\N$, we have that
	\begin{displaymath}
	\int_0^c t^l \ \dd\mu_{X_n}^a = \frac{1}{\abs{X_n}}\tr \big((\phi_{X_n}^a)^l\big).
	\end{displaymath}
\end{lemma}
\begin{proof}
	From the definition of $\mu_{X_n}^a$, it follows that
	\begin{displaymath}
	\int_0^c t^l \ \dd\mu_{X_n}^a =
	\frac{1}{\abs{X_n}}\sum_{\lambda\in\spec \phi_{X_n}^a} \lambda^l
	= \frac{1}{\abs{X_n}}\tr \big((\phi_{X_n}^a)^l\big).
	\end{displaymath}
\end{proof}

\begin{lemma}
	For every $l\in\N$, we have that 
	\begin{displaymath}
	\lim_{n\to\infty}\int_0^c t^l \ \dd\mu_{X_n}^a = \int_0^c t^l \ \dd\mu_{G}^a
	\end{displaymath}
	independent of the sofic approximation $(X_n)_{n\in\N}$.
\end{lemma}
\begin{proof}
	Note that 
	\begin{displaymath}
	\int_0^c t^l \ \dd\mu_{X_n}^a =
	\int_0^c t \ \dd\mu_{X_n}^{a^l},
	\end{displaymath}
	and the analogous result is true for $\mu_G^a$, so we may assume that $l=1$. In light of Lemma~\ref{lem:Integral Luck}, we need to study the limit 
	\begin{displaymath}
	\lim_{n\to\infty}\frac{1}{\abs{X_n}}\tr \phi_{X_n}^a.
	\end{displaymath}
	Further note that if $$a=\sum_{w\in F(S)} a_w w,$$ then
	\begin{displaymath}
	\phi_{X_n}^a=\sum_{w\in F(S)} a_w \phi_{X_n}^w,
	\end{displaymath}
	and so 
	\begin{displaymath}
	\frac{1}{\abs{X_n}}\tr \phi_{X_n}^a=\sum_{w\in F(S)}a_w \frac{1}{\abs{X_n}}\tr \phi_{X_n}^w.
	\end{displaymath}
	Analogously, $$\phi_{G}^a=\sum_{w\in F(S)} a_w \phi_{G}^w,$$
	which implies that 
	\begin{displaymath}
	\int_0^c t\ \dd\mu_G^a= \sum_{w\in F(S)} a_w \scalar[\big]{\phi_{G}^w(1_G),1_G}=\sum_{w\in F(S)}\int_0^c t\ \dd\mu_G^w.
	\end{displaymath}
	We may thus assume that $a=w\in F(S)$. Then, we have that
	\begin{displaymath}
	\tr \phi_{X_n}^w = \abs[\big]{\{x\in X_n\mid xw=x\}},
	\end{displaymath}
	and due to the fact that $(X_n)_{n\in \N}$ is a sofic approximation we get that
	\begin{displaymath}
	\lim _{n\to\infty}\frac{1}{\abs{X_n}}\tr \phi_{X_n}^w = 
	\begin{cases}
	1 & \text{if } w\in N,\\
	0 & \text{if } w\not\in N.
	\end{cases}
	\end{displaymath}
	Finally, 
	\begin{displaymath}
	\int_0^c t\ \dd\mu_G^w = \scalar{w,1_G}=\begin{cases}
	1 & \text{if } w\in N,\\
	0 & \text{if } w\not\in N,
	\end{cases}
	\end{displaymath}
	and so 
	\begin{displaymath}
	\lim _{n\to\infty}\int_0^c t\ \dd\mu_{X_n}^w=\int_0^c t\ \dd\mu_G^w.
	\end{displaymath}
\end{proof}

The Weierstrass Approximation Theorem says that the polynomials form a dense subset of $\mathcal{C}([0,c])$, which along with the Bounded Convergence Theorem and the previous lemma gives us the following result.

\begin{proposition}\label{lem:Weak-* Conv Measures Luck}
	The sequence of probability measures $(\mu_{X_n}^a)_{n\in\N}$ on $[0,c]$ converges weak-$*$ to $\mu_G^a$.
\end{proposition}

Now that we have weak-$*$ convergence of the measures, our goal will be to prove the convergence at the point $0$.

\begin{remark}
	In general, weak-$*$ convergence of the measures does not imply convergence of the value at any given point. As an example, consider the sequence of measures $(\delta_{1/n})_{n\in\N}$ on $[0,1]$, which converges weak-$*$ to $\delta_0$. Nonetheless, $\delta_{1/n}(\{0\})=0$ for all $n\in \N$, whereas $\delta_0(\{0\})=1$.
\end{remark}

In order to prove our result, we will make use of a classical theorem in measure theory.
\begin{theorem}[Portmanteau]\label{thm:Portmanteau Measures}
	Let $\Omega$ be a compact metric space. Let $\mu$ and $\mu_n$ for $n\in\N$ be Borel probability measures on $\Omega$. Then, the following are equivalent:
	\begin{enumerate}
		\item The sequence $(\mu_n)_{n\in\N}$ converges weak-$*$ to $\mu$.
		
		\item For every closed subset $V\subseteq \Omega$ we have that 
		\begin{displaymath}
		\limsup_{n\to\infty} \mu_n(V) \leq \mu(V).
		\end{displaymath}
		
		\item For every open subset $U\subseteq \Omega$ we have that 
		\begin{displaymath}
		\liminf_{n\to\infty} \mu_n(U) \geq \mu(U).
		\end{displaymath}
	\end{enumerate}
\end{theorem}

As a consequence of Proposition~\ref{lem:Weak-* Conv Measures Luck}, the Portmanteau Theorem automatically tells us that 
\begin{equation}\label{eqn:limsup mu Xn 0}
\limsup_{n\to\infty} \mu_{X_n}^a\big(\{0\}\big) \leq \mu_{G}^a\big(\{0\}\big),
\end{equation}
so we are left with the other inequality to prove.

In principle, we have that $a\in\Q[F(S)]$.  Nevertheless, multiplying by some constant if necessary, we can assume that $a\in\Z[F(S)]$. Then, for each $n\in\N$ we define
\begin{displaymath}
{\det}_+ \phi_{X_n}^a = \prod_{\lambda\in\spec\phi_{X_n}^a\setminus\{0\}} \lambda.
\end{displaymath}
Under the assumption that $a\in\Z[F(S)]$, we can easily prove the following result.

\begin{lemma}\label{lem:det+ in N}
	For every $n\in\N$ we have that $\det_+ \phi_{X_n}^a\in\N$. 
\end{lemma}
\begin{proof}
	Because $a\in\Z[F(S)]$, every element in $\spec\phi_{X_n}^a$ is an algebraic integer. Furthermore, $\det_+ \phi_{X_n}^a$ is invariant under the action of the Galois group of the splitting field of the characteristic polynomial of $\phi_{X_n}^a$, so $\det_+ \phi_{X_n}^a\in\Q$. Therefore, $\det_+ \phi_{X_n}^a\in\Z$. Finally, since $\spec\phi_{X_n}^a\subseteq [0,c]$, we have that $\det_+ \phi_{X_n}^a>0$, and so $\det_+ \phi_{X_n}^a\in\N$.
\end{proof}

Given $0<\varepsilon<1$, we can find a uniform bound for $\mu_{X_n}^a((0,\varepsilon))$ with the help of the previous result. 

\begin{lemma}\label{eqn:Luck bound}
	Given $0<\varepsilon<1$, we have that 
	\begin{equation*}
	\mu_{X_n}^a\big((0,\varepsilon)\big)\leq -\frac{\log c}{\log\varepsilon}
	\end{equation*}
	for every $n\in\N$.
\end{lemma}

\begin{proof}
	By Lemma~\ref{lem:det+ in N}, we have that $\det_+ \phi_{X_n}^a\in\N$, and so $\det_+ \phi_{X_n}^a\geq 1$. On the other hand, 
	\begin{align*}
		{\det}_+\phi_{X_n}^a & =\prod_{\lambda\in\spec\phi_{X_n}^a\setminus\{0\}} \lambda \\
		& = \bigg(\prod_{\substack{\lambda\in\spec\phi_{X_n}^a \\ 0<\lambda<\varepsilon}} \lambda\bigg) \cdot \bigg( \prod_{\substack{\lambda\in\spec\phi_{X_n}^a \\ \lambda\geq \varepsilon}} \lambda \bigg) \\
		& \leq \varepsilon^{\abs{X_n}\mu_{X_n}^a\big((0,\varepsilon)\big)} \cdot c^{\abs{X_n}},
	\end{align*}
	and so we have that
	\begin{displaymath}
	\varepsilon^{\abs{X_n}\mu_{X_n}^a\big((0,\varepsilon)\big)} \cdot c^{\abs{X_n}}\geq 1.
	\end{displaymath}
	Taking logarithms, we obtain that
	\begin{displaymath}
	\abs{X_n}\Big(\mu_{X_n}^a\big((0,\varepsilon)\big)\log\varepsilon + \log c\Big) \geq 0.
	\end{displaymath}
	Hence, we have that
	\begin{equation*}
		\mu_{X_n}^a\big((0,\varepsilon)\big)\leq -\frac{\log c}{\log\varepsilon}.
	\end{equation*}
\end{proof}

Applying the Portmanteau Theorem once again along with the bound in Lemma~\ref{eqn:Luck bound}, we obtain that
\begin{align*}
\mu_G^a\big(\{0\}\big) & \leq \mu_g\big([0,\varepsilon)\big) \\
& \leq \liminf_{n\to\infty} \mu_{X_n}^a\big([0,\varepsilon)\big) \\
& \leq \liminf_{n\to\infty}\mu_{X_n}^a\big(\{0\}\big) - \frac{\log c}{\log \varepsilon}.
\end{align*}
Thus, if we make $\varepsilon$ tend to $0$, we get that
\begin{equation}\label{eqn:liminf mu Xn 0}
\mu_G^a\big(\{0\}\big)\leq \liminf_{n\to\infty}\mu_{X_n}^a\big(\{0\}\big).
\end{equation}
Therefore, the inequalities (\ref{eqn:limsup mu Xn 0}) and (\ref{eqn:liminf mu Xn 0}) imply that
\begin{displaymath}
\lim_{n\to \infty}\mu_{X_n}^a\big(\{0\}\big)= \mu_G^a\big(\{0\}\big).
\end{displaymath}
Consequently, we obtain the following result. 

\begin{theorem}[L\"uck]\label{thm:Luck Conj Q}
	Let $G$ be a finitely generated group with $S\subseteq G$ a finite generating set and $a\in\Q[F(S)]$. Then, for every sofic approximation $(X_n)_{n\in\N}$ of $G$ we have that
	\begin{displaymath}
	\lim_{n\to \infty}\frac{\dim_{\Q}\ker \phi_{X_n}^a}{\abs{X_n}} = \dim_G \ker \phi_G^a.
	\end{displaymath}
	In particular, the Sofic L\"uck Approximation Conjecture holds over $\Q$.
\end{theorem}

%% file: Chapter5.tex
\chapter{Convergence of Adelic Measures Associated to Sofic Representations}
\label{chapter:LastChapter}

In this final chapter we will present a generalisation of the Sofic L\"uck Approximation Conjecture, first for discrete valuation domains and then for rings of integers of number fields. Associated to each of the operators appearing in the conjecture, we will construct a probability measure on the space of ideals of our ring, and then study the convergence of these measures for amenable groups, proving that they converge strongly to some limit measure. The main results in this chapter are original. Some of the auxiliary results in the second section have been taken from \cite[\S 8.3, \S 8.4]{jaikin2}.


\section{Approximation of Local Measures}


Let $\OO$ be a discrete valuation domain, i.e. a principal ideal domain with a unique non-zero maximal ideal $\m$, and let $K$ be the field of fractions of $\OO$. 
%
If the ideal $\m$ is generated by the prime element $\pi\in \Oo$, then every non-trivial ideal of $\Oo$ is of the form $\m^i=\pi^i\Oo$ with $i\in\N$. We will denote the set of ideals of $\Oo$ by
\begin{displaymath}
\I=\{0,\Oo,\m,\m^2,\dotsc\}.
\end{displaymath}
Furthermore, given an ideal $\m^i\in\I$ with $i\in\N$ we will write
\begin{displaymath}
[0,\m^i]=\{0,\m^i,\m^{i+1},\dotsc\}.
\end{displaymath}

Let $G$ be a finitely generated amenable group with a finite generating subset $S\subseteq G$ and a sofic approximation $(X_n)_{n\in\N}$. Consider then 
%
%
%
%
%
an element $a\in \Oo[F(S)]$ and the associated linear map of $K$-vector spaces $\phi_{X_n}^a\colon K[X_n] \longrightarrow K[X_n]$ for each $n\in \N$. Then, $\phi_{X_n}^a$ can be associated to a matrix $A_{n}\in \Mat_{\abs{X_n}}(\Oo)$.

\begin{remark}
	We will only consider elements in $\Oo[F(S)]$, for given any element $a\in K[F(S)]$, we can always multiply it by some constant $\lambda\in \OO$ so that $\lambda a\in \OO[F(S)]$.
\end{remark}

We will now use of the Smith normal form of a matrix defined over a principal ideal domain in order to construct a measure associated to $\phi_{X_n}^a$ for each $n\in\N$.

\begin{proposition}[Smith normal form]
	Let $R$ be a principal ideal domain and $A\in\Mat_{k}(R)$. Then, there exist invertible matrices $P,Q\in\GL_k(R)$ 
	and a diagonal matrix 
	\begin{displaymath}
	D=\begin{pmatrix}
	\alpha_1 & & & &\\
	& \ddots & & & &\\
	& & \alpha_t& & & \\
	& & & 0 & & \\
	& & & & \ddots & \\
	& & & & & 0
	\end{pmatrix}\in \Mat_{k}(R)
	\end{displaymath}
	with $\alpha_i\mid \alpha_{i+1}$ for all $i=1,\dotsc,t-1$, such that $A=PDQ$. Furthermore, the elements $\alpha_1,\dotsc,\alpha_t\in R$ are unique up to multiplication by units. The matrix $D$ is called the \emph{Smith normal form} of $A$.
\end{proposition}

For a proof of the existence and uniqueness of the Smith normal form, see \cite{jacobson}.

\begin{remark}\label{rmk:Smith normal form module}
	Given a matrix $A\in\Mat_{k}(R)$ with Smith normal form $D$, if $\alpha_1,\dotsc,\alpha_t\in R$ are the non-zero elements that appear in the diagonal of $D$, then the  $R$-module $R^k/R^k A$ can be written as 
	\begin{displaymath}
	R^k/R^k A\isom R/\alpha_1 R\oplus \dotsb\oplus R/\alpha_t R\oplus R^r,
	\end{displaymath}
	with $r\geq 0$ being the number of zeroes in the diagonal of $D$. This decomposition of $R^k/R^k A$ is the one given by the Structure Theorem of finitely generated modules over principal ideal domains.
\end{remark}

Using the Smith normal form of $A_n$, we can assume that $\phi_{X_n}^a$ is associated to a diagonal matrix of the form 
\begin{displaymath}
D_n=\begin{pmatrix}
\pi^{k_1} & & & &\\
& \ddots & & & &\\
& & \pi^{k_t}& & & \\
& & & 0 & & \\
& & & & \ddots & \\
& & & & & 0
\end{pmatrix}\in\Mat_{\abs{X_n}}(\Oo)
\end{displaymath}
with $k_i\leq k_{i+1}$ for all $i=1,\dotsc, t-1$. This allows us to define a probability measure $$\nu_{X_n}^a=\frac{1}{\abs{X_n}}\Bigg(\sum_{i=1}^t \delta_{\m^{k_i}} + \sum_{i=t+1}^{\abs{X_n}}\delta_0\Bigg)$$ on the space of ideals $\I$ for each $n\in\N$. Our goal will now be to prove that these measures converge at each ideal 
independent of the chosen sofic approximation $(X_n)_{n\in\N}$.

Firstly, observe that the measure at zero is
\begin{displaymath}
\nu_{X_n}^a\big(\{0\}\big) =\frac{\dim_K\ker\phi_{X_n}^a}{\abs{X_n}}= 1-\rk_{X_n}(a) 
\end{displaymath}
for each $k\in\N$. 
Therefore, convergence at zero is equivalent to the Sofic L\"uck Approximation Conjecture, which holds for amenable groups by Theorem~\ref{thm:Luck holds Amenable}. Thus, we have the following result.

\begin{lemma}\label{lem: Measure local zero}
	The limit 
	\begin{displaymath}
	\lim_{n\to \infty}\nu_{X_n}^a\big(\{0\}\big)
	\end{displaymath}
	exists and is independent of the sofic approximation $(X_n)_{n\in \N}$.
\end{lemma}

We will now seek a formula for the measures at each $\m^i\in \I$. In order to do this, we need to introduce the concept of length of a module over a ring.

Let $R$ be a commutative unitary ring. Given an $R$-module $M$, we can define the \emph{length of $M$ over $R$}, which we denote by $L_R (M)$, as the supremum of the lengths of chains of $R$-submodules of the form
\begin{displaymath}
0=M_0\subsetneq M_1 \subsetneq \dotsb \subsetneq M_k=M.
\end{displaymath}
This concept serves as a generalisation for modules of the concept of dimension for vector spaces.

The length function satisfies some key properties.
\begin{itemize}
	\item If the chain of $R$-submodules
	\begin{displaymath}
	0=M_0\subsetneq M_1 \subsetneq \dotsb \subsetneq M_k=M
	\end{displaymath}
	is maximal, i.e. $M_{i-1}$ is a maximal $R$-submodule of $M_{i}$ for $i=1,\dotsc,k$, then 
	$k = L_R (M)$.
	
	

	\item If 
	\begin{displaymath}
	0\longrightarrow M'\longrightarrow M\longrightarrow M'' \longrightarrow 0
	\end{displaymath}
	is a short exact sequence of $R$-modules, then
	\begin{displaymath}
	L_R (M) = L_R (M')+L_R (M'').
	\end{displaymath}
\end{itemize}

Let us now return to the case that we were studying, with $\Oo$ a discrete valuation ring and $\m$ its maximal ideal. For each $i\in \N$, we have the maximal chain of $\Oo/\m^i$-modules
\begin{displaymath}
0\subsetneq \m^{i-1}/\m^i \subsetneq  \dotsb \subsetneq \m/\m^i\subsetneq\Oo/\m^i,
\end{displaymath}
and so $L_{\Oo/\m^i} (\Oo/\m^i)=i$. Furthermore, because $$\Oo/\m\cong (\Oo/\m^i)/(\m/\m^i)$$ as $(\Oo/\m^i)$-modules, we have that $L_{\Oo/\m^i} (\Oo/\m) = 1$. 

Now, for an element $a\in \Oo[F(S)]$ and $n\in \N$, not only can we consider the associated $\Oo$-module homomorphism $\phi_{X_n}^a\colon \Oo[X_n] \longrightarrow \Oo[X_n]$, but also the induced $(\Oo/\m^i)$-module homomorphism $\phi_{X_n,i}^a\colon (\Oo/\m^i)[X_n] \longrightarrow (\Oo/\m^i)[X_n]$ for each $i\in\N$.

We are now going to find a way to compute $\nu_{X_n}^a$ using the lengths of the kernels of these induced homomorphisms.

\begin{lemma}\label{lem:length kernel local}
	For each $n,i\in\N$, we have that 
	\begin{displaymath}
	\frac{L_{\Oo/\m^{i}}(\ker\phi_{X_n,i}^a)}{\abs{X_n}}=\sum_{j\in\N\cup\{\infty\}} \nu_{X_n}^a\big(\{\m^j\}\big)\min\{j,i\},
	\end{displaymath}
	where $\m^{\infty}=0$.
\end{lemma}

\begin{proof}
	Assume that $\phi_{X_n}^a$ is associated to a diagonal matrix $D_n\in\Mat_{\abs{X_n}}(\OO)$ in Smith normal form as before, with $\pi^{k_1},\dotsc,\pi^{k_t}$ the non-zero elements in the diagonal of $D_n$. Then, given $x\in \OO[X_n]$ of the form
	\begin{displaymath}
	x=\alpha_1 x_1 +\dotsb +\alpha_t x_t+\alpha_{t+1}x_{t+1}+\dotsb+\alpha_{\abs{X_n}}x_{\abs{X_n}},
	\end{displaymath}
	we have that 
	\begin{displaymath}
	D_n x=\alpha_1\pi^{k_1} x_1 +\dotsb +\alpha_t\pi^{k_t} x_t.
	\end{displaymath}
	Now, for each$i\in\N$ the induced homomorphism $\phi_{X_n,i}^a$ is associated to the reduction of $D_n$ modulo $\m^i$, which we will denote by $D_{n,i}$. Then,
	\begin{displaymath}
	D_{n,i} x= \alpha_1\pi^{k_1} x_1 +\dotsb +\alpha_r\pi^{k_r} x_r,
	\end{displaymath} 
	where $$r=\max \{1\leq j\leq t\mid k_j<i\}.$$
	Therefore, $x\in \ker\phi_{X_n,i}^a$ if and only if
	$\alpha_j\in \m^{i-k_j}$
	for all $j=1,\dotsc,r$. Hence, 
	\begin{displaymath}
	\ker\phi_{X_n,i}^a\isom(\m^{i-k_1}/\m^i)\oplus\dotsb \oplus (\m^{i-k_r}/\m^i)\oplus \bigoplus_{j=r+1}^{\abs{X_n}}(\Oo/\m^i),
	\end{displaymath}
	and so 
	\begin{displaymath}
	L_{\Oo/\m^i}(\ker\phi_{X_n,i}^a)=\abs{X_n}\sum_{j\in\N\cup\{\infty\}} \nu_{X_n}^a\big(\{\m^j\}\big)\min\{j,i\}.
	\end{displaymath}
\end{proof}

As a direct consequence of this result, we obtain the following one.
\begin{lemma}\label{lem:Measure local intervals}
	For each $n\in\N$, we have that
	\begin{displaymath}
	\nu_{X_n}^a\big([0,\m]\big) = \frac{L_{\Oo/\m}(\ker\phi_{X_n,1}^a)}{\abs{X_n}},
	\end{displaymath}
	and
	\begin{displaymath}
	\nu_{X_n}^a\big([0,\m^i]\big) = \frac{L_{\Oo/\m^{i}}(\ker\phi_{X_n,i}^a)}{\abs{X_n}} - \frac{L_{\Oo/\m^{i-1}}(\ker\phi_{X_n,i-1}^a)}{\abs{X_n}}
	\end{displaymath}
	for $i\geq 2$.
\end{lemma}

Then, we can compute 
\begin{displaymath}
\nu_{X_n}^a\big(\{\m^i\}\big) 
= \nu_{X_n}^a\big([0,\m^i]\big)- \nu_{X_n}^a\big([0,\m^{i+1}]\big) 
\end{displaymath}
for each $n,i\in\N$. Consequently, if we prove the convergence of the measures of intervals, we will also obtain the convergence at each ideal. In order to do so, following Lemma~\ref{lem:Measure local intervals} we will show that 
\begin{displaymath}
\frac{L_{\Oo/\m^i}(\ker\phi_{X_n,i}^a)}{\abs{X_n}}
\end{displaymath}
converges for each $i\in\N$ independent of the sofic approximation $(X_n)_{n\in\N}$. The proof of this is very similar to that of Theorem~\ref{thm:preLuck holds Amenable}.

\begin{proposition}\label{prop: lengths ultralimit local}
	Let $(X_n)_{n\in\N}$ and $(Y_n)_{n\in\N}$ be two sofic approximations of $G$ and $\omega$ be a non-principal ultrafilter on $\N$. Then, 
	\begin{displaymath}
	\lim_{n\to \omega}\frac{L_{\Oo/\m^i}(\ker \phi_{X_n,i}^a)}{\abs{X_k}}= \lim_{n\to \omega}\frac{L_{\Oo/\m^i}(\ker \phi_{Y_n,i}^a)}{\abs{Y_n}}
	\end{displaymath}
	for every $i\in\N$.
\end{proposition}

\begin{proof}
	Using the same argument as in the proof of Theorem~\ref{thm:preLuck holds Amenable}, we may assume without loss of generality that $\abs{X_n}=\abs{Y_n}$ for each $n\in\N$. As a consequence of Theorem~\ref{thm:Sofic Embeddings Conjugate}, for each $n\in\N$ there is some bijection $\sigma_n\colon X_n\longrightarrow Y_n$ such that, if we denote by
	\begin{displaymath}
	X_n'=\big\{x\in X_n\mid (\sigma_n^{-1}\circ\phi_{Y_n}^a\circ\sigma_n)(x)=\phi_{X_n}^a(x)\big\}, \qquad Y_n'=\sigma_n(X_n'),
	\end{displaymath}
	then
	\begin{displaymath}
	\lim_{n\to \omega}\frac{\abs{X_n'}}{\abs{X_n}}=\lim_{n\to \omega}\frac{\abs{Y_n'}}{\abs{Y_n}}=1.
	\end{displaymath}
	Observe that, given $x\in (\Oo/\m^i)[X_n']$, we have that $x\in\ker \phi_{X_n,i}^a$ if and only if $\sigma_n(x)\in\ker\phi_{Y_n,i}^a$. 
	Consider then the restrictions of $\phi_{X_n,i}^a$ to $(\Oo/\m^i)[X_n']$ and of $\phi_{Y_n,i}^a$ to $(\Oo/\m^i)[Y_n']$, which we will denote by $\phi_{X_n',i}^a$ and $\phi_{Y_n',i}^a$, respectively. Thus, we have that
	\begin{displaymath}
	\ker \phi_{X_n',i}^a=\ker \phi_{X_n,i}^a \cap (\Oo/\m^i)[X_n'], 
	\qquad \ker \phi_{Y_n',i}^a =\ker \phi_{Y_n,i}^a \cap (\Oo/\m^i)[Y_n'],
	\end{displaymath}
	and
	\begin{displaymath}
	\ker \phi_{X_n',i}^a \cong \ker \phi_{Y_n',i}^a
	\end{displaymath}
	for every $i\in\N$. Furthermore, the Second Isomorphism Theorem implies that $$\ker \phi_{X_n,i}^a/\ker \phi_{X_n',i}^a \lesssim (\Oo/\m^i)[X_n]/(\Oo/\m^i)[X_n'].$$ 
	As a consequence, from the short exact sequence of $(\Oo/\m^i)$-modules
	\begin{displaymath}
	0\longrightarrow \ker\phi_{X_n',i}\longrightarrow\ker \phi_{X_n,i}\longrightarrow \ker\phi_{X_n,i}/\ker\phi_{X_n',i}\longrightarrow 0
	\end{displaymath}
	we obtain that 
	\begin{align*}
		L_{\Oo/\m^i}(\ker \phi_{X_n,i}^a)&=L_{\Oo/\m^i}(\ker \phi_{X_n',i}^a)+L_{\Oo/\m^i}(\ker \phi_{X_n,i}^a/\ker \phi_{X_n',i}^a) \\
		& \leq L_{\Oo/\m^i}(\ker \phi_{X_n',i}^a)+ L_{\Oo/\m^i}\big((\Oo/\m^i)[X_n]/(\Oo/\m^i)[X_n']\big) \\
		& = L_{\Oo/\m^i}(\ker \phi_{X_n',i}^a) + \abs{X_n}-\abs{X_n'}
	\end{align*}
	for all $n\in\N$. 
	Analogously, we obtain the inequality
	\begin{displaymath}
	L_{\Oo/\m^i}(\ker \phi_{Y_n,i}^a) \leq L_{\Oo/\m^i}(\ker \phi_{Y_n',i}^a) + \abs{Y_n}-\abs{Y_n'}
	\end{displaymath}
	for all $n\in\N$.
	Therefore, since $\ker \phi_{X_n',i}^a \cong \ker \phi_{Y_n',i}^a$, we obtain that 
	\begin{align*}
		\lim_{n\to \omega}\frac{L_{\Oo/\m^i}(\ker \phi_{X_n,i}^a)}{\abs{X_n}} & = \lim_{k\to \omega}\frac{L_{\Oo/\m^i}(\ker \phi_{X_n',i}^a)}{\abs{X_k}} \\
		& = \lim_{n\to \omega}\frac{L_{\Oo/\m^i}(\ker \phi_{Y_n',i}^a)}{\abs{Y_n}} \\
		& = \lim_{n\to \omega}\frac{L_{\Oo/\m^i}(\ker \phi_{Y_n,i}^a)}{\abs{Y_n}}
	\end{align*}
	for each $i\in\N$. 
	%
\end{proof}

As a consequence of this result along with Lemma~\ref{lem: Measure local zero} and Lemma~\ref{lem:Measure local intervals},  we obtain the pointwise convergence of our measures $\nu_{X_n}^a$.

\begin{corollary}
	The limit $$\lim_{n\to \infty}\nu_{X_n}^a\big(\{I\}\big)$$ exists and is independent of the approximation $(X_n)_{n\in\N}$ for every $I\in\I$.
\end{corollary}


We would now want to prove that this gives us strong convergence of the measures $(\nu_{X_n}^a)_{n\in\N}$ to a probability measure $\nu_G^a$ on $\I$, so that we have that
$$\nu_G^a(\Omega)=\lim_{n\to \infty}\nu_{X_n}^a(\Omega)$$ for any subset $\Omega\subseteq \I$. 

Nevertheless, in order to show this in general, we would need some sort of uniform bound on the measures, which we have not found. In the next section, we will work over number fields in order to develop a global version of this construction, and in that case we will be able to obtain a uniform bound that will allow us to prove the strong convergence of the measures constructed.

\section{Approximation of Adelic Measures}

We will now develop a global version of the construction from the previous section. We will work over number fields and, using the structure theory of finitely generated modules over Dedekind domains, we will develop an analogue of the measures constructed in the previous section.

Let $K$ be a number field with ring of integers $\Oo$. Then, $\OO$ is a Dedekind domain, and so every non-zero ideal can be written in a unique way as a product of maximal ideals. We will denote by $\I$ the space of ideals of $\OO$ and by $\Imax\subseteq\I$ the set of maximal ideals.

Because $\OO$ is a Dedekind domain, 
a finitely generated $\OO$-module $M$ can be written as a direct sum $$M\isom M_{\operatorname{tors}}\oplus M/M_{\operatorname{tors}},$$ where $M\tors$ is the torsion submodule of $M$ and $M/M\tors$ is torsion-free. Now, the torsion part is of the form 
\begin{displaymath}
M_{\operatorname{tors}}\isom \OO/I_1\oplus \dotsb\oplus \OO/I_t
\end{displaymath}
with $I_1,\dotsc,I_t\in \I$ non-trivial ideals. Furthermore, it is possible to find such a decomposition with $I_{i+1}\subseteq I_i$ for $i=1,\dotsc,t-1$, in which case the ideals $I_1,\dotsc,I_t$ are unique. On the other hand, the torsion-free part is of the form 
\begin{displaymath}
M/M_{\operatorname{tors}}\isom J_1\oplus\dotsb\oplus J_r
\end{displaymath}
with $J_1,\dotsc,J_r\in \I$ non-zero ideals. 
Hence, we have that
\begin{displaymath}
M\isom \OO/I_1\oplus \dotsb\oplus \OO/I_t \oplus J_1\oplus\dotsb\oplus J_r.
\end{displaymath}
For more information on the structure of finitely generated modules over Dedekind domains, see \cite[\S 1.3]{narkiewicz}.

Given a maximal ideal $\m\in \Imax$ and a non-zero ideal $I\in\I$ we can consider the $\m$-adic valuation of $I$, denoted by $v_{\m}(I)$, which is defined as the biggest integer $n\geq 0$ such that $I\subseteq \m^n$. We also define $v_{\m}(0)=\infty$. Setting $v_{\m}(\alpha)=v_{\m}(\alpha\OO)$ for $\alpha\in\OO$, this defines a discrete valuation on $\OO$ that is extended naturally to $K$. Now, we may consider
\begin{displaymath}
\OO_{\m}=\{\alpha\in K\mid v_{\m}(\alpha)\geq 0\},
\end{displaymath}
the localisation of $\OO$ at $\m$, which is a discrete valuation ring with unique maximal ideal $\m\OO_{\m}$.

Let $G$ be a finitely generated amenable group, with $S\subseteq G$ a finite generating set, $G=F(S)/N$ with $N\normaleq F(S)$ and $(X_n)_{n\in\N}$ a sofic approximation of $G$. Take an element $a\in \Oo[F(S)]$. For each $n\in\N$, we have the induced linear map of $K$-vector spaces $\phi_{X_n}^a\colon K[X_n]\longrightarrow K[X_n]$. 

Consider now the $\Oo$-module $$M_n=\Oo[X_n]/\Oo[X_n]a.$$ Then, $M_n$ can be written in the form
\begin{displaymath}
M_n\cong \Oo/I_1\oplus \dotsb\oplus \Oo/I_t\oplus J_1\oplus\dotsb\oplus J_r,
\end{displaymath}
where $I_i,J_j\in \I$ are non-zero ideals such that $I_{i+1}\subseteq I_i\subsetneq\OO$ for any $i=1,\dotsc,t-1$ and $j=1,\dotsc,r$. Furthermore, the ideals $I_1,\dotsc,I_t$ are unique. 
%
%
%
%
%
%
%
%
%
We can then write 
\begin{displaymath}
M_n\cong (\Oo/\Oo)^s\oplus\Oo/I_1\oplus \dotsb\oplus \Oo/I_t\oplus J_1\oplus\dotsb\oplus J_r
\end{displaymath}
with $s\in\N$ such that 
\begin{displaymath}
s+t+r=\abs{X_n}.
\end{displaymath}

Using this decomposition of the $\OO$-module $M_n$, we can define for each $n\in\N$ a probability measure $\nu_{X_n}^a$ on $\I$ by 
\begin{displaymath}
\nu_{X_n}^a= \frac{1}{\abs{X_n}} \bigg(\sum_{i=1}^s \delta_{\Oo} + \sum_{i=1}^t \delta_{I_i} + \sum_{i=1}^{r}\delta_{0}\bigg).
\end{displaymath}

Observe that 
\begin{displaymath}
K[X_n]/K[X_n]a\isom K\otimes_{\OO}M_n\isom K^r,
\end{displaymath}
meaning that
\begin{displaymath}
\nu_{X_n}^a\big(\{0\}\big) =\frac{\dim_K\ker\phi_{X_n}^a}{\abs{X_n}}.
\end{displaymath}
Therefore, convergence at zero is once again equivalent to the Sofic L\"uck Approximation Conjecture, which holds for amenable groups by Theorem~\ref{thm:Luck holds Amenable}. Thus, we have the following result.

\begin{lemma}\label{lem: Measure adelic zero}
	The limit 
	\begin{displaymath}
	\lim_{n\to \infty}\nu_{X_n}^a\big(\{0\}\big)
	\end{displaymath}
	exists and is independent of the sofic approximation $(X_n)_{n\in \N}$.
\end{lemma}

Given a non-trivial ideal $I\in\I$, for each $n\in \N$ we can consider the induced $(\OO/I)$-module homomorphism $\phi_{X_n,I}^a\colon (\OO/I)[X_n]\longrightarrow(\OO/I)[X_n]$. 
Also, denote by
\begin{displaymath}
[0,I]=\{J\in \I\mid J\subseteq I\}.
\end{displaymath}
We are now going to find a way to compute the measures of intervals using the lengths of the kernels of these induced homomorphisms.

\begin{lemma}\label{lemma:length kernel global}
	For each non-trivial ideal $I\in \I$, we have that 
	\begin{displaymath}
	\frac{L_{\OO/I}(\ker\phi_{X_n,I}^a)}{\abs{X_n}}=\sum_{J\in \I}\nu_{X_n}^a\big(\{J\}\big)\bigg(\sum_{\m\in\Imax}\min\big\{v_{\m}(I),v_{\m}(J)\big\}\bigg).
	\end{displaymath}
\end{lemma}

\begin{proof}
	We have that 
	\begin{equation}\label{eqn:kernel global 1}
	\frac{L_{\OO/I}(\ker\phi_{X_n,I}^a)}{\abs{X_n}}=\sum_{\substack{\m\in\Imax\\  I\subseteq \m}}\frac{L_{\OO\m/I\OO_{\m}}(\ker\phi_{X_n,I\OO_{\m}}^a)}{\abs{X_n}}.
	\end{equation}
	Now, Lemma~\ref{lem:length kernel local} implies that 
	\begin{equation}\label{eqn:kernel global 2}
	\frac{L_{\OO\m/I\OO_{\m}}(\ker\phi_{X_n,I\OO_{\m}}^a)}{\abs{X_n}}= \sum_{j\in\N\cup\{\infty\}}\nu_{X_n,\OO_{\m}}^a\big(\{\m^j\OO_{\m}\}\big)\min\big\{v_{\m}(I),j\big\},
	\end{equation}
	where $\nu_{X_n,\OO_{\m}}^a$ denotes the local measure induced on $\mathcal{I}(\OO_{\m})$ for each maximal ideal $\m\in\Imax$. 
	Now, for each $j\in\N\cup\{\infty\}$ we have that
	\begin{equation}\label{eqn:kernel global 3}
	\nu_{X_n,\OO_{\m}}^a\big(\{\m^j\OO_{\m}\}\big)=\sum_{\substack{J\in\I\\  v_{\m}(J)=j}}\nu_{X_n}^a\big(\{J\}\big).
	\end{equation}
	Hence, from (\ref{eqn:kernel global 1}), (\ref{eqn:kernel global 2}) and (\ref{eqn:kernel global 3}) we obtain that 
	\begin{align*}
	\frac{L_{\OO/I}(\ker\phi_{X_n,I}^a)}{\abs{X_n}} & = \sum_{\m\in\Imax}\sum_{j\in\N\cup\{\infty\}}\sum_{\substack{J\in\I\\  v_{\m}(J)=j}} \nu_{X_n}^a\big(\{J\}\big)\min\big\{v_{\m}(I),v_{\m}(J)\big\} \\
	& = \sum_{J\in \I}\nu_{X_n}^a\big(\{J\}\big)\bigg(\sum_{\m\in\Imax}\min\big\{v_{\m}(I),v_{\m}(J)\big\}\bigg).
	\end{align*}
\end{proof}

Now, if we take a maximal ideal $\m\in\Imax$ and $d\in\N$, we have that
\begin{align*}
\frac{L_{\Oo/\m^{d}}(\ker\phi_{X_n,\m^d}^a)}{\abs{X_n}} - \frac{L_{\Oo/\m^{d-1}}(\ker\phi_{X_n,\m^{d-1}}^a)}{\abs{X_n}} & = \sum_{J\in[0,\m^d]}\nu_{X_n}^a(J) \\
& = \nu_{X_n}^a\big([0,\m^d]\big).
\end{align*}
Take now two maximal ideals $\m_1,\m_2\in\Imax$ and exponents $d_1,d_2\in\N$. Then, if we write $I=\m_1^{d_1}\m_2^{d_2}$ and $I'=\m_1^{d_1-1}\m_2^{d_2-1}$, we have that 
\begin{align*}
\frac{L_{\Oo/I}(\ker\phi_{X_n,I}^a)}{\abs{X_n}} - \frac{L_{\Oo/I'}(\ker\phi_{X_n,I'})}{\abs{X_n}} & = \\
& \!\!\!\!\!\!\!\!\!\!\!\!\!\!\!\!\!\!\!\!\!\!\!\!\!\!\!\!\!\!\!\!\!\!\!\!\!\!\!\!\!\!\!\!\!\!\!\!\!\!\!\!\!\!=  
\sum_{\substack{J\in[0,\m_1^{d_1}]\\ J\not\in[0,\m_2^{d_2}]}
} \nu_{X_n}^a(J)+\sum_{\substack{J\not\in[0,\m_1^{d_1}]\\ J\in[0,\m_2^{d_2}]}}\nu_{X_n}^a(J) +\sum_{J\in[0,\m_1^{d_1}\m_2^{d_2}]}2\nu_{X_n}^a(J) \\
& \!\!\!\!\!\!\!\!\!\!\!\!\!\!\!\!\!\!\!\!\!\!\!\!\!\!\!\!\!\!\!\!\!\!\!\!\!\!\!\!\!\!\!\!\!\!\!\!\!\!\!\!\!\!=  \nu_{X_k}^a\big([0,\m_1^{d_1}]\big)+\nu_{X_k}^a\big([0,\m_2^{d_2}]\big)+2\nu_{X_k}^a\big([0,\m_1^{d_1}\m_2^{d_2}]\big).
\end{align*}
We can thus show that, if $I=\m_1^{d_1}\dotsm\m_n^{d_n}$ and $I'=\m_1^{d_1-1}\dotsm\m_n^{d_n-1}$ with $\m_1,\dotsc,\m_n\in\Imax$ distinct maximal ideals and $d_1,\dotsc,d_n\in\N$, then 
\begin{align*}
\frac{L_{\Oo/I}(\ker\phi_{X_n,I}^a)}{\abs{X_n}} - \frac{L_{\Oo/I'}(\ker\phi_{X_n,I'})}{\abs{X_n}} & =\\
& \!\!\!\!\!\!\!\!\!\!\!\!\!=\sum_{j=1}^{n} \ \sum_{1\leq i_1<\dotsb <i_j\leq n}j
 \nu_{X_n}^a\big([0,\m_{i_1}^{d_{i_1}}\dotsm\m_{i_j}^{d_{i_j}}]\big),
\end{align*}
and so we can inductively write $\nu_{X_n}^a([0,I])$ in terms of the lengths of the kernels of $\phi_{X_n,J}^a$ for $J\in\I$ with $I\subseteq  J$.

We can now prove the following result, which is analogous to Proposition~\ref{prop: lengths ultralimit local}.

\begin{proposition}\label{prop:limit length ker global}
	Let $(X_n)_{n\in\N}$ and $(Y_n)_{n\in\N}$ be two sofic approximations of $G$ and $\omega$ be a non-principal ultrafilter on $\N$. Then, 
	\begin{displaymath}
	\lim_{n\to \omega}\frac{L_{\Oo/I}(\ker \phi_{X_n,I}^a)}{\abs{X_n}}= \lim_{n\to \omega}\frac{L_{\Oo/I}(\ker \phi_{Y_n,I}^a)}{\abs{Y_n}}
	\end{displaymath}
	for every non-trivial ideal $I\in \I$.
\end{proposition}
\begin{proof}
	Formula (\ref{eqn:kernel global 1}) tells us that we can write
	\begin{align*}
	\frac{L_{\OO/I}(\ker\phi_{X_n,I}^a)}{\abs{X_n}}=\sum_{\substack{\m\in\Imax\\  I\subseteq \m}}\frac{L_{\OO\m/I\OO_{\m}}(\ker\phi_{X_n,I\OO_{\m}}^a)}{\abs{X_n}}, \\
	\frac{L_{\OO/I}(\ker\phi_{Y_n,I}^a)}{\abs{Y_nk}}=\sum_{\substack{\m\in\Imax\\  I\subseteq \m}}\frac{L_{\OO\m/I\OO_{\m}}(\ker\phi_{Y_n,I\OO_{\m}}^a)}{\abs{Y_n}}.
	\end{align*}
	We can then apply Proposition~\ref{prop: lengths ultralimit local} to obtain that
	\begin{displaymath}
	\lim_{n\to \omega}\frac{L_{\OO\m/I\OO_{\m}}(\ker \phi_{X_n,I\OO_{\m}}^a)}{\abs{X_n}}= \lim_{n\to \omega}\frac{L_{\OO\m/I\OO_{\m}}(\ker \phi_{Y_n,I\OO_{\m}}^a)}{\abs{Y_n}}
	\end{displaymath}
	for every $\m\in\Imax$ with $I\subseteq \m$. Because there are only finitely many of these summands, this leads us to conclude that 
	\begin{displaymath}
	\lim_{n\to \omega}\frac{L_{\Oo/I}(\ker \phi_{X_n,I}^a)}{\abs{X_n}}= \lim_{n\to \omega}\frac{L_{\Oo/I}(\ker \phi_{Y_n,I}^a)}{\abs{Y_n}}.
	\end{displaymath}
\end{proof}

As a consequence of Lemma~\ref{lemma:length kernel global} and Proposition~\ref{prop:limit length ker global}, we obtain 
the convergence of the measures of intervals.

\begin{proposition}\label{prop:convergence nu intervals global}
	Given an ideal $I\in\I$, the limit
	\begin{displaymath}
	\lim_{n\to\infty}\nu_{X_n}^a\big([0,I]\big)
	\end{displaymath}
	exists and is independent of the approximation $(X_n)_{n\in\N}$.
\end{proposition}

Now that we have obtained convergence of our measures for intervals, we will work towards proving point-wise convergence. 
In order to do this, we will seek to write the measure of an ideal in terms of the measures of a finite number of intervals.

Observe that, given a non-zero ideal $I\in \I$, we can write
\begin{displaymath}
[0,I]=\{I\}\cup\bigg(\bigcup_{\m\in\I}[0,\m I]\bigg).
\end{displaymath}
From this, we can obtain that 
\begin{displaymath}
\nu_{X_n}^a\big(\{I\}\big)=\nu_{X_n}^a\big([0,I]\big)-\nu_{X_n}^a\bigg(\bigcup_{\m\in\I}[0,\m I]\bigg)
\end{displaymath}
for each $n\in\N$. 
Our goal will now be to approximate the measures of this union by the measures of a finite union of intervals. We will do this by studying the sizes of the ideals that can appear in the decompositions of the modules 
$M_n=\OO[X_n]/\OO[X_n]a$.

If $I\in\I$ is a non-zero ideal, we can consider its norm, defined as
\begin{displaymath}
N(I)=\abs{\OO/I}.
\end{displaymath}
Then, for $\alpha\in \OO$ we have that $$N(\alpha\OO)=\abs[\big]{N_{K:\Q}(\alpha)},$$
where
\begin{displaymath}
N_{K:\Q}(\alpha)=\bigg(\prod_{i=1}^{k}\alpha_i\bigg)^{\abs{K:\Q}/k}
\end{displaymath}
with $\alpha_1,\dotsc,\alpha_k\in\bar{\Z}$  the roots of the minimal polynomial of $\alpha$ over $\Q$. Furthermore, it is well-known that for any constant $\lambda>0$ there are only finitely many non-trivial ideals $I\in\I$ such that $N(I)<\lambda$.

Now, given any matrix $A\in\Mat_k(\OO)$, we can consider the finitely generated $\OO$-module $$M_A=\OO^k/\Oo^k A.$$
Then, we can define
\begin{displaymath}
\detp(A)=\abs[\big]{(M_A)\tors}.
\end{displaymath}
In particular,  identifying the $\OO$-module homomorphism $\phi_{X_n}^a$ with its associated matrix, we can write 
\begin{displaymath}
\detp(\phi_{X_n}^a)=\abs[\big]{(M_n)\tors}.
\end{displaymath}


\begin{lemma}\label{lem:det+(M) minors}
	Let $A\in\Mat_{k}(\OO)$ with $\dim_K \im A=t$ and $I\in \I$ the ideal generated by all the non-zero $t\times t$ minors of $A$. Then, we have that $$\detp(A)=N(I).$$
\end{lemma}
\begin{proof}
	First, observe that if we write 
	\begin{displaymath}
	(M_A)\tors\isom \OO/I_1\oplus\dotsb\oplus \OO/I_t
	\end{displaymath}
	with $I_1,\dotsc I_t\in\I$ non-trivial ideals, 
	then
	\begin{displaymath}
	\OO_{\m}^k/\OO_{\m}^k A\isom \OO_{\m}\otimes_{\OO}M_A,
	\end{displaymath}
	and so
	\begin{align*}
	(\OO_{\m}^k/\OO_{\m}^k A)\tors & \isom \OO_{\m}\otimes_{\OO}(M_A)\tors \\
	& \isom \Oo_{\m}/I_1\Oo_{\m}\oplus\dotsb\oplus \OO_{\m}/I_t\Oo_{\m}
	\end{align*}
	for any maximal ideal $m\in\Imax$. But in the local case, we have that 
	\begin{displaymath}
	\abs[\big]{(\OO_{\m}^k/\OO_{\m}^k A)\tors}=\abs{\OO_{\m}/I\OO_{\m}}
	\end{displaymath}
	due to the existence of the Smith normal form and Remark~\ref{rmk:Smith normal form module}, and so we can compute 
	\begin{align*}
	\detp(A) & = \abs[\big]{(M_A)\tors} \\
	& = \abs{\Oo/I_1\oplus\dotsb\oplus \OO/I_t} \\
	& = \prod_{\m\in \Imax}\abs{\Oo_{\m}/I_1\Oo_{\m}\oplus\dotsb\oplus \OO_{\m}/I_t\Oo_{\m}} \\
	& = \prod_{\m\in \Imax} \abs[\big]{\OO_{\m}\otimes_{\OO}(M_A)\tors} \\
	& = \prod_{\m\in \Imax} \abs{\OO_{\m}/I\OO_{\m}} \\
	& = \abs{\OO/I} \\
	& = N(I).
	\end{align*}
\end{proof}

Given an element $\alpha\in\OO$, let $\alpha_1,\dotsc,\alpha_k\in\bar\Z$ be the roots of the minimal polynomial of $\alpha$ over $\Q$. Then, define  
\begin{displaymath}
\ceil{\alpha}=\max_{i=1,\dotsc,k}\abs{\alpha_i}.
\end{displaymath}

\begin{remark}\label{rmk:ceiling sum product ineq}
	Given $\alpha,\beta\in \OO$, we can check that 
	\begin{displaymath}
	\ceil{\alpha+\beta}\leq \ceil{\alpha}+\ceil{\beta}
	\end{displaymath}
	and
	\begin{displaymath}
	\ceil{\alpha\beta}\leq \ceil{\alpha}\ceil{\beta}.
	\end{displaymath}
\end{remark}

More generally, given a non-zero matrix $A=(a_{ij})\in\Mat_k(\OO)$, we can define
\begin{displaymath}
\ceil{A}=\max_{j=1,\dotsc,k}\sum_{i=1}^k\ceil{a_{ij}},
\end{displaymath}
and set $\ceil{0}=1$.

\begin{lemma}
	Given a non-zero $\alpha\in\OO$, we have that  
	\begin{displaymath}
	N(\alpha\OO)\leq \ceil{\alpha}^{\abs{K:\Q}}.
	\end{displaymath}
\end{lemma}
\begin{proof}
	Let $\alpha_1,\dotsc,\alpha_k\in\bar\Z$ be the roots of the minimal polynomial of $\alpha$ over $\Q$. Then, we have that 
	\begin{align*}
	\detp(\alpha) & =\abs[\big]{N_{K:\Q}(\alpha)} \\
	& =\bigg(\prod_{i=1}^{k}\abs{\alpha_i}\bigg)^{\abs{K:\Q}/k} \\
	& \leq \ceil{\alpha}^{\abs{K:\Q}}.
	\end{align*}
\end{proof}

This bound can now be generalised to matrices.

\begin{lemma}\label{lem:det+(A)<= ceiling A}
	Given $A\in\Mat_{k}(\OO)$, we have that  
	\begin{displaymath}
	\detp(A)\leq \ceil{A}^{k\abs{K:\Q}}.
	\end{displaymath}
\end{lemma}
\begin{proof}
	Let $t=\dim_K\im A$ and $I\in \I$ be the ideal generated by all the non-zero $t\times t$ minors of $A$. By Lemma~\ref{lem:det+(M) minors}, we have that $\detp(A)=N(I)$. In particular, if $\beta\in\OO$ is a non-zero $t\times t$ minor of $A$, we have that 
	\begin{displaymath}
	\detp(A)=\abs{\OO/I}\leq\abs{\OO/\beta\OO} =N(\beta\OO).
	\end{displaymath}
	Then, by Remark~\ref{rmk:ceiling sum product ineq} we have that 
	\begin{displaymath}
	\ceil{\beta}\leq\ceil{A}^t\leq \ceil{A}^k,
	\end{displaymath}
	and so 
	\begin{displaymath}
	\detp(A)\leq N(\beta\OO)\leq\ceil{\beta}^{\abs{K:\Q}}\leq \ceil{A}^{k\abs{K:\Q}}.
	\end{displaymath}
\end{proof}

Now, given $a\in\OO[F(S)]$ of the form $$a=\sum_{w\in F(S)}a_w w,$$
we set 
\begin{displaymath}
\ceil{a}=\sum_{w\in F(S)}\ceil{a_w}.
\end{displaymath}

\begin{lemma}\label{lem:ceil phi^a<=ceil a}
	Given $a\in\OO[F(S)]$, we have that
	\begin{displaymath}
	\ceil{\phi_{X_n}^a}\leq \ceil{a}
	\end{displaymath}
	for all $n\in\N$.
\end{lemma}

\begin{proof}
	Assume that
	\begin{displaymath}
	X_n=\{x_1,\dotsc,x_k\}
	\end{displaymath}
	and $\phi_{X_n}^a$ is associated to the matrix $A_n=(a_{ij})$. 
	Then, for each $i=1,\dotsc,k$ we have that 
	\begin{equation}\label{eqn:ceil1}
	\phi_{X_n}^a(x_i)= a_{i1}x_1+\dotsb +a_{in}x_k.
	\end{equation}
	On the other hand, if we write 
	$$a=\sum_{w\in F(S)}a_ww,$$
	then for each $i=1,\dotsc,k$ we have that 
	\begin{equation}\label{eqn:ceil2}
	\phi_{X_n}^a(x_i)=x_i a=\sum_{\substack{w\in F(S)\\ x_iw=x_1}}a_w x_1 +\dotsb +\sum_{\substack{w\in F(S)\\ x_iw=x_k}}a_w x_k.
	\end{equation}
	Therefore, combining (\ref{eqn:ceil1}) and (\ref{eqn:ceil2}) and applying Remark~\ref{rmk:ceiling sum product ineq}, we obtain that 
	\begin{align*}
	\sum_{i=1}^k\ceil{a_{ij}} & = \sum_{i=1}^k\ceil[\bigg]{\sum_{\substack{w\in F(S)\\ x_iw=x_j}}a_w} \\
	& \leq \sum_{i=1}^k\sum_{\substack{w\in F(S)\\ x_iw=x_j}}\ceil{a_w} \\
	& \leq \sum_{w\in F(S)}\ceil{a_w} \\
	& = \ceil{a}
	\end{align*}
	for all $j=1,\dotsc,k$. 
	As a consequence,
	\begin{displaymath}
	\ceil{\phi_{X_n}^a}=\max_{j=1,\dotsc,k}\sum_{i=1}^k\ceil{a_{ij}}\leq \ceil{a}
	\end{displaymath}
	for all $n\in\N$.
\end{proof}

This result allows us to give a uniform bound for $\detp(\phi_{X_n}^a)$.

\begin{corollary}\label{cor:bound det+ phi ceil a}
	Given $a\in\OO[F(S)]$, we have that 
	\begin{displaymath}
	\detp(\phi_{X_n}^a)\leq \ceil{a}^{\abs{X_n}\abs{K:\Q}}
	\end{displaymath}
	for all $n\in\N$.
\end{corollary}
\begin{proof}
	Applying both Lemma~\ref{lem:det+(A)<= ceiling A} and Lemma~\ref{lem:ceil phi^a<=ceil a}, we obtain that
	\begin{displaymath}
	\detp(\phi_{X_n}^a)\leq \ceil{\phi_{X_n}^a}^{\abs{X_n}\abs{K:\Q}}\leq \ceil{a}^{\abs{X_n}\abs{K:\Q}}
	\end{displaymath}
	for any $n\in\N$.
\end{proof}

Consequently, we can bound 
\begin{displaymath}
\detp(\phi_{X_n}^a)\leq c^{\abs{X_n}}
\end{displaymath}
with some $c>0$ for every $n\in\N$. This allows us to bound the measures of sets of big ideals.

\begin{lemma}\label{lem:bound nu Omega norm}
	Given $\Omega\subseteq \I\setminus\{\OO,0\}$ and $c=\ceil{a}^{\abs{K:\Q}}$, we have that
	\begin{displaymath}
	\nu_{X_n}^a(\Omega)\leq \frac{1}{\log_c\min_{I\in\Omega} \big\{N(I)\big\}}
	\end{displaymath}
	for every $n\in\N$.
\end{lemma}
\begin{proof}
	If $k=\abs{X_n}\nu_{X_n}^a(\Omega)$, then we must have 
	\begin{displaymath}
	(M_{n})\tors\isom (\OO/I_1\oplus\dotsb\oplus\OO/I_k)\oplus (\Oo/I_{k+1}\oplus\dotsb\oplus \Oo/I_t)
	\end{displaymath}
	with $I_1,\dotsc,I_k\in\Omega$ and $I_{k+1},\dotsc,I_t\not\in\Omega$, and so 
	\begin{align*}
	\detp(\phi_{X_n}^a) & = \abs[\big]{(M_{n})\tors} \\
	& \geq N(I_1)\dotsm N(I_k) \\
	& \geq \min_{I\in\Omega}\big\{N(I)^k\big\}.
	\end{align*}
	Applying now Corollary~\ref{cor:bound det+ phi ceil a}, we obtain that
	\begin{displaymath}
	\min_{I\in\Omega}\big\{N(I)^k\big\}\leq c^{\abs{X_n}},
	\end{displaymath}
	which taking logarithms gives us that
	\begin{displaymath}
	\nu_{X_n}^a(\Omega)\leq \frac{1}{\log_c\min_{I\in\Omega} \big\{N(I)\big\}}
	\end{displaymath}
	for all $n\in\N$.
\end{proof}

Now, given $\lambda>0$ we will denote by
\begin{displaymath}
\I_{\lambda}=\big\{J\in\I\setminus\{\OO,0\}\mid N(J)> \lambda\big\}.
\end{displaymath}
Using the previous result, we can show that the measures of $\I_{\lambda}$ are small for large $\lambda$.

\begin{proposition}\label{prop:nu(I lambda) < epsilon}
	Given $\varepsilon>0$, there exists some $\lambda>0$ such that
	$$\nu_{X_n}^a\big(\I_{\lambda}\big)<\varepsilon$$ 
	for all $n\in\N$.
\end{proposition}
\begin{proof}
	If we take $\lambda>c^{1/\varepsilon}$, then Lemma~\ref{lem:bound nu Omega norm} implies that
	\begin{displaymath}
	\nu_{X_n}^a\big(\I_{\lambda}\big)\leq \frac{1}{\log_c\lambda}<\varepsilon
	\end{displaymath}
	for all $n\in\N$.
\end{proof}

\begin{remark}
The previous result can be summed up by saying that for big $\lambda>0$ the value $\nu_{X_n}^a(\I_{\lambda})$ is uniformly small. It can be seen as an analogue to Lemma~\ref{eqn:Luck bound}, which said that for small $\varepsilon>0$ the value $\mu_{X_n}^a((0,\varepsilon))$ was uniformly small.
\end{remark}

Let $I\in \I$ be a non-zero ideal. Then, we have that 
\begin{displaymath}
[0,I]=\{I\}\cup\bigg(\bigcup_{\m\in\Imax}[0,\m I]\bigg),
\end{displaymath}
Given $\varepsilon>0$, as a consequence of Proposition~\ref{prop:nu(I lambda) < epsilon} there is some $\lambda>0$ such that $\nu_{X_n}^a(\I_{\lambda})<\varepsilon$ for all $n\in\N$. Furthermore, there are only finitely many  distinct maximal ideals $\m_1,\dotsc,\m_k\in\Imax$ with $N(\m_i I)\leq \lambda$ for $i=1,\dotsc,k$. Then, we have that
\begin{displaymath}
\bigcup_{\m\in\Imax}[0,\m I]=\bigg(\bigcup_{i=1}^k[0,\m_i I]\bigg) \cup \big(\I_{\lambda}\cap[0,I]\big).
\end{displaymath}
Using now that $$[0,\m_iI]\cap[0,\m_jI]=[0,\m_i\m_jI]$$ for any $i\not=j$, we can apply the inclusion–exclusion principle to compute
\begin{displaymath}
\nu_{X_n}^a\bigg(\bigcup_{i=1}^k[0,\m_i I]\bigg)=  
\sum_{j=1}^{k}(-1)^{j-1} \sum_{1\leq i_1<\dotsb<i_j\leq k}\nu_{X_n}^a\big([0,\m_{i_1}\dotsm\m_{i_j} I]\big),
\end{displaymath}
which converges independent of the approximation $(X_n)_{n\in\N}$ because the measures of the intervals converge by Proposition~\ref{prop:convergence nu intervals global}.  
Furthermore, if we write $$\I_{\lambda}'=[0,I]\setminus\Bigg(\{I\}\cup\bigg(\bigcup_{i=1}^k[0,\m_i I]\bigg)\Bigg)\subseteq \I_{\lambda}\cap[0,I],$$
then
\begin{displaymath}
\nu_{X_n}^a\big(\I_{\lambda}'\big)\leq \nu_{X_n}^a\big(\I_{\lambda}\big)<\varepsilon.
\end{displaymath}
Therefore, 
\begin{displaymath}
\nu_{X_n}^a\bigg(\bigcup_{i=1}^k[0,\m_i I]\bigg)\leq \nu_{X_n}^a\bigg(\bigcup_{\m\in\Imax}[0,\m I]\bigg) <\nu_{X_n}^a\bigg(\bigcup_{i=1}^k[0,\m_i I]\bigg)+\varepsilon
\end{displaymath}
for every $n\in\N$, which implies that $$\nu_{X_n}^a\bigg(\bigcup_{\m\in\Imax}[\m I,0]\bigg)$$
converges independent of the approximation $(X_n)_{n\in\N}$. As a consequence, 
\begin{displaymath}
\nu_{X_n}^a\big(\{I\}\big)=\nu_{X_n}^a\big([I,0]\big)-\nu_{X_n}^a\bigg(\bigcup_{\m\in\Imax}[\m I,0]\bigg)
\end{displaymath}
converges as well for any non-zero $I\in\I$ independent of the approximation $(X_n)_{n\in\N}$.

More generally, the uniform bound in Proposition~\ref{prop:nu(I lambda) < epsilon} allows us to approximate uniformly the measures of each subset  $\Omega\subseteq \I$ by the measures of a finite number of ideals. This implies that for each $\Omega\subseteq \I$ the limit
\begin{displaymath}
	\lim_{n\to \infty}\nu_{X_n}^a(\Omega)
\end{displaymath}
exists and is independent of the approximation $(X_n)_{n\in \N}$. Indeed, because $\I$ is a countable space, we can write
\begin{displaymath}
\nu_{X_n}^a(\Omega)=\sum_{I\in\Omega}\nu_{X_n}^a\big(\{I\}\big)
\end{displaymath}
and, because of Proposition~\ref{prop:nu(I lambda) < epsilon}, this series converges uniformly. Therefore, taking limits commutes with the sum, and 
\begin{displaymath}
\lim_{n\to \infty}\nu_{X_n}^a(\Omega)= \sum_{I\in \Omega} \lim_{n\to \infty}\nu_{X_n}^a\big(\{I\}\big)
\end{displaymath}

As a consequence, we can consider the limit probability measure $\nu^a_G$ on $\I$, given by $$\nu^a_G(\Omega)=\lim_{n\to \infty}\nu_{X_n}^a(\Omega)$$ for any subset $\Omega\subseteq \I$. This can be summed up in the following result. 

\begin{theorem}
	Let $G$ be a finitely generated amenable group with $S\subseteq G$ a finite generating set and $(X_n)\ninN$ a sofic approximation. Let $K$ be a number field with ring of integers $\OO$, and $a\in \OO[F(S)]$.  Then, the sequence of probability measures $(\nu_{X_n}^a)_{n\in \N}$ converges strongly to some probability measure $\nu_G^a$ on $\I$, independent of the sofic approximation $(X_n)_{n\in \N}$.
\end{theorem}

\section{Adelic L\"uck Approximation}

Throughout this chapter, we have always worked with amenable groups, which allowed us to use the characterisation of amenability in Theorem~\ref{thm:Sofic Embeddings Conjugate} to prove the convergence of our measures. Nevertheless, the constructions of the measures themselves are not dependent on whether our group is amenable or not.

We would then like to finish by conjecturing that the measures constructed in the last section converge in general, even for non-amenable groups. This conjecture, which we will call the Adelic L\"uck Approximation Conjecture, serves as a generalisation of the Sofic L\"uck Approximation Conjecture.

\begin{conjecture}
	Let $G$ be a finitely generated sofic group with $S\subseteq G$ a finite generating set, $(X_n)_{n\in\N}$ a sofic approximation of $G$, $K$  a number field with ring of integers $\OO$, and $a\in\OO[F(S)]$. Consider for each $n\in \N$ the measure $\nu_{X_n}^a$ as before. Then, the sequence of probability measures $(\nu_{X_n}^a)_{n\in \N}$ converges strongly to some probability measure $\nu_G^a$ on $\I$, independent of the sofic approximation $(X_n)_{n\in \N}$.
\end{conjecture}